\documentclass[preprint]{imsart}

\RequirePackage[OT1]{fontenc}
\RequirePackage{amsthm,amsmath,bm,amssymb}
\RequirePackage[numbers]{natbib}
\RequirePackage[colorlinks,citecolor=blue,urlcolor=blue]{hyperref}
\RequirePackage{graphicx} 
\RequirePackage{epsfig}
\RequirePackage{multirow}

\startlocaldefs
\numberwithin{equation}{section}
\theoremstyle{plain}

\newtheorem{theorem}{Theorem}[section]

\newtheorem{remark}{Remark}[section]

\endlocaldefs

\begin{document}

\begin{frontmatter}
\title{High-dimensional varying index coefficient quantile regression model}
\runtitle{QUANTILE VARYING INDEX COEFFICIENT MODEL}

\begin{aug}
\author{\fnms{JIALIANG} \snm{LI}\ead[label=e1]{stalj@nus.edu.sg}\thanksref{t1}}
\and
\author{\fnms{JING} \snm{LV}\corref{}\ead[label=e2]{lvjing@swu.edu.cn}\thanksref{t2}}
\thankstext{t2}{Corresponding author. }
\runauthor{J. LI AND J. LV}

\affiliation{National University of Singapore\thanksmark{t1} and Southwest University\thanksmark{t2} }

\address{JIALIANG LI\\
Department of Statistics and Applied Probability\\ National University of Singapore\\
Singapore, 119077\\
\printead{e1}\\}

\address{JING LV\\
School of Mathematics and Statistics\\ Southwest University\\Chongqing, 400715, China\\
\printead{e2}\\}

\end{aug}

\begin{abstract}
Statistical learning evolves quickly with more and more sophisticated models proposed to incorporate the complicated data structure from modern scientific and business problems. Varying index coefficient models extend varying coefficient models and single index models, becoming the latest state-of-the-art for semiparametric regression. This new class of models offers greater flexibility to characterize complicated nonlinear interaction effects in regression analysis. To safeguard against outliers and extreme observations, we consider a robust quantile regression approach to estimate the model parameters in this paper. High-dimensional loading parameters are allowed in our development under reasonable theoretical conditions. In addition, we propose a regularized estimation procedure to choose between linear and non-linear forms for interaction terms. We can simultaneously select significant non-zero loading parameters and identify linear functions in varying index coefficient models, in addition to estimate all the parametric and nonparametric components consistently. Under technical assumptions, we show that the proposed procedure is consistent in variable selection as well as in linear function identification, and the proposed parameter estimation enjoys the oracle property. Extensive simulation studies are carried out to assess the finite sample performance of the proposed method. We illustrate our methods with an environmental health data example.
\end{abstract}


\end{frontmatter}

\section{Introduction}

Semiparametric regression models are powerful statistical learning approaches and become more and more popular in scientific and business research studies since they can enjoy the merits of both parametric and nonparametric models. We consider the varying index coefficient model (VICM) recently proposed in the literature (\cite{MS15}). This new class of models extends varying coefficient models (\cite{FZ99}), single-index models (\cite{XTLZ02}), single index coefficient models (\cite{XP13}) and almost all other familiar semiparametric models, thus becoming the latest state of the art. To safeguard against outliers and extreme observations, we consider a robust quantile regression approach to fit the VICM in this paper. Specifically, for a given quantile level $\tau \in (0,1)$, varying index coefficient quantile regression models are given by
\begin{equation}\label{eq1}
Q_{\tau}\left( Y | \bm X, \bm Z \right)= \sum\limits_{l = 1}^d {{m_{\tau,l}}} ({\bm Z^T}{\bm \beta _{\tau,l}}){X_{l}} ,
\end{equation}
where $\bm X=\left(X_1,...,X_d\right)^T$ with $X_1\equiv 1$ and $\bm Z=\left(Z_1,...,Z_p\right)^T$  are covariates for the response variable $Y \in {R}$, ${\bm \beta _{\tau,l}} = {({\beta _{\tau,l1}},...,{\beta _{\tau,l{p}}})^T}$ are unknown loading parameters for the $l$th covariate $X_l$ and $m_{\tau,l}(\cdot)$ are unknown nonparametric functions, $l=1,...,d$. Let $\varepsilon_{\tau}=Y-Q_\tau(Y|{\bm X},{\bm Z})$ be the model error with an unspecified conditional density function $f_{\tau}\left(\cdot|\bm X, \bm Z \right)$ and conditional cumulative distribution function $F_{\tau}(\cdot|\bm X,\bm Z)$ given $\left( \bm X, \bm Z\right)$. Please note that $\varepsilon_{\tau}$'s conditional $\tau$th quantile equals zero, that is, $P\left(\varepsilon_{\tau}<0|\bm X, \bm Z\right)=\tau$. In the rest of the article, we drop the subscript $\tau$ from ${\bm \beta _{\tau,l}}$, $m_{\tau,l}(\cdot)$, $\varepsilon _{\tau}$, $f_{\tau}\left(\cdot|\bm X, \bm Z \right)$ and $F_{\tau}(\cdot|\bm X,\bm Z)$ to simplify the notations, but it is helpful to bear in mind that all those quantities are $\tau$-specific. For the sake of identifiability, we assume that $\bm \beta  = {\left( {\bm{\beta } _{1}^T,...,\bm{\beta} _{d}^T} \right)^T}$ belongs to the following parameter space:
\[\Theta  = \left\{ {\bm \beta  = {{\left( {\bm \beta _{l}^T:1 \le l \le d} \right)}^T}:\left\| {{\bm \beta _{l}}} \right\|_2 = 1,{\beta _{l1}} > 0,{\bm \beta _{l}} \in {{\mathcal{R}}^{p}}} \right\},\]
where $\left\|.\right\|_2$ denotes the $L_2$ norm such that ${\left\| \bm\xi  \right\|_2} = {\left( {\xi _1^2 + ... + \xi _s^2} \right)^{{1 \mathord{\left/
 {\vphantom {1 2}} \right.
 \kern-\nulldelimiterspace} 2}}}$ for any vector $\bm \xi  = {\left( {{\xi _1},...,{\xi _s}} \right)^T} \in {{\mathcal{R}}^s}$.
 Model (\ref{eq1}) is quite general and includes many other existing models as special cases. For example, (i) when $m_l(\cdot)$ are assumed to be constant or linear function, it reduces to the linear regression model with interactions; (ii) when $d=1$ and $X_l=1$, it is the single index model; (iii) when $m_l(\cdot)$ are set as constant for $l\geq 2$ and $X_1=1$, it is the partial linear single-index model; (iv) when common coefficient vector $\bm{\beta}_l$ are used, it is the single index coefficient model; (v) when $p=1$ and by the definition of $\Theta$ we have $\beta_l=1$, it reduces to the varying coefficient model. The VICM is very flexible to model and assess nonlinear interaction effects between the covariate $\bm X$ and $\bm Z$. Our main interest is to make statistical inference on both the loading coefficients $\bm \beta_{l}$ and the nonparametric functions $m_l \left( \cdot\right)$.

  \cite{MS15} proposed a profile least squares estimation procedure for the VICM and established its theoretical properties. Their work focused on mean regression, which is most suitable for nicely distributed data such as Gaussian and may perform badly in the presence of outliers and heavy-tailed errors. Our model (\ref{eq1}) imposes different assumptions on the error structure and thus produces a novel and robust framework applicable for wider applications. The estimation methods and the associated asymptotic theories are thus totally different from \cite{MS15}.

  Since the seminal work of \cite{KB78}, quantile regression has emerged as an important alternative to mean regression. It is well understood that inference based on quantile regression is more robust against distribution contamination (\cite{K05}). A full range of quantile analysis can provide more complete description of the conditional distribution. It is now widely acknowledged that quantile regression based analysis may lead to more appropriate findings. For example, climatologists often pay close attention to how the high quantiles of tropical cyclone intensity change over time (\cite{EKJ2008}), as it not only generates strong winds and waves, but also often results in heavy rain and storm surges, causing serious disasters. In another health sciences example, medical scientists often study the influences of maternal behaviors on the low quantiles of the birth weight distributions (\cite{A2001}). Consider one more case study from business and economics: petroleum is a primary source of non-renewable energy and has important influence on industrial production, electric power generation and transportation (\cite{MRT2009}). Thus, most analysts particularly focus on the high quantiles of oil prices, as oil price fluctuations have considerable impacts on economic activity. The quantile regression framework considered in this paper may impact all these important fields where direct application of mean regression is inappropriate.

  In recent decades the classical parametric quantile regression has been integrated with semiparametric models to produce more flexible inference tools. Here we only list a few relevant works among the abundant developments. For single index models, \cite{WYY10} developed a robust minimum average variance estimation procedure based on the familiar quantile regression. \cite{KX12} combined quantile regression and a penalty function to develop an adaptive quantile estimation algorithm. \cite{MH16} considered a pseudo-profile likelihood approach, which enables a straight forward statistical inference on the index coefficients. \cite{CA16} proposed a non-iterative quantile estimation algorithm for heteroscedastic data, and provided the asymptotic properties of the proposed approach. For varying coefficient model, \cite{TWZ13} developed a new variable selection procedure by utilizing basis function approximation and a class of group versions of the adaptive LASSO penalty. \cite{PXK14} considered a shrinkage estimator under quantile regression. For single index coefficient models, \cite{JQ16} considered a new estimation procedure to reduce the computing cost of existing back-fitting algorithm. \cite{ZLL18} developed a bias-corrected quantile estimating equations
and presented the fixed-point algorithm for fast and accurate computation of the parameter estimates. Other related works about semiparametric quantile regression include \cite{JPCL12,JZQC13,LP15,SPMM16, WSZ12,WZL13,WZ11,ZL17,ZHL12}, among many others.

Another important contribution of this paper is that we consider the high-dimensional learning issues for the new VICM. In fact, recent advances in technologies for cheaper and faster data acquisition and storage have led to an explosive growth of data complexity in a variety of scientific areas such as medicine, economics and environmental science. We have to consider a realistic solution facing the ``large $n$, diverging $p$'' data setting. Specifically we will allow the dimension of the covariates $\bm Z$ to increase to infinity as the sample size increases. Many penalty-based estimation methods are proposed in modern statistical community to address the high dimensional issue (\cite{ BV11,FL06,G15,HTW15}). This framework can effectively reduce the model bias and improve the prediction performance of the fitted model. \cite{FP04} first studied nonconcave penalized likelihood estimation when the number of covariates increases with the sample size. \cite{WZQ12} extended the method to generalized linear models for longitudinal measurements. High dimensional issue has also been investigated for semiparametric models. \cite{LY12}combined basis function approximation with the SCAD penalty to propose a variable selection procedure for generalized varying coefficient partially linear models with diverging number of parameters, and also established the consistency and oracle property of their method. \cite{WW15} applied the SCAD penalty to perform variable selection for single index prediction models with a diverging number of index parameters. \cite{FLL17} presented a penalized empirical likelihood approach for high dimensional semiparametric models.

 Variable selection for model  (\ref{eq1}) is challenging since the high-dimensional loading parameter is structured within the unknown nonparametric function coefficients. We adopt a spline basis approximation to the estimation of $m_l(\cdot)$ and consequently estimate the unknown loading parameters vector $\bm{\beta}_l$ under the sparsity assumption. In addition, we tackle the problem of correctly identifying the linear interaction effects between covariates. That is, we want to decide whether it is necessary to model $m_l(\cdot)$ nonparametrically for all the $d$ varying index functions. \cite{MS15} constructed a generalized likelihood ratio statistic to test whether there exists a linear interaction effect between covariates. Although this test approach works very well for low-dimensional problems, it is computationally infeasible when the number of covariate is large. To this end, we develop a group penalization method that can quickly and effectively differentiate linear functions from nonparametric functions. The theoretical justification is also non-trivial for this complicated setting.

The rest of the paper is organized as follows. In Sect. \ref{sect.2}, using the B-spline basis approximation, we construct robust quantile estimating equations for loading parameters and obtain the estimators of unknown nonparametric functions by minimizing the quantile loss function. Asymptotic properties of the proposed estimators are also established in this section. In Sect. \ref{sect.3}, we consider high-dimensional issues and describe the variable selection procedure for loading parameters. Theoretical results are also presented including the estimation convergence rate, selection consistency and oracle property of estimators. In Sect. \ref{sect.4}, a group penalized method is proposed to identify linear functional effects along with the theoretical properties. In Sect. \ref{sect.5}, simulation studies and real data analysis are provided to illustrate our methods. In Sect. \ref{sect.6}, we conclude with some remarks. All technical proofs are given in the Appendix B.

\section{Quantile Regression Estimation of Functions and Loadings in VICM}\label{sect.2}
\subsection{Estimation procedures}\label{subsect.2.1}
Suppose that $\{ (\bm{X}_i,\bm{Z}_i,Y_i),1 \le i \le n\} $ is an independent and identically distributed sample from model (1). Without loss of generality, we assume that $ {\bm Z_i^T}{\bm \beta _{l}}$ is confined in a compact set $[0,1]$.
B-spline basis functions are commonly used to approximate the unknown smooth functions owing to its desirable numerical stability in practice (\cite{De01}). We thus adopt such a nonparametric approach to estimate the index functions. More specifically, let ${\bm B}(u) = {\left( {B_{s}}(u):1\leq s\leq J_n\right)^T}$ be a set of B-spline basis functions of order $q$ ($q\geq 2$) with $N_n$ internal knots and $J_n=q+N_n$.  We then approximate $m_l(\cdot)$ by a linear combination of B-spline basis functions ${m _l}(\cdot) \approx {\bm {B}}{(\cdot)^T}{\bm {\lambda} _l }, $ where $\bm \lambda = {\left( \bm \lambda _1^T,...,\bm \lambda _d^T \right)^T}$ is the spline coefficient vector with ${\bm \lambda _l} = {\left( {{\lambda _{ls}}:1 \le s \le {J_n}} \right)^T}$ for $l=1,...,d$.

Let ${\rho _\tau }(u) = u\left\{ {\tau  - I(u \le 0)} \right\}$ be the quantile loss function where $I(\cdot)$ is an indicator function. We obtain the estimators of the spline coefficients $\bm \lambda$ and the loading parameters $\bm \beta$ by minimizing
\begin{equation}\label{eq2}
{\mathcal{L}_{\tau n}}\left( {\bm\lambda ,\bm\beta} \right) =\sum\limits_{i = 1}^n {{\rho _\tau }\left\{ {{Y_i} - \sum\limits_{l = 1}^d {{\bm B}{\left( \bm Z_i^T\bm\beta _{l}\right)^T}{\bm\lambda _{l}}X_{il}} } \right\}}
\end{equation}
subject to the constraint $\left\| {{\bm \beta _{l}}} \right\|_2 = 1$ and $\beta _{l1}>0$. Minimizing (\ref{eq2}) with respect to all unknown quantities requires a very non-standard nonlinear programming and the solution is usually hard to obtain directly. To address this computing difficulty, we consider an iterative procedure to estimate $\bm\beta_{l}$ and $m_l(\cdot)$. The detailed steps are given below.

\emph{Step 0.} Initialization step: Obtain an initial value $\hat {\bm\beta} ^{\left( 0 \right)}$ with $\left\|\hat {\bm\beta} ^{\left( 0 \right)} \right\|_2 = 1$. For example, one may use the profile least squares estimation proposed by \cite{MS15}.

\emph{Step 1.} For a given ${\bm\beta}$, $\bm{\hat\lambda}({\bm\beta})$ can be attained by
$\hat {\bm\lambda }\left( {\bm\beta}  \right) = \mathop {\arg \min }\limits_{\bm\lambda  \in {\mathbb{R}^{d{J_n}}}} {\mathcal{L}_{\tau n}}\left( {\bm\lambda ,{\bm\beta}} \right)$. This leads to ${{\hat m}_l}(\cdot,{\bm\beta} ) = {\bm B}{(\cdot)^T}{{\bm{\hat \lambda }}_l}({\bm\beta})$. The first-order derivative $\dot{m_l}(\cdot)$ can be approximated by the spline functions of one order lower than that of ${m_l}(\cdot)$. That is, ${{\hat  {\dot{m}}}_l}(\cdot,{\bm\beta} ) = {\bm {\dot{B}}}{(\cdot)^T}{{\bm{\hat \lambda }}_l}({\bm\beta} )$ where ${\bm {\dot{B}}}$ is the first order derivative of the basis function ${\bm {{B}}}$.

In order to estimate $\bm\beta$, we consider leaving one component of ${{\bm \beta _{l}}}$ out to acknowledge the constraint $\left\| {{\bm \beta _{l}}} \right\|_2 = 1$. Let $\bm{\beta}_{l,-1}  = {\left( {{\beta _{l2}},...,{\beta _{lp}}} \right)^T}$ be
a $p-1$ dimensional parameter vector after removing $\beta_{l1}$ in $\bm\beta_{l}$. The original loading parameter $\bm\beta_{l}$ can be rewritten as
\begin{equation}\label{eq3}
\bm\beta _{l} = \bm\beta_{l} (\bm{\beta}_{l,-1} ) = {(\sqrt {1 - {{\left\| \bm{\beta}_{l,-1}  \right\|_2}^2}} ,{\bm{\beta}_{l,-1}^T })^T},~~\left\|\bm{ \beta}_{l,-1} \right\|_2^2 < 1.
\end{equation}
It is obvious that $\bm\beta_{l}$ is infinitely differentiable with respect to $\bm\beta_{l,-1} $ and the Jacobian matrix is given by
\[{\bm J_l} \left(\bm\beta_{l,-1}\right) = \frac{{\partial \bm \beta_{l} }}{{\partial \bm{\beta}_{l,-1}}} = \left( \begin{array}{l}
  - {{{\bm{\beta}_{l,-1}^T}} \mathord{\left/
 {\vphantom {{{\bm{\beta}_{l,-1}^T}} {\sqrt {1 - {{\left\|\bm{ \beta}_{l,-1}  \right\|_2}^2}} }}} \right.
 \kern-\nulldelimiterspace} {\sqrt {1 - {{\left\| \bm{\beta}_{l,-1}  \right\|_2^2}}} }} \\
~~~~~~~~~ {\bm I_{{p} - 1}} \\
 \end{array} \right),\]
where $\bm I_{p}$ is a ${p} \times {p}$ identity matrix. Denote ${\bm \beta _{ - 1}} = {\left( {\bm \beta _{1, - 1}^T,...,\bm \beta _{d, - 1}^T} \right)^T}$. Then ${\bm \beta _{ - 1}}$ belongs to
\[{\Theta _{ - 1}}  = \left\{ {\bm \beta_{-1}  = {{\left( {\bm \beta _{l,-1}^T:1 \le l \le d} \right)}^T}: {{\left\| {{\bm \beta _{l, - 1}}} \right\|_2^2}} < 1 ,{\bm \beta _{l,-1}} \in {{R}^{{p}-1}}} \right\}.\]

\emph{Step 2.} Let $\bm\beta  = \bm\beta (\bm{\beta}_{-1} )$ with the aforementioned definition $\bm\beta _{l} = \bm\beta_{l} (\bm{\beta}_{l,-1} )$ for $1\leq l \leq d$. Combining the estimators $\bm{\hat\lambda} _l$, ${{\hat { m}}_l}$ and ${{\hat {\dot m}}_l}$ from Step 1, we may construct the quantile regression estimating equations for $\bm\beta_{-1}$ by setting
${{\partial {\mathcal{L}_{\tau n}}\left( {\hat {\bm\lambda} ,\bm\beta } \right)} \mathord{\left/
 {\vphantom {{\partial {L_{\tau n}}\left( {\hat {\bm\lambda} ,\bm\beta } \right)} {\partial {\bm\beta _{-1}}}}} \right.
 \kern-\nulldelimiterspace} {\partial {\bm\beta _{-1}}}}=\bm 0$. However, the equations involve the discontinuous function ${\psi _\tau }\left( u \right) = {\dot{\rho} _\tau }\left( u \right) = \tau  - I\left( {u \le 0} \right)$. This adds difficulty in computation despite there is linear programming solver (eg. \cite{JLWY03}). To achieve faster and more stable estimation, we consider an induced smoothing method via approximating ${\psi _\tau }\left( \cdot \right)$ by a smooth function ${\psi _{\tau{h}} }\left( \cdot \right)$ (\cite{BW07,CKY15,W06}).  We introduce ${G_h}\left( x \right) = G\left( {{x \mathord{\left/
 {\vphantom {x h}} \right.
 \kern-\nulldelimiterspace} h}} \right)$, where $G\left( x \right) = \int_{u < x} {K\left( u \right)} du$, $K\left( \cdot \right)$ is a kernel function and $h$ is a bandwidth.
 Thus, we construct the approximation function
 ${\psi _{\tau{h}} }\left( \cdot \right)= \tau-1+{G_{h}}\left(\cdot\right)  $. Consequently the smoothed estimating equations are given as
\begin{equation}\label{eq4}
\begin{array}{l}
{\mathcal{R}_{\tau n h}}\left( {\bm\beta_{-1} } \right) =\sum\limits_{i = 1}^n {{\psi _{\tau h} }\left\{ {{Y_i} - \sum\limits_{l = 1}^d {{\bm B}{{\left( \bm Z_i^T\bm\beta _{l}\right)}^T}{\bm{\hat\lambda} _l}\left( \bm\beta  \right){X_{il}}} } \right\}}
\\
~~~~~~~~~~~~~~\times \left[ \begin{array}{l}
\left\{ {{{\hat {\dot m}}_1}\left(\bm Z_i^T\bm\beta _1 ,\bm \beta \right){X_{i1}}\bm J_1^T{\bm Z_i} + \left( {{{\partial \bm{\hat \lambda} {{(\bm \beta )}^T}} \mathord{\left/
 {\vphantom {{\partial \bm{\hat \lambda }{{(\bm \beta )}^T}} {\partial {\bm \beta _{1, - 1}}}}} \right.
 \kern-\nulldelimiterspace} {\partial {\bm \beta _{1, - 1}}}}} \right){\bm D_i}(\bm \beta )} \right\} \\
 ~~~~~~~~~~~~~~~~~~~~~~~~~~~ ~~~~~~~\vdots  \\
 \left\{ {{{\hat {\dot m}}_d}\left( \bm Z_i^T\bm\beta _{d},\bm \beta  \right){X_{id}}\bm J_d^T{\bm Z_i} + \left( {{{\partial \bm{\hat \lambda} {{(\bm \beta )}^T}} \mathord{\left/
 {\vphantom {{\partial \bm{\hat \lambda} {{(\bm \beta )}^T}} {\partial {\bm \beta _{d, - 1}}}}} \right.
 \kern-\nulldelimiterspace} {\partial {\bm \beta _{d, - 1}}}}} \right){\bm D_i}(\bm \beta )} \right\} \\
 \end{array} \right]=\bm 0, \\
 \end{array}
\end{equation}
where ${\bm D_i}(\bm \beta ) = \left( {{D_{i,sl}}({\bm \beta _{l}}),1 \le s \le {J_n},1 \le l \le d} \right)^T$ with ${D_{i,sl}}({\bm \beta _{l}}) = {B_s}(\bm Z_i^T \bm \beta_l){X_{il}}$. Then we may employ Fisher scoring algorithm to solve the equations to obtain the estimates. That is
\begin{equation}\label{eq5}
\begin{array}{l}
\bm\beta _{ - 1}^{(k+1)} = \bm\beta _{ - 1}^{(k)} -{\left[ {{{\partial {{\rm{{\cal R}}}_{\tau nh}}\left( {{\bm\beta _{ - 1}}} \right)} \mathord{\left/
 {\vphantom {{\partial {{\rm{{\cal R}}}_{\tau nh}}\left( {{\bm\beta _{ - 1}}} \right)} {\partial {\bm\beta _{ - 1}}}}} \right.
 \kern-\nulldelimiterspace} {\partial {\bm\beta _{ - 1}}}}} \right]^{ - 1}}{{\rm{{\cal R}}}_{\tau nh}}\left( {{\bm\beta _{ - 1}}} \right){\mid _{{\bm\beta _{ - 1}} = \bm\beta _{- 1}^{(k)}}}.
 \end{array}
 \end{equation}

\emph{Step 3.} Repeat Steps 1 and 2 until convergence, and denote the final estimators as $\hat{\bm \beta}_{-1} $ and $\hat{\bm \lambda}  $. Then, we may apply the formula (\ref{eq3}) to obtain $\hat{\bm \beta} $, and construct the estimators of the nonparametric functions ${m_l}(\cdot)$ as ${{\hat m}_l}(\cdot,\hat{\bm\beta} ) = {\bm B}{(\cdot)^T}{{\bm{\hat \lambda }}_l}(\hat{\bm\beta}), l=1,...,d$.

 \begin{remark}\label{remark:2.0}
Another merit of the induced smoothing method is that we can quickly obtain covariance matrix estimation of $\bm{\hat \beta}$ by utilizing the sandwich formula, which can effectively avoid to estimate the density function of the random error. While in step 1, the minimization is a standard quantile regression problem and thus can be solved very easily using the R function ``rq'' from the package ``quantreg''. Therefore, it is not necessary to consider a smooth objective function in the first step.
\end{remark}

\subsection{Theoretical properties}\label{subsect.2.2}
To establish asymptotic normality and the convergence rate of the proposed estimators, we need some assumptions and notations. First, let ${\bm \beta ^0} = {\left\{ {{{\left( {\bm \beta _{1}^0} \right)}^T},...,{{\left( {\bm\beta _{d}^0} \right)}^T}} \right\}^T}$ be the true parameters in model (\ref{eq1}), where $\bm\beta _{l}^0 = {\left\{ {\beta _{l1}^0,{{\left( {\bm \beta _{l, - 1}^0} \right)}^T}} \right\}^T}$ and $\bm\beta _{l, - 1}^0 = {(\beta _{l2}^0,...,\beta _{l{p_n}}^0)^T}$ for $1\le l \le d$. Here the subscript $n$ in $p_n$ is used to make it explicit that the dimension of loading parameters $p_n$ may depend on $n$. Write $\left\|  g \right\|_2 ={\left\{ {\int {{g^2}(x)dx} } \right\}^{{1 \mathord{\left/
 {\vphantom {1 2}} \right.
 \kern-\nulldelimiterspace} 2}}}$ to be the $L_2$ norm of a function $g$, and we focus on the space $\mathcal{M}$ as a collection of functions with finite $L_2$ norm on $[0,1]^d\times {R}^d$ given by
\[\mathcal{M}=\left\{ {g(\bm u,\bm x) = \sum\limits_{l = 1}^d {{g_l}({u_l}){x_l},E{g_l^2}{{(\bm Z^T\bm \beta_l)}} < \infty } } \right\}\]
with $\bm u=(u_1,...,u_d)^T$ and $\bm x=(x_1,...,x_d)^T.$ For $1\le k \le p_n$, we assume that $g_k^0$ is a minimizer in $\mathcal{M}$ for the following optimization problem,
\[\begin{array}{l}
\mathbb{P}({Z_k}) = g_k^0\left( {\bm U({\bm\beta ^0}),\bm X} \right)\\
~~~~~~~~= \sum\limits_{l = 1}^d {g_{l,k}^0(\bm Z^T\bm\beta _{l}^0){X_l}}  \\
~~~~~~~~=\arg {\min _{g \in \mathcal{M}}}E{ \left[f\left(0|\bm X, \bm Z\right)\left( {{Z_k} - g\left( {\bm U\left( {{\bm \beta ^0}} \right),\bm X} \right)} \right)^2\right]} ,
\end{array}\]
where $\bm U({\bm\beta ^0})=(\bm Z^T\bm \beta_1^0,...,\bm Z^T\bm \beta_d^0)^T$. This defines a quadratic projection of ${\bm Z}$ to be $\mathbb{P}({\bm Z})=\{\mathbb{P}({Z_{1}}),\cdots,\mathbb{P}({Z_{p_n}})\}^T $. Next, let $\bm {\tilde{Z}}=\bm Z-\mathbb{P}({\bm {Z}})$, ${\bm A^{ \otimes 2}} = \bm A{\bm A^T}$ for any matrix $\bm A$, ${\bm M_n}({\bm \beta_{-1} ^0}) = {\sum\limits_{i = 1}^n {{{\left[ {\left( {{{\dot m}_l}\left( \bm Z_i^T \bm\beta_{l}^0\right){X_{il}}\bm J_l^{0T}{{\tilde {\bm Z}}_i}} \right)_{l = 1}^d} \right]}^{ \otimes 2}}} }$ and
\[{\bm H_n}({\bm \beta_{-1} ^0}) = {\sum\limits_{i = 1}^n {{{\left[ f\left(0|\bm X_i,\bm Z_i\right) {\left( {{{\dot m}_l}\left( \bm Z_i^T \bm\beta_{l}^0\right){X_{il}}\bm J_l^{0T}{{\tilde {\bm Z}}_i}} \right)_{l = 1}^d} \right]}^{ \otimes 2}}} }\]
with $\bm J_l^0=\bm J_l\left(\bm\beta^0_{l,-1}\right)$ for $1\leq l \leq d$. Suppose that ${\lim _{n \to \infty }}n^{-1}{ \bm M_n}({\bm \beta_{-1} ^{0}})=\mathbb{M} $ and ${\lim _{n \to \infty }}n^{-1}{ \bm H_n}({\bm \beta_{-1} ^{0}})= \mathbb{H}$ are positive definite. Let $r$ ($r \geq 2$) be the order of smoothness of the nonparametric functions $m_l (\cdot)$ given in condition (C2) of the Appendix. We denote $a_n\ll b_n$ if $a_n/b_n=o(1)$. We first present the consistency and asymptotic normality of $\hat{\bm\beta}_{-1}$.

\begin{theorem}\label{theorem:2.1}
 Suppose that conditions (C1)--(C7) in the Appendix B hold, and ${n^{{1 \mathord{\left/
 {\vphantom {1 {(2r + 2)}}} \right.
 \kern-\nulldelimiterspace} {(2r + 2)}}}}  \ll {J_n}  \ll {n^{{1 \mathord{\left/
 {\vphantom {1 4}} \right.
 \kern-\nulldelimiterspace} 4}}}$. If $n^{-1}p_n^{3}=o(1)$, then $\forall \bm e_n\in {R}^{d(p_n-1)}$ such that $\bm e_n^T \bm e_n=1$, we have

\noindent (i) $\left\| {{{\bm{\hat \beta} }_{ - 1}} - {\rm{ }}\bm{\beta} _{ - 1}^0} \right\|_2 = {O_p}\left( \sqrt {{p_n}/n}  \right)$;

\noindent (ii) $\bm e _n^T \bm M_n^{-1/2}\left( {\bm \beta_{-1} ^0} \right)\bm H_n\left( {\bm \beta_{-1} ^0} \right) \left( {{{\hat {\bm \beta} }_{ - 1}} - \bm \beta _{ - 1}^0} \right)\mathop  \to \limits^d N\left( 0, \tau(1-\tau) \right)$,
  \end{theorem}
 \noindent where $\mathop  \to \limits^d$ means the convergence in distribution.

 Following \cite{De01}, for any nonparametric function $m_l$ satisfying Condition (C2) in the Appendix B, there exists a best spline approximation function $m_l^0(u)=\bm B(u)^T \bm \lambda_l^0 $ such that $\mathop {\sup }\limits_{u\in [0,1]} \left| {m _l}(u) - m_l^0(u) \right| =O({J_n^{ - r}})$ for some integer $r\ge 2$. Let $\bm \Psi_n({\bm \beta ^0})={n^{ - 1}}\sum\nolimits_{i = 1}^n { {\bm D_i}({\bm \beta ^0})} {\bm D_i}{({\bm \beta ^0})^T}$ and ${\bm \Omega _n}({\bm \beta ^0}) = {n^{ - 1}}\sum\nolimits_{i = 1}^n f(0\left| {\bm X_i, \bm Z_i} \right.){{\bm D_i}({\bm \beta ^0})} {\bm D_i}{({\bm \beta ^0})^T}$. Suppose ${\lim _{n \to \infty }}{\bm \Omega _n}({\bm \beta ^0}) =\bm\Omega$ and ${\lim _{n \to \infty }}{\bm \Psi _n}({\bm \beta ^0}) =\bm\Psi$ are positive definite.
 Let $\bm e_l$ be the $d\times 1$ vector with the $l$-th element being 1 and other elements being 0, $ {\mathbb{B}}(\bm u) =diag\left(\bm B{(u_1)}^T ,...,\bm B{(u_d)}^T \right)_{d \times d{J_n}}$ with $\bm u=(u_1,...,u_d)^T$ and
\begin{equation}\label{eq6}
\sigma _{nl}^2(u_l) ={n^{ - 1}}\tau(1-\tau)\bm e_l^T\mathbb{B}\left(\bm u \right){\bm \Omega _n^{-1}}({\bm \beta ^0}) \bm \Psi_n({\bm \beta ^0}) {\bm \Omega _n^{-1}}({\bm \beta ^0}) {\mathbb{B}^T}\left( \bm u \right){\bm e_l}.
\end{equation}
The following theorem provides the asymptotic results for the nonparametric estimates.

\begin{theorem}\label{theorem:2.2}
Under conditions (C1)--(C7) in the Appendix B, and ${n^{{1 \mathord{\left/
 {\vphantom {1 {(2r + 2)}}} \right.
 \kern-\nulldelimiterspace} {(2r + 2)}}}}  \ll {J_n}  \ll {n^{{1 \mathord{\left/
 {\vphantom {1 4}} \right.
 \kern-\nulldelimiterspace} 4}}}$, we have for each $1\leq l\leq d$,

 (i)~$|{\hat m}_l( u_l,\bm {\hat\beta} )-{ m}_l( u_l )|={O_p}\left( {\sqrt {{{{J_n}} \mathord{\left/
 {\vphantom {{{J_n}} n}} \right.
 \kern-\nulldelimiterspace} n}}  + J_n^{ - r}} \right)$ uniformly for any $u_l \in [0,1]$;

 (ii) under ${n^{{1 \mathord{\left/
 {\vphantom {1 {(2r + 1)}}} \right.
 \kern-\nulldelimiterspace} {(2r + 1)}}}}  \ll {J_n}  \ll {n^{{1 \mathord{\left/
 {\vphantom {1 4}} \right.
 \kern-\nulldelimiterspace} 4}}}$, $\sigma _{nl}^{ - 1}\left( {{u_l}} \right)\left\{{{\hat m}_{l }}\left( {{u_l},\bm {\hat\beta} } \right)-{m_l^0}({u_l}) \right\}\mathop  \to \limits^d N\left( {0,1} \right)$.
   \end{theorem}

Define $\hat{\mathbb{P}}_n({Z_{ik}})=\bm D_i\left( \hat{\bm\beta}\right)^T \left\{\sum\limits_{i = 1}^n {\hat w_i\bm D_i(\hat{\bm\beta}) \bm D_i(\hat{\bm\beta})^T}\right\}^{-1}  \sum\limits_{i = 1}^n {\hat w_i\bm D_i(\hat{\bm\beta}) Z_{ik}}$ with $\hat w_i=h^{-1}K\left(\hat\varepsilon_i/h \right)$ and $\hat\varepsilon_i=Y_i-\sum\limits_{l = 1}^d {{\hat m_l}} ({\bm Z_i^T}{\hat{\bm \beta} _l}){X_{il}}$. Let $\hat{\mathbb{P}}_n({\bm Z}_i)=\{\hat{\mathbb{P}}_n({Z_{i1}}),...,\hat{\mathbb{P}}_n({Z_{ip_n}})\}^T $, $\hat {\bm Z}_i=\bm Z_i-\hat{\mathbb{P}}_n({\bm {Z}}_i)$, ${\mathbb{H}_{\tau n}}\left(\hat{\bm\beta}_{-1}\right) =\sum\limits_{i = 1}^n \hat w_i {\left[ {\left(\hat {\dot m}_l\left( \bm Z_i^T \hat{\bm\beta}_{l}\right){X_{il}}\hat{\bm J}_l^{T}\hat {\bm Z}_i \right)_{l = 1}^d} \right]^{ \otimes 2}}$ and ${\mathbb{M}_{\tau n}}  \left(\hat{\bm\beta}_{-1}\right) =\sum\limits_{i = 1}^n {\left[ {\psi _{\tau} }\left\{ \hat\varepsilon_i \right\}{\left(\hat {\dot m}_l\left( \bm Z_i^T \hat{\bm\beta}_{l}\right){X_{il}}\hat{\bm J}_l^{T}\hat {\bm Z}_i \right)_{l = 1}^d} \right]^{ \otimes 2}}$ with $\hat{\bm J}_l= \bm J_l\left(\hat{\bm\beta}_{l,-1}\right)$.
Under the same conditions of Theorem \ref{theorem:2.1}, we can show $\frac{1}{n}{\mathbb{H}_{\tau n}}\left(\hat{\bm\beta}_{-1}\right) \mathop  \to \limits^p \mathbb{H}$ and $\frac{1}{n}{\mathbb{M}_{\tau n}}\left(\hat{\bm\beta}_{-1}\right) \mathop  \to \limits^p \tau \left( {1 - \tau } \right)\mathbb{M}$ as $n\rightarrow \infty$, where $\mathop  \to \limits^p$ denotes the convergence in probability.

 \begin{remark}\label{remark:2.1}

Based on the iterative formula (\ref{eq5}) and above results, we apply the following sandwich formula to consistently estimate the asymptotic covariance of $\hat{\bm\beta}_{-1}$
\begin{equation}\label{eq7}
Cov \left(\hat{\bm\beta}_{-1}\right)= {\mathbb{H}_{\tau n}^{-1}}\left(\hat{\bm\beta}_{-1}\right) {\mathbb{M}_{\tau n}}\left(\hat{\bm\beta}_{-1}\right){\mathbb{H}_{\tau n}^{-1}}\left(\hat{\bm\beta}_{-1}\right).
\end{equation}
 Furthermore, we define $\mathbb{\hat J} =  \oplus _{l = 1}^d{\bm {\hat J}_l} = {\rm{diag}}({\bm {\hat J}_1},...,{\bm {\hat J}_d})$ as the direct sum of Jacobian matrices ${\bm {\hat J}_1},...,{\bm {\hat J}_d}$ with dimension $dp_n\times d(p_n-1)$. Then, we can obtain the estimated asymptotic covariance of $\hat{\bm\beta} $ by $Cov \left(\hat{\bm\beta}\right)= \hat{\mathbb{J}}Cov \left(\hat{\bm\beta}_{-1}\right)\hat{\mathbb{J}}^T$.
\end{remark}

\begin{remark}\label{remark:2.2}
Define $\mathbb{D}_{\tau n}=\sum\limits_{i = 1}^n  {{\psi _{\tau}^2 }\left\{ \hat\varepsilon_i \right\} \bm D_i(\hat{\bm\beta})\bm D_i(\hat{\bm\beta})^T} $ and ${\mathbb{C}_{\tau n}} =\sum\limits_{i = 1}^n \hat w_i \bm D_i(\hat{\bm\beta})\bm D_i(\hat{\bm\beta})^T $. Under the conditions of Theorem \ref{theorem:2.2}, we can show $n^{-1}{\mathbb{C}_{\tau n}}\mathop  \to \limits^p \bm\Omega$ and $n^{-1}{\mathbb{D}_{\tau n}} \mathop  \to \limits^p \tau \left( {1 - \tau } \right)\bm\Psi$ as $n\rightarrow \infty$. Thus, variance of $\hat m_l (u_l, \hat {\bm\beta})$ can be consistently estimated by
\begin{equation}\label{eq8}
Var \left( \hat m_l (u_l, \hat {\bm\beta})\right)=\bm e_l^T\mathbb{B}\left( \bm u \right) {\mathbb{C}_{\tau n}^{-1}} {\mathbb{D}_{\tau n}} {\mathbb{C}_{\tau n}^{-1}} {\mathbb{B}^T}\left( \bm u \right){\bm e_l}.
\end{equation}

\end{remark}

\section{Penalized Estimation for High-dimensional Loading Parameters}\label{sect.3}
So far all covariates $\bm Z$ in model (\ref{eq1}) are assumed to be important for predicting the response variable. However, the true model is often unknown. On one hand, the fitted models may be seriously biased and non-informative if important predictors are omitted; on the other hand, including spurious covariates may unnecessarily increase the complexity and further reduce the estimation efficiency. Thus, it is a fundamental issue to select variables for the VICM when there is no prior knowledge of the true model form. In particular, we consider estimation when facing a diverging number of loading parameters. As usual we assume the model sparsity in the sense that most of the components of $\bm\beta$ are essentially zero. For selecting important variables and estimating them simultaneously, penalized robust estimating equations are developed as
\begin{equation}\label{eq9}
{\mathcal{R}_{\tau nh}}\left( {\bm\beta_{-1} } \right)-n{\bm b_{{\alpha _1}}}\left( {{\bm \beta _{ - 1}}} \right)=\bm 0,
\end{equation}
where $\bm b_{\alpha _1}\left(\bm \beta _{- 1} \right)=  [\dot{p}_{\alpha _1}\left( {\left| \beta _{12} \right|} \right){\rm sgn} \left( \beta _{12}\right),...,\dot{p}_{\alpha _1}\left( {\left| \beta _{1p_n} \right|} \right){\rm sgn} \left( \beta _{1p_n}\right),...,\dot{p}_{\alpha _1}\left( {\left| \beta _{dp_n} \right|} \right)$

\noindent $\times{\rm sgn} \left( \beta _{dp_n}\right)]$ is a $d(p_n-1)$ vector with $\mathop{\rm sgn}\left( t \right) = I\left( {t > 0} \right) - I\left( {t < 0} \right)$ and $\dot p _{\alpha _1}(\cdot)$ is the first order derivative of the SCAD penalty function, defined by
\[{\dot{p}_{{\alpha _1}}}(x) = {\alpha _1}\left\{ {I\left( {x \le {\alpha _1}} \right) + \frac{{{{\left( {a{\alpha _1} - x} \right)}_ + }}}{{\left( {a - 1} \right){\alpha _1}}}I\left( {x > {\alpha _1}} \right)} \right\},\]
where $a > 2$, ${{p}_{{\alpha _1}}}(0)=0$ and $\alpha _1$ is a nonnegative penalty parameter which regulates the complexity of the model.  It is easy to see that $\dot{p}_{\alpha _1}(\left| x \right|) $ is close to zero if $\left| x \right|$ is large. Thus little extra bias is introduced by the penalty term. Meanwhile, $\dot{p}_{\alpha _1}(\left| x \right|) $ should be large when $\left| x \right|$ is close to zero, which results in these small components being shrunk to zero.  An iterative majorize-minorize (MM) algorithm proposed by \cite{HL05} can be incorporated to estimate $\bm\beta_{-1}$ in estimating equations (\ref{eq9}). Specifically, for a fixed $\alpha_1$, we can obtain the estimate $\bar{\bm\beta} _{\alpha_1, - 1} $ of $\bm\beta _{ - 1} $ using the following iterative procedure
\begin{equation}\label{eq10}
\begin{array}{l}
\bm\beta _{\alpha_1, - 1}^{(k+1)} = \bm\beta _{\alpha_1, - 1}^{(k)} - \left[ {\partial {\mathcal{R}_{\tau nh}}\left( \bm\beta _{ - 1} \right)/\partial {\bm\beta _{ - 1}} - n{\bm\Delta _{{\alpha _1}}}} \right]^{ - 1} \times \\
~~~~~~~~~~~~~~~~~~~~~~~~~\left(  {\mathcal{R}_{\tau nh}}\left( {\bm\beta_{-1} } \right)-n{\bm b_{{\alpha _1}}}\left( {{\bm \beta _{ - 1}}} \right)  \right)\mid _{\bm\beta _{ - 1}=\bm\beta _{\alpha_1, - 1}^{(k)}}, \\
 \end{array}
 \end{equation}
where ${\bm\Delta _{{\alpha _1}}} = diag\left( {\frac{{{{\dot p}_{{\alpha _1}}}\left( {\left| {{\beta _{12}}} \right|} \right)}}{{\kappa  + \left| {{\beta _{12}}} \right|}},...,\frac{{{{\dot p}_{{\alpha _1}}}\left( {\left| {{\beta _{1{p_n}}}} \right|} \right)}}{{\kappa  + \left| {{\beta _{1{p_n}}}} \right|}},...,\frac{{{{\dot p}_{{\alpha _1}}}\left( {\left| {{\beta _{d{p_n}}}} \right|} \right)}}{{\kappa  + \left| {{\beta _{d{p_n}}}} \right|}}} \right)$ and $\kappa$ is a small number such as $10^{-6}$.
The above iterative formula is similar to the MM algorithm of \cite{HL05}, and its convergence can be similarly justified using their proposition 3.3 under the stationary and continuity assumptions.

In general, we define the true coefficients as $\bm\beta_{l,-1}^{0}= \left(\left(\bm\beta_{l,-1}^{0(1)}\right)^T, \left(\bm\beta_{l,-1}^{0(2)}\right)^T\right)^T $ with $\bm\beta_{l,-1}^{0(1)}={\left( \beta _{l2}^{0},...,\beta _{ls_l} ^{0}\right)^T}$ and $\bm\beta_{l,-1}^{0(2)}={\left( \beta _{l(s_l+1)}^{0},...,\beta _{lp_n} ^{0}\right)^T}$, $\beta_{lj}^{0}\neq 0, j=2,...,s_l$ and $\beta_{lj}^{0}=0,j=s_l+1,...,p_n$, $\bm\beta _{ - 1}^{0(1)}=\left(\left(\bm\beta _{1, - 1}^{0(1)}\right)^T,...,\left(\bm\beta _{d, - 1}^{0(1)}\right)^T\right)^T$ and $\bm\beta _{- 1}^{0(2)}=\left(\left(\bm\beta _{1, - 1}^{0(2)}\right)^T,...,\left(\bm\beta _{d, - 1}^{0(2)}\right)^T\right)^T$. Correspondingly, we also divide $\bar{\bm\beta}_{\alpha_1l,-1}$ into two parts, $\bar{\bm\beta}_{\alpha_1l,-1} = \left((\bar{\bm\beta}_{\alpha_1l,-1}^{(1)})^T,(\bar{\bm\beta}_{\alpha_1l,-1}^{(2)})^T\right)^T$ with $\bar{\bm\beta} _{\alpha_1l, - 1}^{(1)} = {\left( {{\bar{\beta} _{\alpha_1l2}},...,{\bar{\beta} _{\alpha_1ls_l}}} \right)^T}$ and $\bar{\bm\beta}_{\alpha_1l, - 1}^{(2)} = {\left( {{\bar{\beta} _{\alpha_1l(s_l + 1)}},...,{\bar{\beta} _{\alpha_1lp_n}}} \right)^T}$. Define $\bar{\bm\beta} _{\alpha_1, - 1}^{(1)}=\left(\left(\bar{\bm\beta} _{\alpha_11, - 1}^{(1)}\right)^T,...,\left(\bar{\bm\beta} _{\alpha_1d, - 1}^{(1)}\right)^T\right)^T $ and $\bar{\bm\beta} _{\alpha_1, - 1}^{(2)}=\left(\left(\bar{\bm\beta} _{\alpha_11, - 1}^{(2)}\right)^T,...,\left(\bar{\bm\beta} _{\alpha_1d, - 1}^{(2)}\right)^T\right)^T$. Here we assume the number of nonzero components in $\bm\beta_l$ is fixed for $l=1,...,d$, namely, $s_l$ does not vary with $n$.

We establish the following main results for the penalized estimation.

\begin{theorem}\label{theorem:3.1}
Under conditions (C1)--(C12) in the Appendix B, and ${n^{{1 \mathord{\left/
 {\vphantom {1 {(2r + 2)}}} \right.
 \kern-\nulldelimiterspace} {(2r + 2)}}}}  \ll {J_n}  \ll {n^{{1 \mathord{\left/
 {\vphantom {1 4}} \right.
 \kern-\nulldelimiterspace} 4}}}$. If $n^{-1}p_n^{3}=o(1)$ as $n\rightarrow \infty$, we have $\left\| {{{\bm{\bar \beta} }_{\alpha _1, - 1}} - \bm{\beta} _{- 1}^0} \right\|_2 = {O_p}\left( \sqrt {{p_n}} \left( {{n^{{{ - 1} \mathord{\left/
 {\vphantom {{ - 1} 2}} \right.
 \kern-\nulldelimiterspace} 2}}} + {a_n}} \right) \right)$.
  \end{theorem}

Let ${\lim _{n \to \infty }}\frac{1}{n}{ \bm M_n}({\bm \beta_{-1} ^{0(1)}})=\mathbb{M}^{(1)} $ and ${\lim _{n \to \infty }}\frac{1}{n}{ \bm H_n}({\bm \beta_{-1} ^{0(1)}}) =\mathbb{H}^{(1)}$, where ${\bm {M}_n}\left( {\bm \beta _{ - 1}^{0(1)}} \right)$ and ${\bm {H}_n}\left( {\bm \beta _{ - 1}^{0(1)}} \right)$ are $\sum\limits_{l = 1}^d {({s_l}-1)}  \times \sum\limits_{l = 1}^d {({s_l}-1)} $ sub-matrices of $\bm{M}_n\left( {\bm \beta _{ - 1}^{0}} \right)$ and $\bm{H}_n\left( {\bm \beta _{ - 1}^{0}} \right)$ corresponding to ${\bm\beta}_{-1}^{0(1)}$.

\begin{theorem}\label{theorem:3.2}
 Under conditions (C1)--(C12) in the Appendix B, and ${n^{{1 \mathord{\left/
 {\vphantom {1 {(2r + 2)}}} \right.
 \kern-\nulldelimiterspace} {(2r + 2)}}}}  \ll {J_n}  \ll {n^{{1 \mathord{\left/
 {\vphantom {1 4}} \right.
 \kern-\nulldelimiterspace} 4}}}$. If $\alpha_1 \to 0 $, $\sqrt {{n \mathord{\left/
 {\vphantom {n {{p_n}}}} \right.
 \kern-\nulldelimiterspace} {{p_n}}}}\alpha_1 \to \infty$ and $n^{-1}p_n^{3}=o(1)$ as $n\rightarrow \infty$, with probability tending to one, the consistent estimator $\bar{\bm\beta}_{ {\alpha _1},- 1}$ satisfies

\noindent (i) $\bar{\bm\beta}_{\alpha _1l,-1}^{(2)}=\bm 0$ for $1\leq l\leq d$;

\noindent (ii) $\sqrt n \left(\bar{\bm \beta}_{\alpha _1, - 1}^{(1)}-\bm \beta _{ - 1}^{0(1)}\right) \mathop  \to \limits^d   N\left( \bm {0},\tau(1-\tau)(\mathbb{H} ^{(1)})^{-1}\mathbb{M}^{(1)}  (\mathbb{H} ^{(1)})^{-1}\right). $
  \end{theorem}

Now we define $\mathbb{J}^0 =  \oplus _{l = 1}^d{\bm J_l^0} = {\rm{diag}}({\bm J_1^0},...,{\bm J_d^0})$ as the direct sum of Jacobian matrices ${\bm J_1^0},...,{\bm J_d^0}$ with dimension
$dp_n\times d(p_n-1)$. For $1\le l \le d$, $\bm{ \beta}_{l}$ can be estimated by $\bm{\bar \beta}_{\alpha_1,l}=(\bar \beta_{\alpha_1,l1},...,\bar \beta_{\alpha_1,lp_n})^T$ with ${\bar \beta}_{\alpha_1,l1}={\left( {1 - \sum\limits_{k = 2}^{p_n} {\bar \beta _{\alpha_1,lk}^2} } \right)^{{1 \mathord{\left/
 {\vphantom {1 2}} \right.
 \kern-\nulldelimiterspace} 2}}}$. Based on Theorem~\ref{theorem:3.2} (ii), we can use the multivariate delta method to obtain the asymptotic normality of $\bm{\bar\beta}_{\alpha_1}^{(1)}=(\bm{\bar\beta}_{\alpha_1,1}^{(1)T},...,\bm{\bar\beta}_{\alpha_1,d}^{(1)T})^T$
 with $\bar{\bm\beta} _{\alpha_1,l}^{(1)} = {\left( {{\bar{\beta} _{\alpha_1l1}}, {\bar{\beta} _{\alpha_1l2}},...,{\bar{\beta} _{\alpha_1ls_l}}} \right)^T}$. That is,
\[\sqrt n (\bm{\bar \beta}_{\alpha_1}^{(1)} - \bm \beta^{0(1)})\mathop  \to \limits^d N\left( {\bm 0},\tau(1-\tau)\mathbb{J}^{0(1)}(\mathbb{H} ^{(1)})^{-1}\mathbb{M}^{(1)}  (\mathbb{H} ^{(1)})^{-1}\mathbb{J}^{0(1)T}\right),\]
where $\bm \beta^{0(1)}={\left( \beta _{11}^{0},...,\beta _{1s_1} ^{0},...,\beta _{d1} ^{0},...,\beta _{ds_d} ^{0}\right)^T}$ and $\mathbb{J}^{0(1)}$ is sub-matrix of $\mathbb{J}^{0}$ corresponding to ${\bm\beta}^{0(1)}$.

\begin{remark}\label{remark:3.1}
 Theorem~\ref{theorem:3.1} provides the convergence rate of $\bm{\bar \beta} _{\alpha _1, - 1}$. Theorem~\ref{theorem:3.2} indicates that $\bm{\bar \beta} _{\alpha _1, - 1}$ is consistent in variable selection and has the oracle property when the number of loading parameters diverges. These results provide a theoretical guarantee for the application of our proposed estimation for high-dimensional quantile regression VICM. Based on the iterative procedure (\ref{eq10}), we obtain the following
sandwich formula to estimate the asymptotic covariance matrix of $\bar{\bm\beta}_{\alpha_1,-1}$ by
 \begin{equation}\label{eq11}
 Cov \left(\bar{\bm\beta}_{\alpha_1,-1}\right)= {\bar{\mathbb{H}}_{\tau n}^{-1}}\left(\bar{\bm\beta}_{\alpha_1,-1}\right) {{\mathbb{M}}_{\tau n}}\left(\bar{\bm\beta}_{\alpha_1,-1}\right){\bar{\mathbb{H}}_{\tau n}^{-1}}\left(\bar{\bm\beta}_{\alpha_1,-1}\right)
 \end{equation}
where $\bar{\mathbb{H}}_{\tau n}\left(\bar{\bm\beta}_{\alpha_1,-1}\right) =\mathbb{H}_{\tau n}\left(\bar{\bm\beta}_{\alpha_1,-1}\right)+n\bm\Delta_{\alpha_1}$, $\mathbb{M}_{\tau n}$ and ${{\mathbb{H}}_{\tau n}} $ are defined in (\ref{eq7}).
 \end{remark}

\section{Identification of linear components in quantile regression VICM}\label{sect.4}
In varying index coefficient models, identification of linear interaction components is also an important issue. A hypothesis test may be conducted to distinguish linear functions from nonparametric functions. \cite{MS15} proposed a generalized likelihood ratio test for this purpose. However, the classical significance tests may not be so desirable in high dimensional settings for computational and theoretical concerns. In this paper, we develop a penalized procedure based on the SCAD penalty to investigate whether there is a linear interaction effect between $\bm Z^T\bm\beta_l$ and $X_l$.

Let $\ddot{m}_l$ be the second derivative of ${m}_l$. It is obvious that ${\left\| {{{\ddot m}_l}} \right\|_2}=0$ if ${m}_l$ is a linear function for $1\leq l \leq d$. Thus, by shrinking ${\left\| {{{\ddot m}_l}} \right\|_2}$ towards zero, we can automatically identify the linear and non-linear components in model (\ref{eq1}). Note that ${\left\| {{{\ddot m}_l}} \right\|_2}={\left\{ {\int {\ddot m_l^2\left( x \right)} dx} \right\}^{1/2}}$ after the basis approximation can be equivalently written as $\sqrt {\bm \lambda _l^T{\bm D}{\bm\lambda _l}}\equiv {\left\| {{\bm\lambda _l}} \right\|_{{\bm D}}}$ due to the well-known algebraic property of the B-spline approximation, where $\bm D$ is a $J_n \times J_n$ matrix with the $(k, k')$ entry being $\int_0^1 {{\ddot{B}_{k}}\left( x \right){\ddot{B}_{k'}}\left( x \right)dx}$. This consideration leads to solving the following minimization problem
\begin{equation}\label{eq12}
\bar {\bm\lambda }=   \mathop {\arg \min }\limits_{\bm\lambda  \in {\mathbb{R}^{d{J_n}}}} \mathcal{L}_{\tau n}^*\left( {\bm\lambda ,\bar{\bm \beta}_{\alpha _1} } \right)
\equiv
 \mathop {\arg \min }\limits_{\bm\lambda  \in {\mathbb{R}^{d{J_n}}}} \left\{ {{\mathcal{L}_{\tau n}}\left( {\bm\lambda ,\bar{\bm \beta}_{\alpha _1} } \right) +n\sum\limits_{l = 1}^d {{p_{{\alpha _2}}}\left( {{{\left\| {{\bm \lambda _l}} \right\|}_{{\bm D}}}} \right)}} \right\},
\end{equation}
where $p_{\alpha _2}(\cdot)$ is the SCAD penalty with a penalty parameter $\alpha _2$ and $\bar{\bm \beta}_{\alpha _1}$ is given in sect. \ref{sect.3}. This is a still complicated nonlinear programming problem and we use the ``ucminf'' function in R software to find the minimum of (\ref{eq12}) using numerical computing methods. This R function was developed by Hans Bruun Nielsen and Stig Bousgaard Mortensen for
general-purpose unconstrained non-linear optimization. The algorithm is of quasi-Newton type with BFGS updating of the inverse Hessian and soft line search with a trust region monitoring method. Using this numerical computing routine we may not need to apply the induced smoothing any more.

\begin{remark}\label{remark:4.0}

We may combine two types of penalties in the objective function (\ref{eq2}) to perform variable selection for loading parameters and detect linear/nonlinear simultaneously, that is ${\mathcal{Q}_{\tau n}}\left( {\bm\lambda ,\bm\beta} \right)={\mathcal{L}_{\tau n}}\left( {\bm\lambda ,\bm\beta} \right)+n\sum\limits_{l = 1}^d {\sum\limits_{j = 2}^{{p_n}} {{p_{{\alpha _1}}}\left( {\left| {{\beta _{lj}}} \right|} \right)} } +n\sum\limits_{l = 1}^d {{p_{{\alpha _2}}}\left( {{{\left\| {{\bm \lambda _l}} \right\|}_{{\bm D}}}} \right)}$. However, $\bm\lambda$ depends on $\bm\beta$, which indicates that we can not simultaneously obtain the estimators of $\bm\lambda$ and $\bm\beta$ by minimizing ${\mathcal{Q}_{\tau n}}\left( {\bm\lambda ,\bm\beta} \right)$. To address this difficulty, an iterative procedure is proposed to select the loading parameter and detect the linear/nonlinear components. That is, for a given $\lambda$, we develop the penalized robust estimating equations (\ref{eq9}) to select the loading parameters, and minimize (\ref{eq12}) to detect the linear/nonlinear components for a given $\bm\beta$.

\end{remark}

Let $\bm{\bar \lambda }=\left(\bm{\bar \lambda }_{1}^T,...,\bm{\bar \lambda }_{d}^T\right)^T$ be the minimizer of $\mathcal{L}_{\tau n}^*\left( {\bm\lambda ,\bar{\bm \beta}_{\alpha _1} } \right)$. Consequently, the estimator of ${{\bar m}_{l}}(\cdot )$ is ${{\bar m}_{l}}(\cdot  ) = {\bm B}{(\cdot )^T}{{\bm{\bar \lambda }}_{l}}$ for $1 \leq l \leq d$. Without loss of generality, we suppose that $m_l$ is truly nonlinear for $1\leq l \leq d_1$ and is linear for $d_1 +1\leq l \leq d$. We have the following theoretical results.

\begin{theorem}\label{theorem:4.1}
 Suppose that conditions (C1)--(C12) in the Appendix B hold, together with ${n^{{1 \mathord{\left/
 {\vphantom {1 {(2r + 2)}}} \right.
 \kern-\nulldelimiterspace} {(2r + 2)}}}}  \ll {J_n}  \ll {n^{{1 \mathord{\left/
 {\vphantom {1 4}} \right.
 \kern-\nulldelimiterspace} 4}}}$ and $\alpha_2 \rightarrow 0$, we have for each $1\leq l\leq d$, $|{\bar m}_l( u_l,\bm {\bar\beta}_{\alpha_1} )-{ m}_l( u_l )|={O_p}\left( {\sqrt {{{{J_n}} \mathord{\left/
 {\vphantom {{{J_n}} n}} \right.
 \kern-\nulldelimiterspace} n}}  + J_n^{ - r}} \right)$ uniformly for any $u_l \in [0,1]$.
  \end{theorem}

\begin{theorem}\label{theorem:4.2}
In addition to the conditions in Theorem 4.1,
and we further assume ${( {{{ \sqrt {{{{J_n}} \mathord{\left/
 {\vphantom {{{J_n}} n}} \right.
 \kern-\nulldelimiterspace} n}} }} + J_n^{ - r}} )^{ - 1}}{\alpha _2} \to \infty $. Then with probability approaching 1, $\|{\bar{\bm\lambda}}_l\|_{\bm D}=0$ and $\bar m_{l}$ is a linear function for $1+d_1 \leq l \leq d$.
  \end{theorem}

\begin{remark}\label{remark:4.1}
When $\alpha_2 \rightarrow 0$, Theorem \ref{theorem:4.1} and Theorem \ref{theorem:2.2} indicates that the nonparametric function estimates $\hat m_l(\cdot)$ and $\bar m_l(\cdot)$ attain the same convergence rate, which shows that the addition of a penalty term does not impact the asymptotic properties of $\bar m_l(\cdot)$. Theorem \ref{theorem:4.2}  establishes that the proposed method can identify linear functions consistently. These results provide a solid support for our proposed identification procedure.
\end{remark}

\section{Numerical illustration}\label{sect.5}
In this section, simulation studies and real data analysis are provided to assess the finite sample performance of the proposed estimation methods.

\subsection{Selection of tuning parameters}\label{subsect.5.1}

\begin{figure}
\centering
\includegraphics[scale=0.32]{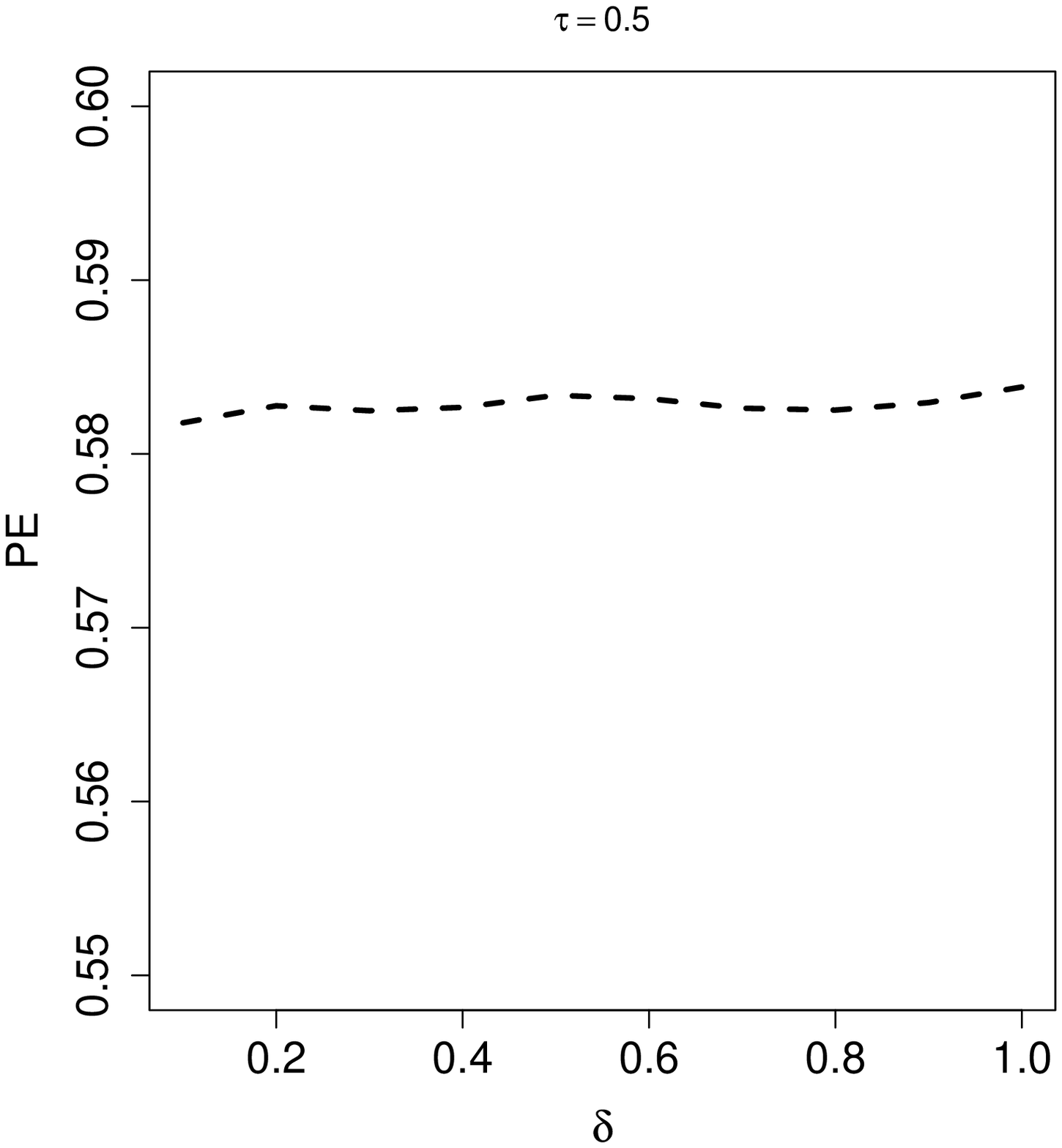}
\includegraphics[scale=0.32]{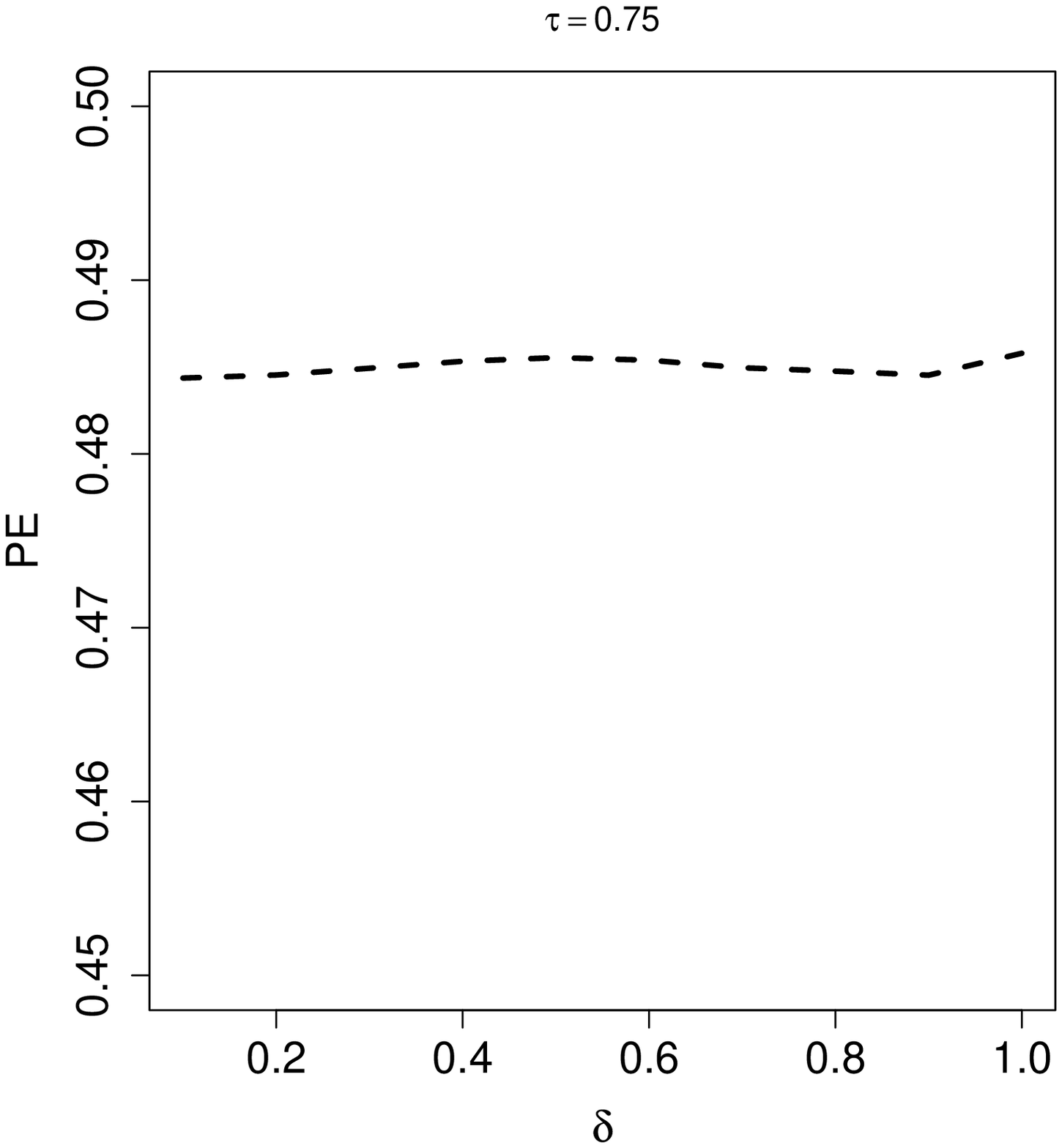}
\caption{Prediction error from 5-fold cross-validation with different bandwidth $h=n^{-\delta}$ with $\delta=0.1,0.2,...,1$.}
\label{figure1}
\end{figure}

For all numerical studies we use the cubic spline ($q = 4$) to approximate nonparametric functions $m_l (\cdot)$ in our simulations. We choose the number of interior knots as ${N_n} = \left[ {{n^{1/(2q + 1)}}} \right]$ that satisfies theoretical requirement, where $[a]$ stands for the largest integer not greater than $a$. The kernel function $K\left( \cdot \right)$ given in sect.~\ref{sect.2} is set as the second-order Bartlett kernel ($\nu = 2$), that is,
\[K\left( u \right) = \frac{3}{{4\sqrt 5 }}\left( {1 - {{{u^2}} \mathord{\left/
 {\vphantom {{{u^2}} 5}} \right.
 \kern-\nulldelimiterspace} 5}} \right)I\left( {\left| u \right| \le \sqrt 5 } \right).\]
To examine the dependence on the bandwidth, we conduct a sensitivity analysis for the selection of $h$ for a selected setting in the following. Let $\{T^v, v=1,\cdots,5\}$ be a random partitioning with size $n/5$ of the full data set $T=({T-T^v})\bigcup {T^v}$ and set $T-T^v$ and set ${T^v}$ be the cross validated training and test sets respectively for $v=1,...,5$. The prediction error (PE) from the 5-fold cross-validation is given by
\[{\rm{PE}} = n^{-1}\sum\limits_{v = 1}^5 {\sum\limits_{\left( {{Y_{i}},{\bm{X}_{i}}} ,\bm{Z}_{i}\right) \in {T^v}} {\rho _{\tau}}\left(  Y_i -   \sum\limits_{l = 1}^d {{\hat m_l^{(v)}}} ({\bm Z_i^T}{\hat{\bm \beta} _l^{(v)}}){X_{il}} \right)  },\]
where $\hat m_l^{(v)}$ and $\hat{\bm\beta}_l^{(v)}$ are estimators of $m_l$ and $\beta_l$ using the training set $T - {T^v}$ for $l=1,...,d$. For quantile levels $\tau= 0.5, 0.75$, we conduct 500 replicates in example 1 given below with normal error distribution. Fig. \ref{figure1} depicts the prediction error from the 5-fold cross-validation with different bandwidth $h=n^{-\delta}, \delta=0.1,0.2,...,1$. It is easy to see that PE does not vary much with different $h$ used, which indicates that the proposed method is not sensitive to the bandwidth $h$. Thus, we fix $h=n^{-0.3}$ in simulation studies to reduce the computational burden. This choice also satisfies the theoretical requirement $n{h^{2 \nu}} \to 0$ with $\nu=2$.

Finally, we utilize a data driven method to select $\alpha_1$ and $\alpha_2$ in the SCAD penalty function. The tuning parameter $\alpha_1$ is used to control the sparsity of the solution and the tuning parameter $\alpha_2$ is to identify the linear functions. Under fixed dimensions, \cite{L12} demonstrated that the Schwartz information criterion (SIC) is consistent in variable selection in a penalized  quantile regression with the SCAD penalty. However, the traditional SIC may not work very well for diverging number of parameters. In this paper, we adopt the following modified SIC (MSIC) to select $\alpha_1$
\[{\rm {MSIC}}\left( {{\alpha _1}} \right) = \log \left( {{\mathcal{L}_{\tau n}}\left( {\hat {\bm\lambda} ,{{\bar {\bm\beta} }_{\alpha _1}}} \right)} \right) + {{d{f_1}{C_n}\log \left( n \right)} \mathord{\left/
 {\vphantom {{d{f_1}{C_n}\log \left( n \right)} {\left( {2n} \right)}}} \right.
 \kern-\nulldelimiterspace} {\left( {2n} \right)}},\]
where $\bar {\bm\beta}_{\alpha _1}$ is the estimated parameter for a given $\alpha _1$, $\hat {\bm\lambda}$ is the unpenalized estimator given in section 2, $df_1$ is the number of nonzero coefficients in $\bar {\bm\beta}_{\alpha _1}$ and $C_n$ is required to be diverging. In our simulations and applications, we choose $C_n$ as $C_n=\max\left\{1,\log \left(\log(dp_n)\right) \right\} $ (\cite{CC08,WLL09}). The optimal $\hat\alpha _1^{opt}$ is defined as
$\hat \alpha _1^{opt} = \mathop {\min }\limits_{{\alpha _1}} \rm {MSIC}\left( {{\alpha _1}} \right)$.
Similarly, for $\alpha_2$, we consider
\[{\rm{MSIC}}\left( {{\alpha _2}} \right) = \log \left( {{\mathcal{L}_{\tau n}}\left( {\bar {\bm\lambda}_{\alpha _2} ,{\bar {\bm\beta} _{\hat\alpha _1^{opt}}}} \right)} \right) + {{d{f_2}{J_n}\log \left( n \right)} \mathord{\left/
 {\vphantom {{d{f_1}{C_n}\log \left( n \right)} {\left( {2n} \right)}}} \right.
 \kern-\nulldelimiterspace} {\left( {2n} \right)}},\]
where $\bar {\bm\lambda}_{{\alpha _2}}$ is the estimated parameter for a given $\alpha _2$, $df_2$ is the number of nonlinear components. Then, we have $\hat \alpha _2^{opt} = \mathop {\min }\limits_{{\alpha _2}} \rm {MSIC}\left( {{\alpha _2}} \right)$. Note that every $m_l(\cdot)$ is characterized by a spline coefficient vector $\bm\lambda_l$ whose dimension is $J_n$. So $df_2 J_n$ is regarded as the dimension of nonlinear function coefficients. Simulation results will confirm that the proposed two MSIC criteria work well for variable selection and identification of linear components.

\subsection{Simulation studies}\label{subsect.5.2}

\tabcolsep=7pt
\begin{table}\scriptsize
\caption{Simulation results ($\times 10^{-2}$) of Bias, MAD, ESD and ASD for SN with $\tau=0.5$ in example 1. LS stands for the method by \cite{MS15} and QR is the proposed quantile regression.}
\label{table1}
\begin{tabular}{cccccccccccc} \noalign{\smallskip}\hline
\multicolumn{1}{c}{method}
&&&\multicolumn{2}{c}{$n=500$}
&&&&\multicolumn{2}{c}{$n=1500$}
\\
\cline{3-6}
\cline{8-11}
&&\multicolumn{1}{c}{Bias}&\multicolumn{1}{c}{MAD}&\multicolumn{1}{c}{ESD}&\multicolumn{1}{c}{ASD}
&&\multicolumn{1}{c}{Bias}&\multicolumn{1}{c}{MAD}&\multicolumn{1}{c}{ESD}&\multicolumn{1}{c}{ASD}

\\
\hline
LS&$\beta_{11}$&0.102&  2.991& 3.770& 3.009&&-0.092&  1.581& 2.061& 1.699\\
&$\beta_{12}$&-0.091 & 4.375 &5.497& 4.224&&0.034 & 2.224 &2.834& 2.363\\
&$\beta_{13}$&-0.362 & 2.129 &2.736 &2.098&&-0.039 & 1.134 &1.461& 1.170\\

&$\beta_{21}$&0.017  &1.389 &1.754 &1.667&&-0.031 & 0.819& 1.037 &0.963\\
&$\beta_{22}$&-0.269 & 2.257& 2.785& 2.598&&0.157 & 1.331& 1.663& 1.492\\
&$\beta_{23}$&0.145 & 2.119 &2.701& 2.959&&-0.339 & 1.232& 1.551 &1.690\\

&$\beta_{31}$&-0.001 & 0.959& 1.211& 1.091&&0.016 & 0.537 &0.661& 0.609\\
&$\beta_{32}$&-0.056 & 0.785& 0.984& 0.931&&-0.016&  0.426& 0.539& 0.516\\
&$\beta_{33}$&0.090 & 1.040& 1.330 &1.322&&-0.005 & 0.567& 0.702 &0.732\\
\hline
QR&$\beta_{11}$&0.182& 3.302& 4.235& 3.595&&-0.078&  1.702& 2.163& 1.851\\
&$\beta_{12}$&-0.037 & 4.775 &6.016& 4.999&&0.057&  2.314& 2.918 &2.566\\
&$\beta_{13}$&-0.503&  2.333 &2.980 &2.484&&-0.064 & 1.216 &1.538 &1.226\\

&$\beta_{21}$&-0.041 & 1.609& 2.036 &2.099&&-0.039&  0.901& 1.150& 1.180\\
&$\beta_{22}$&-0.273 & 2.540& 3.144& 3.267&&0.101&  1.474& 1.846& 1.821\\
&$\beta_{23}$&0.227  &2.429 &3.082& 3.661&&-0.229 & 1.360 &1.721& 2.021\\

&$\beta_{31}$&0.051&  1.096& 1.386& 1.329&&0.039 & 0.590& 0.724& 0.719\\
&$\beta_{32}$&-0.135&  0.941& 1.166& 1.144&&-0.036&  0.475& 0.596& 0.609\\
&$\beta_{33}$&0.192 & 1.309 &1.640& 1.609&&0.001 & 0.633 &0.798& 0.848\\
\hline
\end{tabular}
\end{table}

\tabcolsep=7pt
\begin{table}\scriptsize
\caption{Simulation results ($\times 10^{-2}$) of Bias, MAD, ESD and ASD for $t_3$ with $\tau=0.5$ in example 1. LS stands for the method by \cite{MS15} and QR is the proposed quantile regression.}
\label{table2}
\begin{tabular}{cccccccccccc} \noalign{\smallskip}\hline
\multicolumn{1}{c}{method}
&&&\multicolumn{2}{c}{$n=500$}
&&&&\multicolumn{2}{c}{$n=1500$}
\\
\cline{3-6}
\cline{8-11}
&&\multicolumn{1}{c}{Bias}&\multicolumn{1}{c}{MAD}&\multicolumn{1}{c}{ESD}&\multicolumn{1}{c}{ASD}
&&\multicolumn{1}{c}{Bias}&\multicolumn{1}{c}{MAD}&\multicolumn{1}{c}{ESD}&\multicolumn{1}{c}{ASD}
\\
\hline
LS&$\beta_{11}$&-0.707&  5.128& 6.658& 4.716&&-0.279&  2.520& 3.365 &2.685\\
&$\beta_{12}$&0.003 & 7.000& 9.397& 6.485&&-0.061 & 3.452 &4.442& 3.653\\
&$\beta_{13}$&-0.509  &3.554 &4.902 &3.305&&-0.022 & 1.841& 2.376 &1.830\\

&$\beta_{21}$&-0.420&  2.606& 4.129& 2.775&&0.046&  1.452& 1.843 &1.617\\
&$\beta_{22}$&0.125 & 3.916 &5.142& 4.194&&-0.113 & 2.377& 3.016& 2.532\\
&$\beta_{23}$&-0.222&  3.615 &4.727& 4.932&&-0.281 & 2.071& 2.661& 2.817\\

&$\beta_{31}$&0.079 & 1.637& 2.069 &1.732&&0.139 & 0.891 &1.142& 1.027\\
&$\beta_{32}$&-0.166&  1.340& 1.694 &1.464&&-0.139  &0.697& 0.890& 0.873\\
&$\beta_{33}$&0.117 & 1.736& 2.185 &2.093&&0.074 & 0.890& 1.153& 1.209\\
\hline
QR&$\beta_{11}$&-0.782&  4.049& 5.087& 3.990&&-0.276&  2.007& 2.536& 2.074\\
&$\beta_{12}$&-0.041&  5.636& 7.222& 5.455&&-0.050 & 2.745& 3.465 &2.838\\
&$\beta_{13}$&-0.031&  2.731& 3.504& 2.734&&0.066 & 1.383& 1.736& 1.405\\

&$\beta_{21}$&-0.049&  1.880& 2.354& 2.353&&0.100 & 1.034& 1.311 &1.298\\
&$\beta_{22}$&-0.157&  2.906& 3.626& 3.638&&-0.171 & 1.700 &2.126& 2.025\\
&$\beta_{23}$&-0.102&  2.657& 3.377& 4.159&&-0.152  &1.631& 2.042& 2.227\\

&$\beta_{31}$&0.031 & 1.311& 1.691& 1.450&&0.095 & 0.686& 0.848 &0.791\\
&$\beta_{32}$&-0.116 & 1.124& 1.411 &1.231&&-0.089 & 0.542 &0.664& 0.669\\
&$\beta_{33}$&0.134 & 1.467& 1.811& 1.779&&0.040&  0.706& 0.887 &0.936\\
\hline
\end{tabular}
\end{table}

\tabcolsep=11pt
\begin{table}\scriptsize
\caption{Simulation results of RASE for $m_1$, $m_2$, $m_3$ with $\tau=0.5$ and $n = 500$ in example 1. LS stands for the method proposed by \cite{MS15} and QR is the proposed quantile regression.}
\label{table3}
\begin{tabular}{ccccccccccccccc} \noalign{\smallskip}\hline
\multirow{1}{*}{$n$}
&\multirow{1}{*}{Error}
&\multicolumn{3}{c}{LS}
&&\multicolumn{3}{c}{QR}
\\
\cline{3-5}
\cline{7-9}
&&\multicolumn{1}{c}{$m_1$}&\multicolumn{1}{c}{$m_2$}&\multicolumn{1}{c}{$m_3$}
&&\multicolumn{1}{c}{$m_1$}&\multicolumn{1}{c}{$m_2$}&\multicolumn{1}{c}{$m_3$}
\\
\hline
500&SN&0.097&0.127&0.099&&0.110& 0.145&0.121\\
&$t_3$&0.144&0.177 &0.156&&0.127 &0.154& 0.138\\
&MN&0.150&0.183&0.172&&0.124&0.152&0.138\\
&LA&0.121&0.151&0.134&&0.117&0.146&0.129\\
\hline
1500&SN&0.074 &0.097& 0.061&& 0.087& 0.107 &0.070\\
&$t_3$&0.095& 0.120& 0.095&& 0.088& 0.112 &0.083\\
&MN&0.098 &0.131& 0.100 &&0.087& 0.118& 0.079\\
&LA&0.084 &0.108 &0.077&& 0.081& 0.105& 0.070\\
\hline
\end{tabular}
\end{table}

\textsc{Example 1.} \label{example.1} In this example, our goal is to compare the proposed quantile regression estimator (QR) with the least-squares estimator (LS) for VICM. We generate the random samples from the following model
\begin{equation}\label{eq13}
Y_i=\sum\limits_{l = 1}^d {m_l\left(\bm Z_i^T\bm \beta_l\right)X_{il}}+\sigma\epsilon_i
\end{equation}
where $\sigma=0.5$, $d=p=3$, $X_{i1}=1$, $ (X_{i2},X_{i3})^T$ and $\bm Z_i =(Z_{i1},Z_{i2},Z_{i3})^T$ follow the multivariate normal distributions with mean 0, variance 1 and constant correlation coefficient 0.5. Here we set the true loading parameters as ${\bm \beta _1} = \frac{1}{{\sqrt {14} }}{(2,1,3)^T},$ ${\bm \beta _2} = \frac{1}{{\sqrt {14} }}{(3,2,1)^T}$ and ${\bm \beta _3} = \frac{1}{{\sqrt {14} }}{(2,3,1)^T}$ and set the true coefficient functions as ${m_1}({u_1}) = {{\exp ({u_1})} \mathord{\left/
 {\vphantom {{\exp ({u_1})} 5}} \right.
 \kern-\nulldelimiterspace} 5},{m_2}({u_2}) = \sin (0.5\pi {u_2})$ and ${m_3}({u_3}) = u_3^2$.
In order to investigate the effect of relatively heavy tail error distributions or outliers, we consider the following four different error distributions of $\epsilon_i$: standard normal distribution (SN), $t$-distribution with freedom degree 3 ($t_3$), Laplace distribution (LA) with location parameter 0 and shape parameter 1 and mixed normal distribution ($\rm {MN}(\rho, \sigma_1,\sigma_2)$) which is a mixture of $\rm {N}(0,\sigma_1^2)$ and $\rm {N}(0, \sigma_2^2)$ with weights $1-\rho$ and $\rho$, respectively. In this example, we consider $\rho=0.1$, $\sigma_1=1$ and $\sigma_2=5$. In this example, for the purpose of comparison, we consider $\tau=0.5$ and the sample size $n=500$ and 1500 with 500 simulation replications. For a fixed $\tau=0.5$, we have $Q_{0.5}\left( Y | \bm X, \bm Z \right)=E\left( Y | \bm X, \bm Z \right)=\sum\limits_{l = 1}^d {m_l\left(\bm Z_i^T\bm \beta_l\right)X_{il}}$ since the median of $\epsilon_i $ is zero under the four distributions, $Q_{0.5}\left( \epsilon_i\right)=0$. That is, model (\ref{eq13}) is a special case of model (\ref{eq1}). Therefore, it is fair to compare the proposed quantile regression estimate with the least squares estimate under this setting.

For parametric part, we report the bias (Bias), empirical standard deviation (ESD), calculated as the sample standard deviation of 500 estimates, estimated asymptotic standard deviation (ASD) based on the sandwich formula (\ref{eq7}) and mean absolute deviation (MAD), calculated as the mean absolute deviation of 500 estimates. We compute the root average squared errors (RASE) to measure the accuracy of nonparametric estimators ${ \hat m_l}$
\[{\rm{RASE}}({\hat m}_l) = \sqrt {\frac{1}{n}\sum\limits_{i = 1}^{n}{{{( {{{\hat m}_l}\left( {{u_{il}}} \right) - {m_{l}}\left( {{u_{il}}} \right)} )}^2}} },u_{il}=\bm{Z_i}^T \hat {\bm{\beta}}_l, l=1,2,3.\]
The corresponding results of the proposed quantile regression estimator with $\tau=0.5$ and that of the least-squares estimator are reported in Tables \ref{table1}--\ref{table3}. For space consideration, the simulation results of estimators $\bm{\hat\beta}_l$ for MN and LA are listed in Tables S1 and S2 in the Appendix A of Supplementary Materials. Both mean regression and median regression in this example are consistent to the true parameters and functions as we observe very small errors. Eyeballing the Tables, the performance of the proposed estimation procedure is much more stable than that of the least-squares estimator especially in the cases with non-normal errors, demonstrating the robust feature of this approach. When the underlying distribution for data is different from normal, it may be more reliable to implement our procedure to fit the VICM in practice. Finally, we can see that the estimated asymptotic standard deviation (ASD) is very close to the empirical standard deviation (ESD), especially for $n=1500$. This demonstrates the sandwich covariance formula (\ref{eq7}) works reasonably well.

\tabcolsep=6pt
\begin{table}\scriptsize
\caption{Simulation results of Bias, MAD, ESD and ASD for $\bm \beta_{\tau}$ with $\tau=0.5,0.75$ and $n = 500,1500$ in example 2.}
\label{table6}
\begin{tabular}{cccccccccccc} \noalign{\smallskip}\hline
\multicolumn{1}{c}{$n$}
&&&\multicolumn{2}{c}{$\tau=0.5$}
&&&\multicolumn{2}{c}{$\tau=0.75$}
\\
\cline{3-6}
\cline{8-11}
&&\multicolumn{1}{c}{Bias}&\multicolumn{1}{c}{MAD}&\multicolumn{1}{c}{ESD}&\multicolumn{1}{c}{ASD}
&&\multicolumn{1}{c}{Bias}&\multicolumn{1}{c}{MAD}&\multicolumn{1}{c}{ESD}&\multicolumn{1}{c}{ASD}
\\
\hline
500&$\beta_{\tau,11}$&-0.003&  0.064 &0.082& 0.083&&0.082&  0.092& 0.071& 0.063\\
&$\beta_{\tau,12}$&-0.004&  0.067& 0.085& 0.089&&-0.014 & 0.061& 0.075& 0.068\\
&$\beta_{\tau,13}$&-0.008&  0.053& 0.066& 0.066&&-0.055 & 0.063 &0.057& 0.051\\

&$\beta_{\tau,21}$&0.011 & 0.183& 0.220& 0.144&&-0.043 & 0.173& 0.206& 0.119\\
&$\beta_{\tau,22}$&-0.039&  0.166& 0.211& 0.133&&0.032 & 0.154& 0.191& 0.112\\
&$\beta_{\tau,23}$&-0.069&  0.151& 0.197& 0.120&&-0.068&  0.132& 0.166& 0.093\\

&$\beta_{\tau,31}$&-0.005 & 0.082 &0.106& 0.083&&-0.068&  0.090& 0.090& 0.075\\
&$\beta_{\tau,32}$&-0.001&  0.111& 0.141& 0.112&&-0.012&  0.097& 0.123& 0.097\\
&$\beta_{\tau,33}$&-0.042&  0.113& 0.144& 0.116&&0.077 & 0.114& 0.116 &0.098\\
\hline
1500&$\beta_{\tau,11}$&-0.008&  0.032& 0.040& 0.043&&0.080&  0.080& 0.037& 0.035\\
&$\beta_{\tau,12}$&-0.001&  0.037 &0.046& 0.047&&-0.008 & 0.031 &0.038& 0.038\\
&$\beta_{\tau,13}$&0.003 & 0.025 &0.032& 0.034&&-0.050 & 0.051& 0.030& 0.029\\

&$\beta_{\tau,21}$&-0.030&  0.108& 0.132& 0.094&&-0.076&  0.114& 0.122& 0.082\\
&$\beta_{\tau,22}$&-0.032&  0.104& 0.124& 0.089&&0.023&  0.088 &0.111 &0.079\\
&$\beta_{\tau,23}$&0.008 & 0.080& 0.100& 0.072&&-0.001&  0.069& 0.089& 0.058\\

&$\beta_{\tau,31}$&0.001 & 0.041& 0.051& 0.049&&-0.062&  0.066 &0.048& 0.045\\
&$\beta_{\tau,32}$&0.000 & 0.055& 0.069& 0.065&&0.001 & 0.048& 0.061& 0.059\\
&$\beta_{\tau,33}$&-0.014 & 0.055& 0.072& 0.067&&0.054 & 0.090& 0.061& 0.059\\
\hline
\end{tabular}
\end{table}

\begin{figure}
\centering
\includegraphics[scale=0.2]{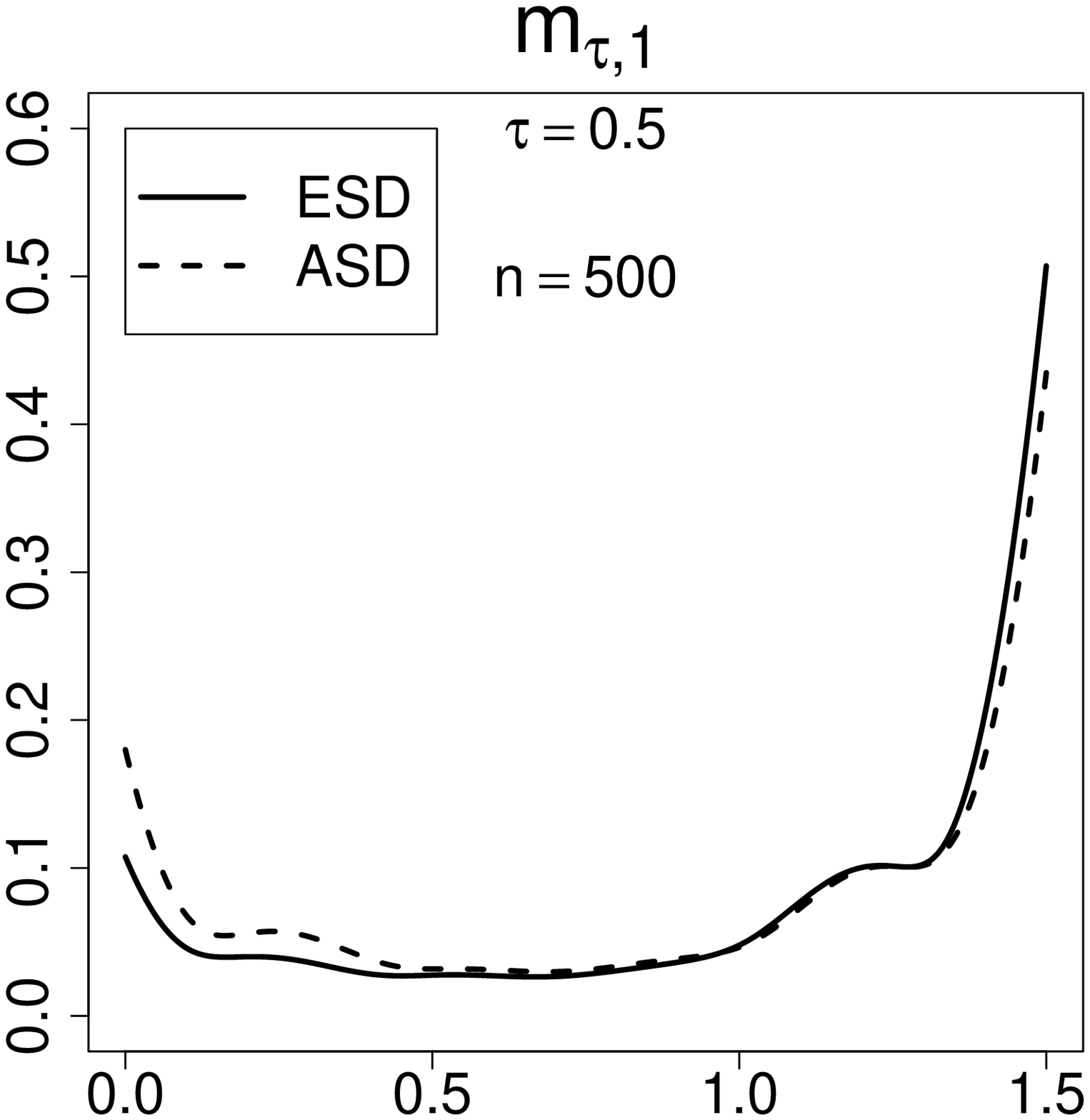}
\includegraphics[scale=0.2]{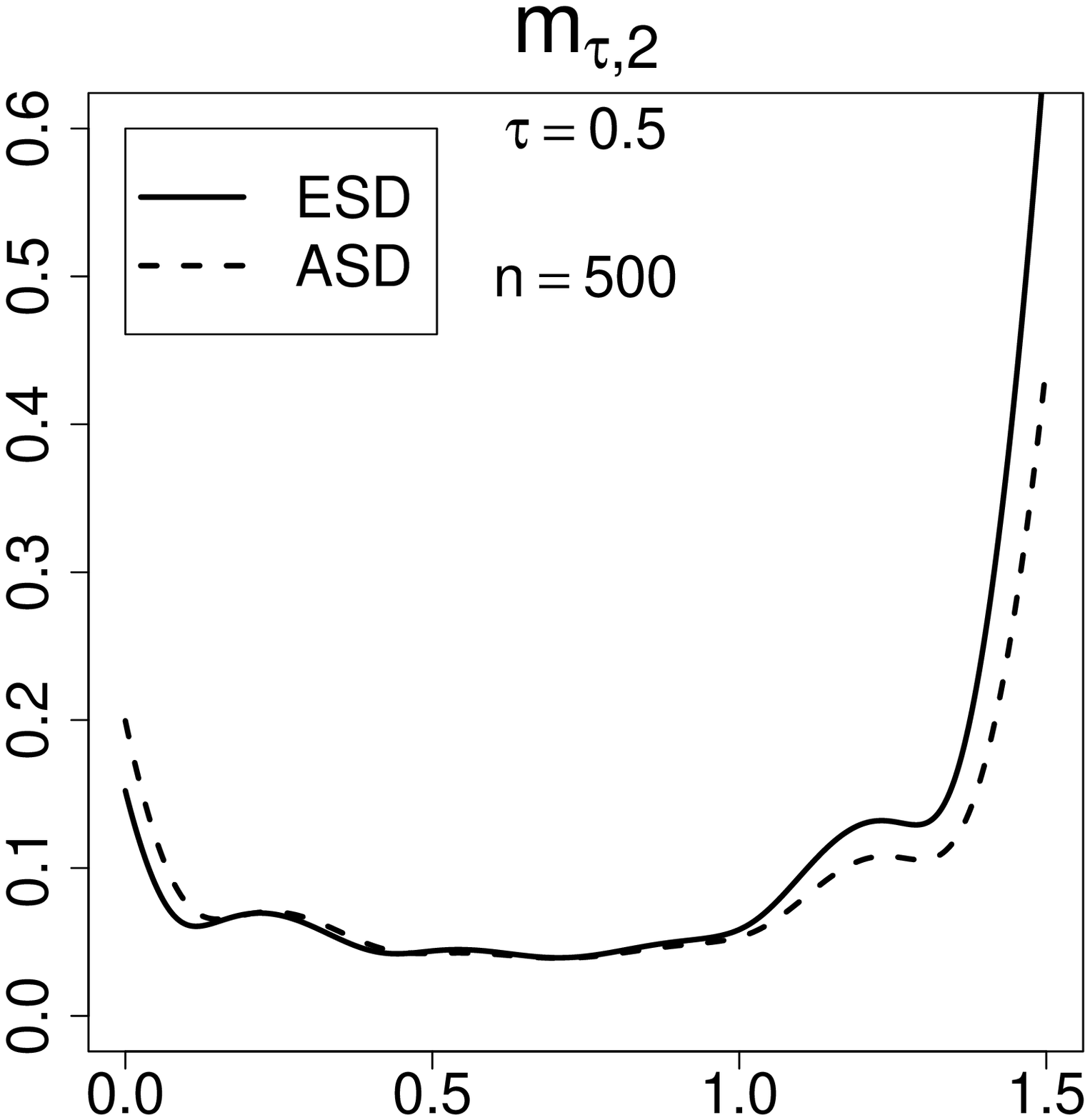}
\includegraphics[scale=0.2]{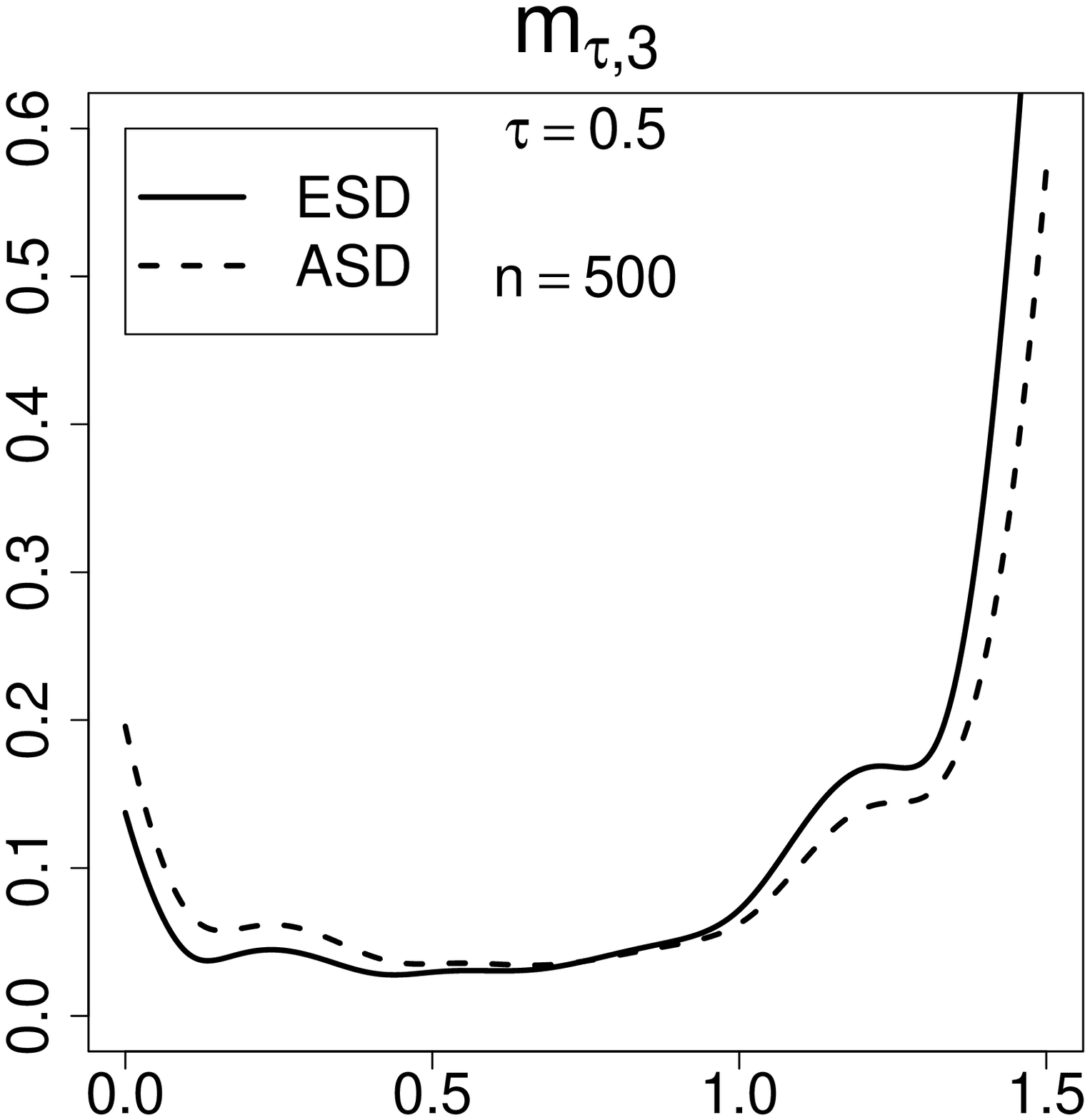}
\includegraphics[scale=0.2]{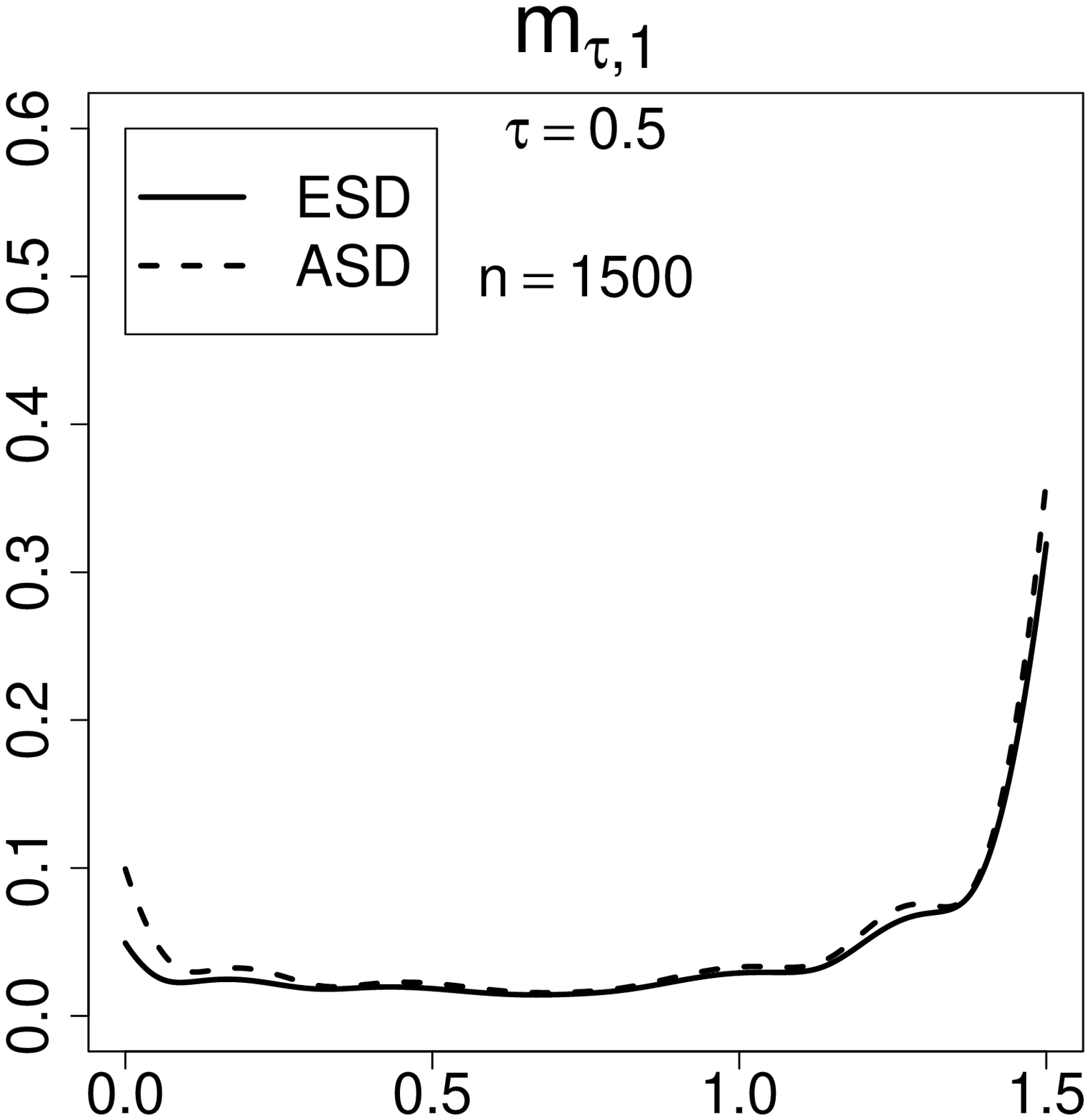}
\includegraphics[scale=0.2]{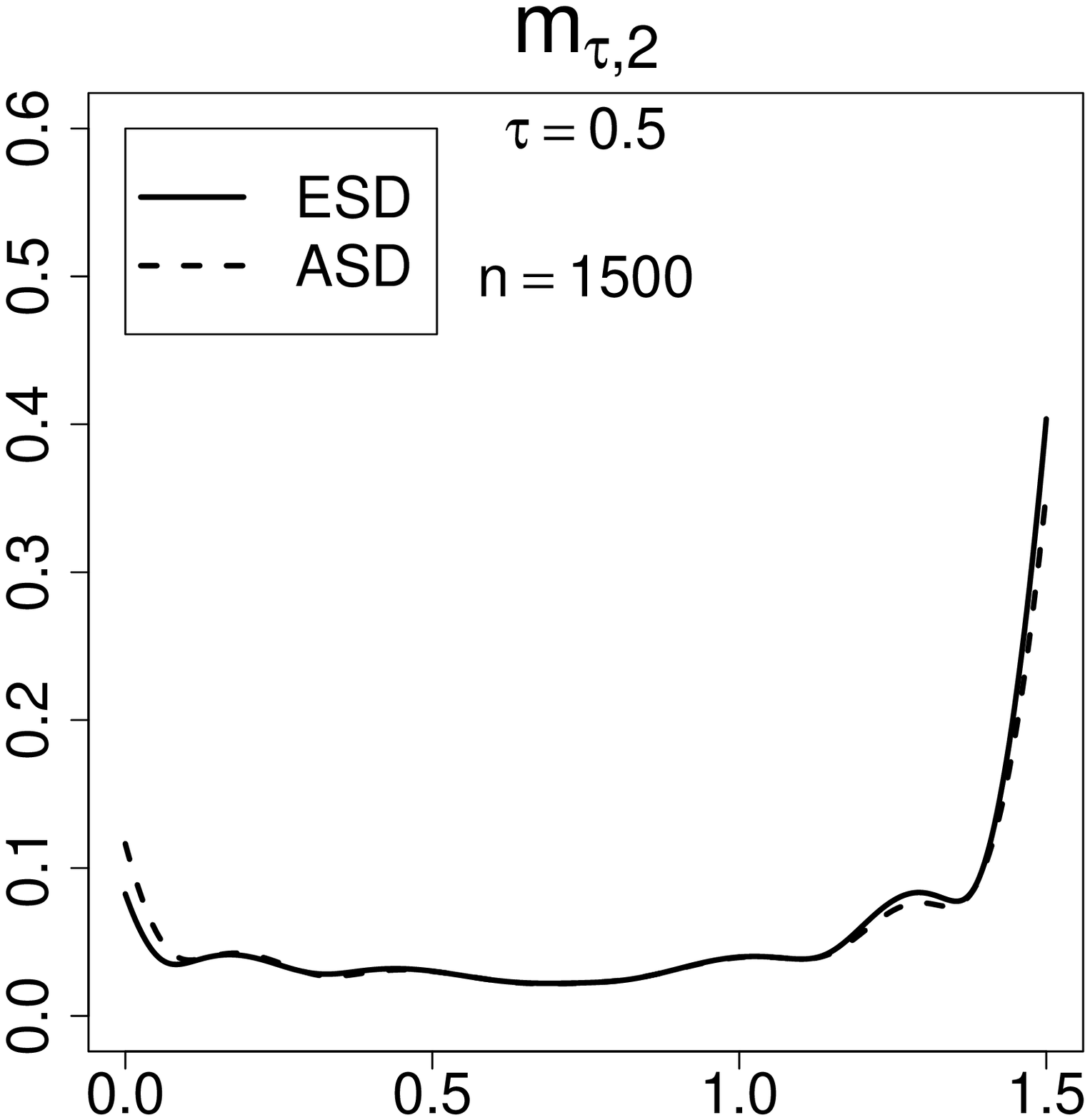}
\includegraphics[scale=0.2]{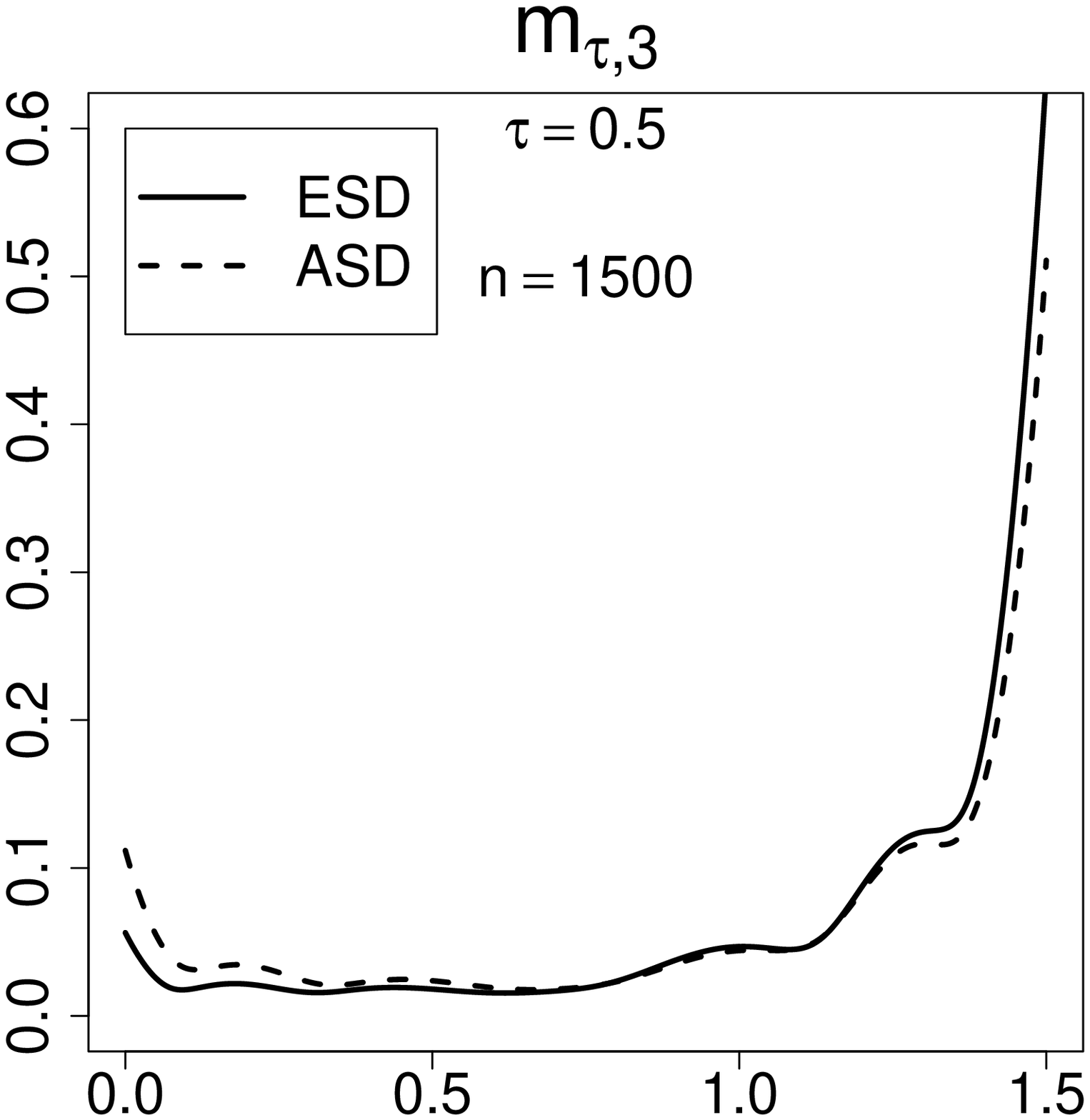}
\includegraphics[scale=0.2]{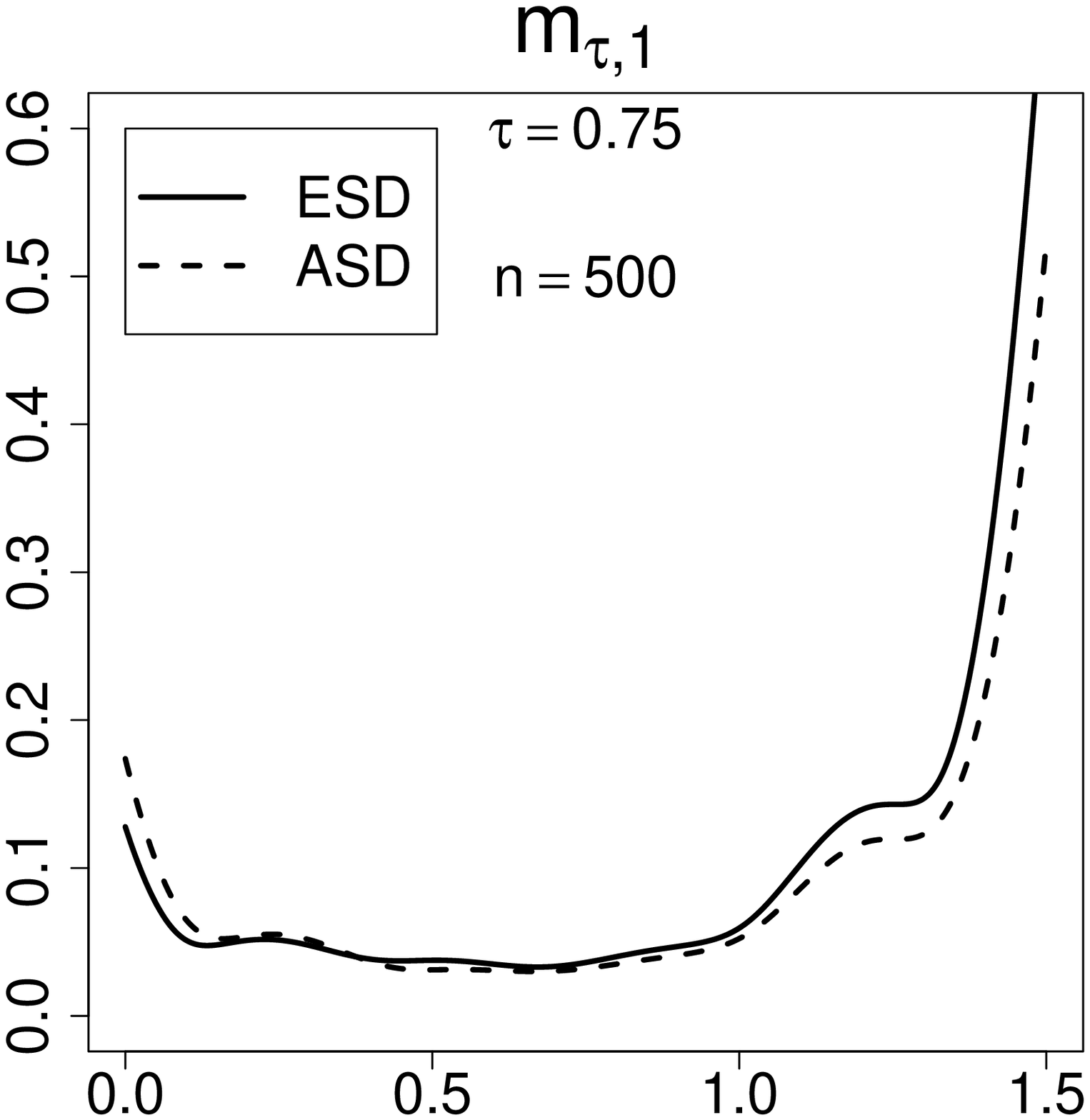}
\includegraphics[scale=0.2]{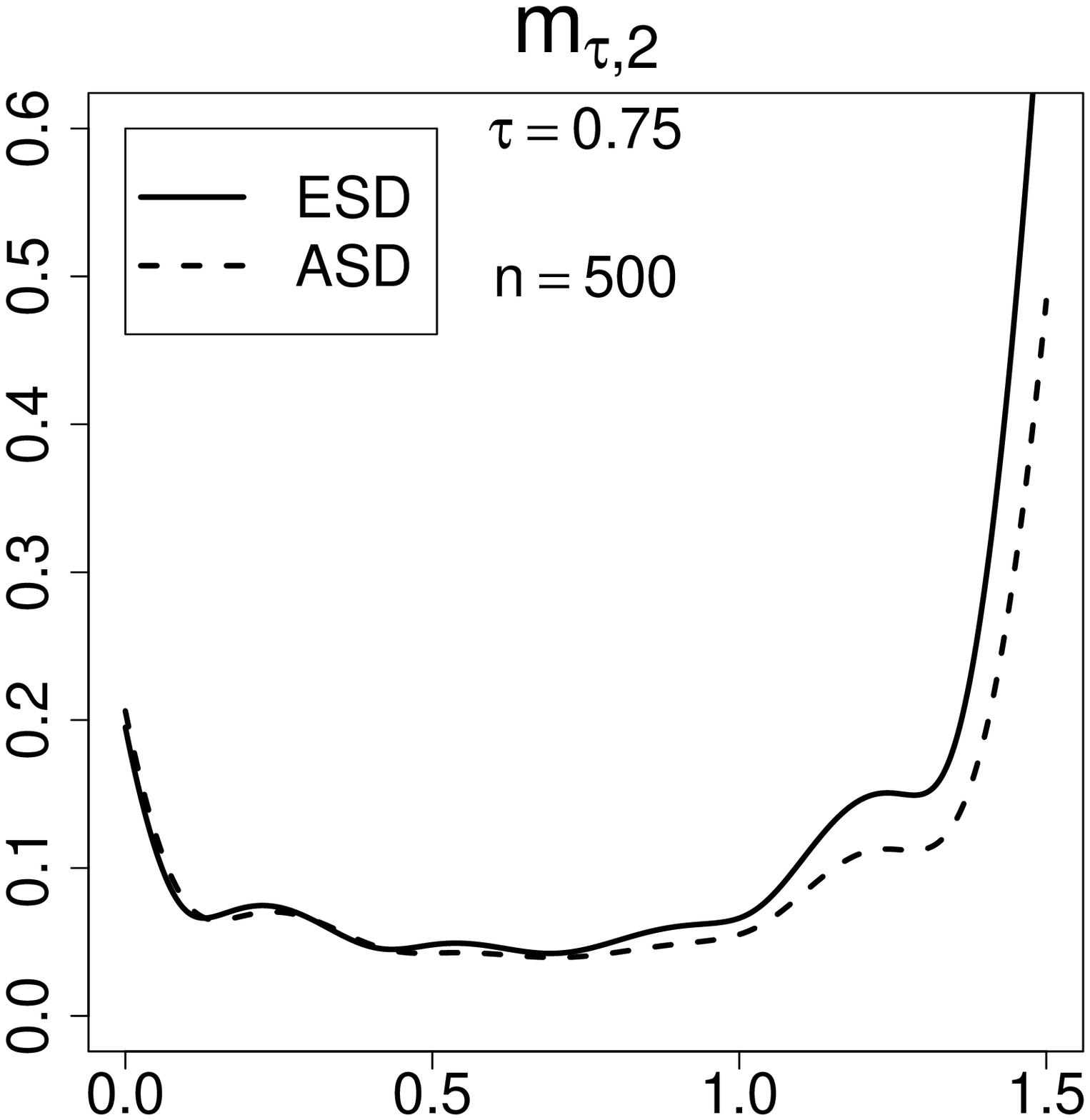}
\includegraphics[scale=0.2]{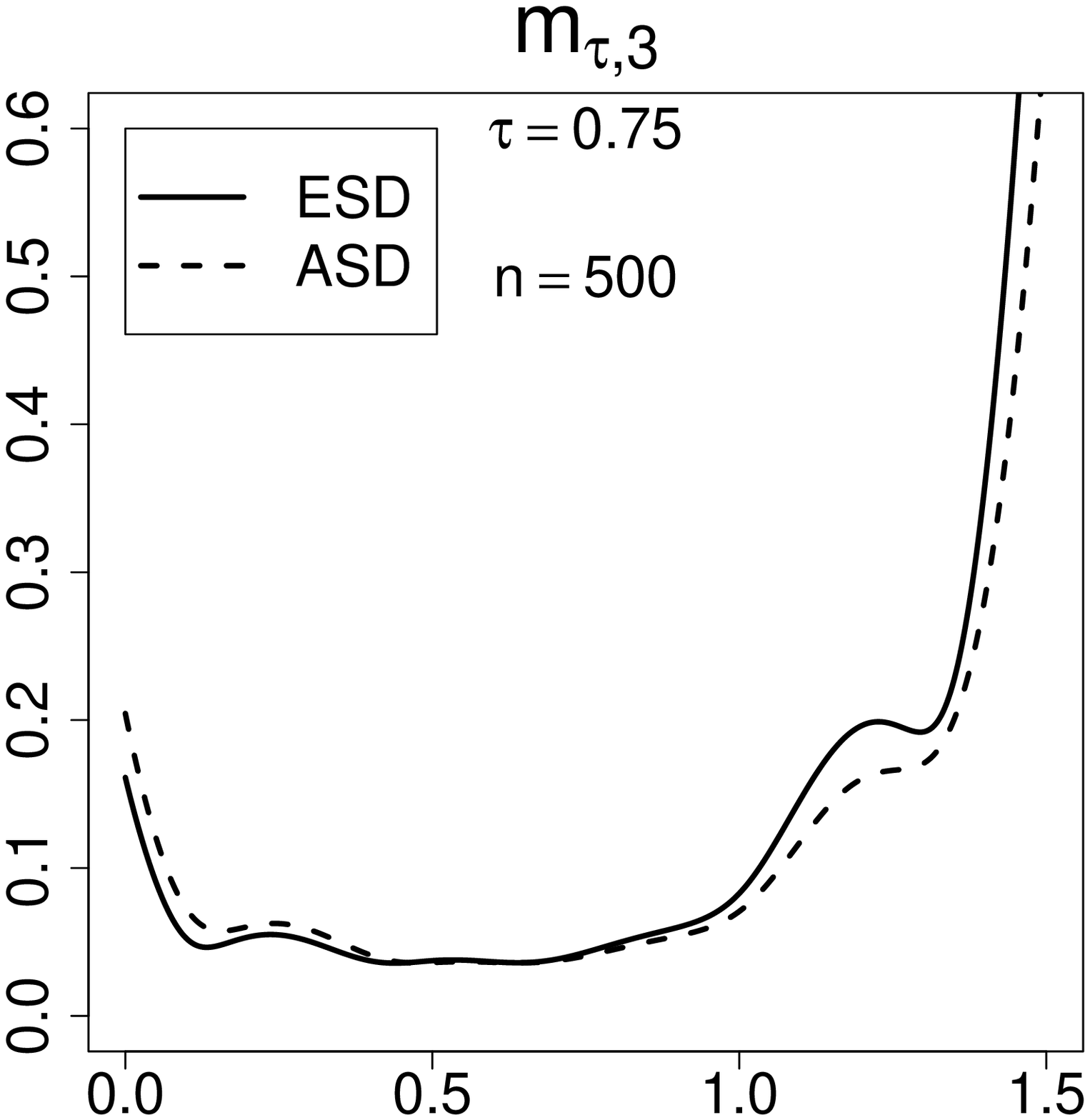}
\includegraphics[scale=0.2]{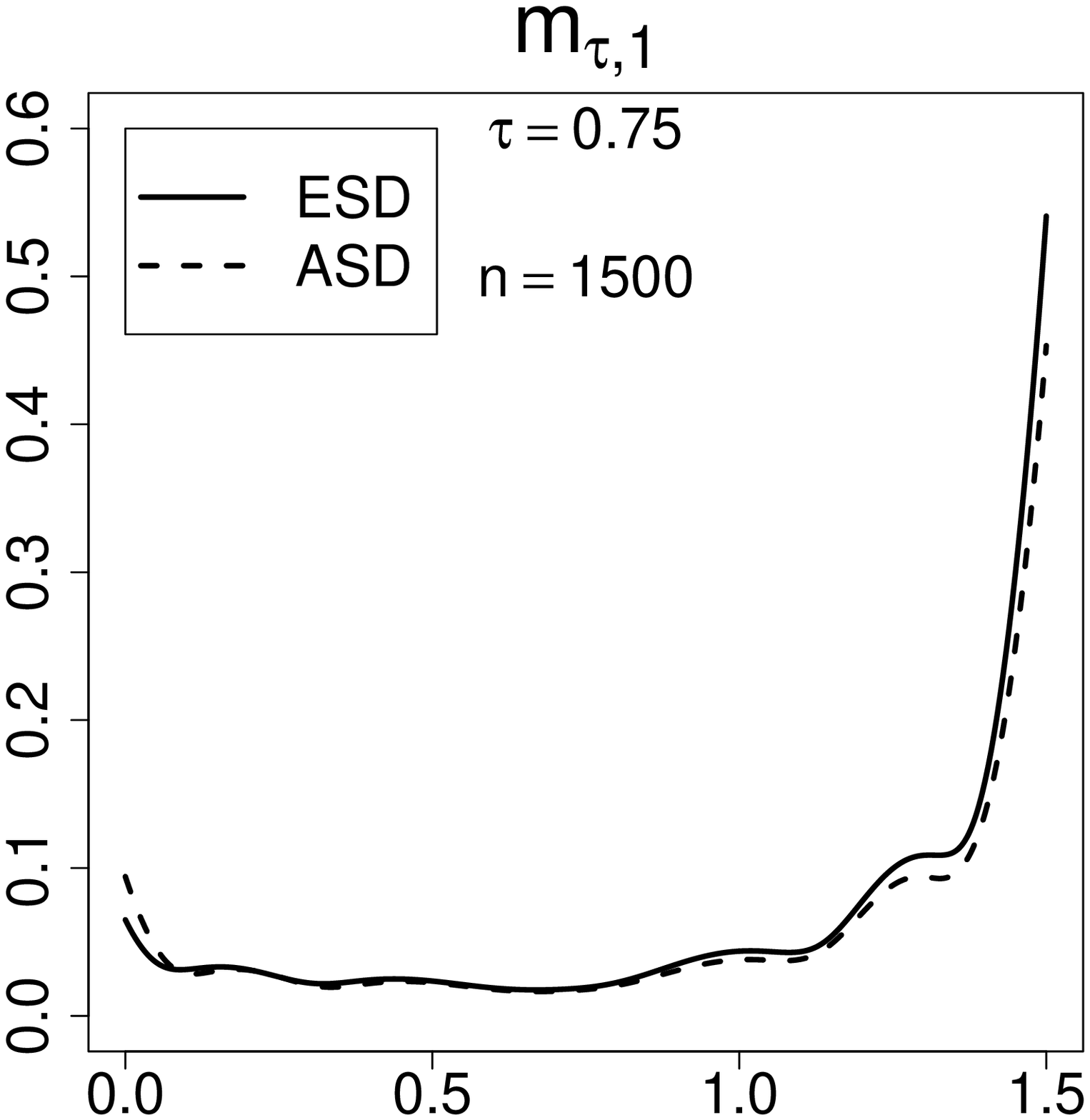}
\includegraphics[scale=0.2]{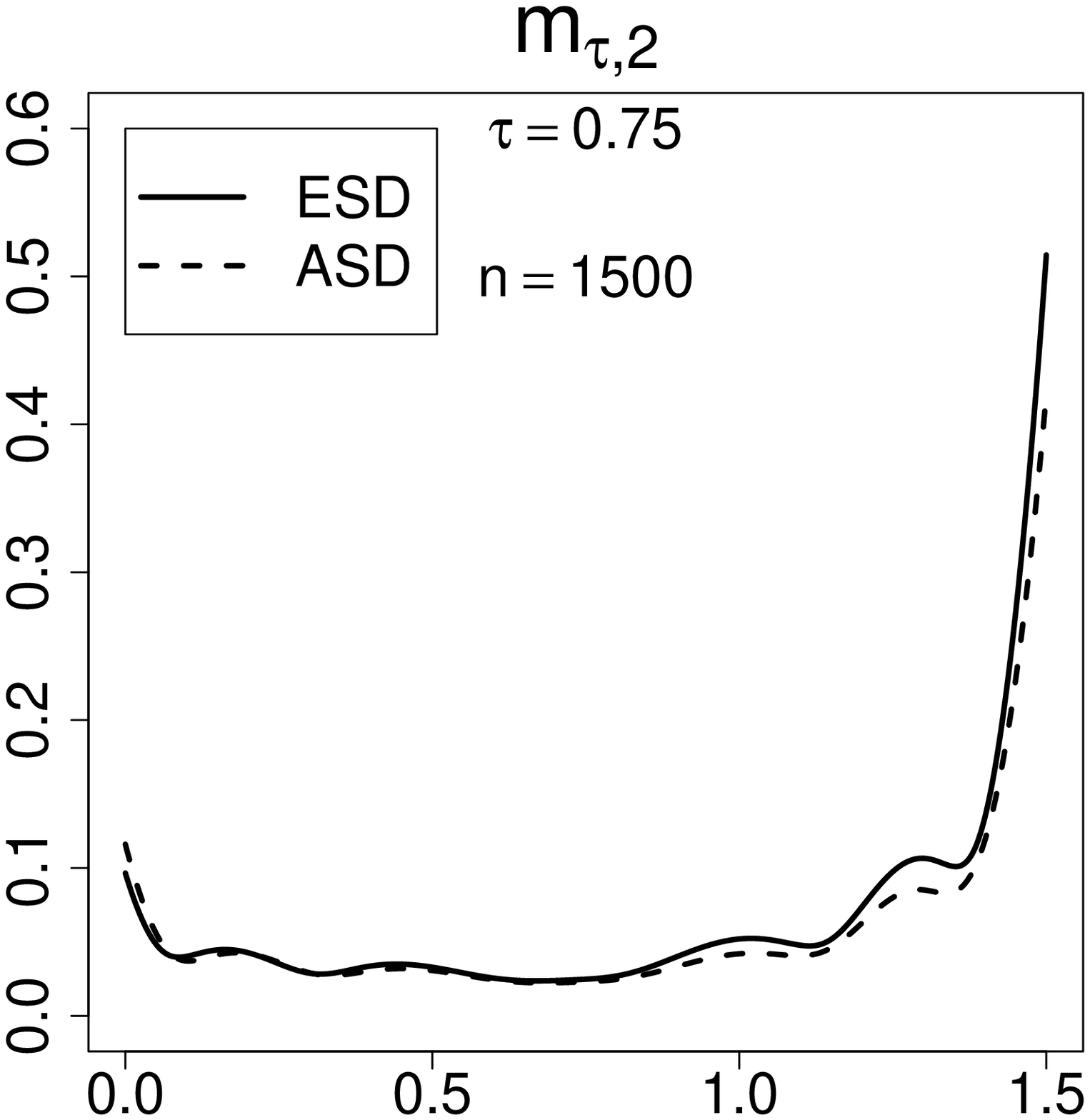}
\includegraphics[scale=0.2]{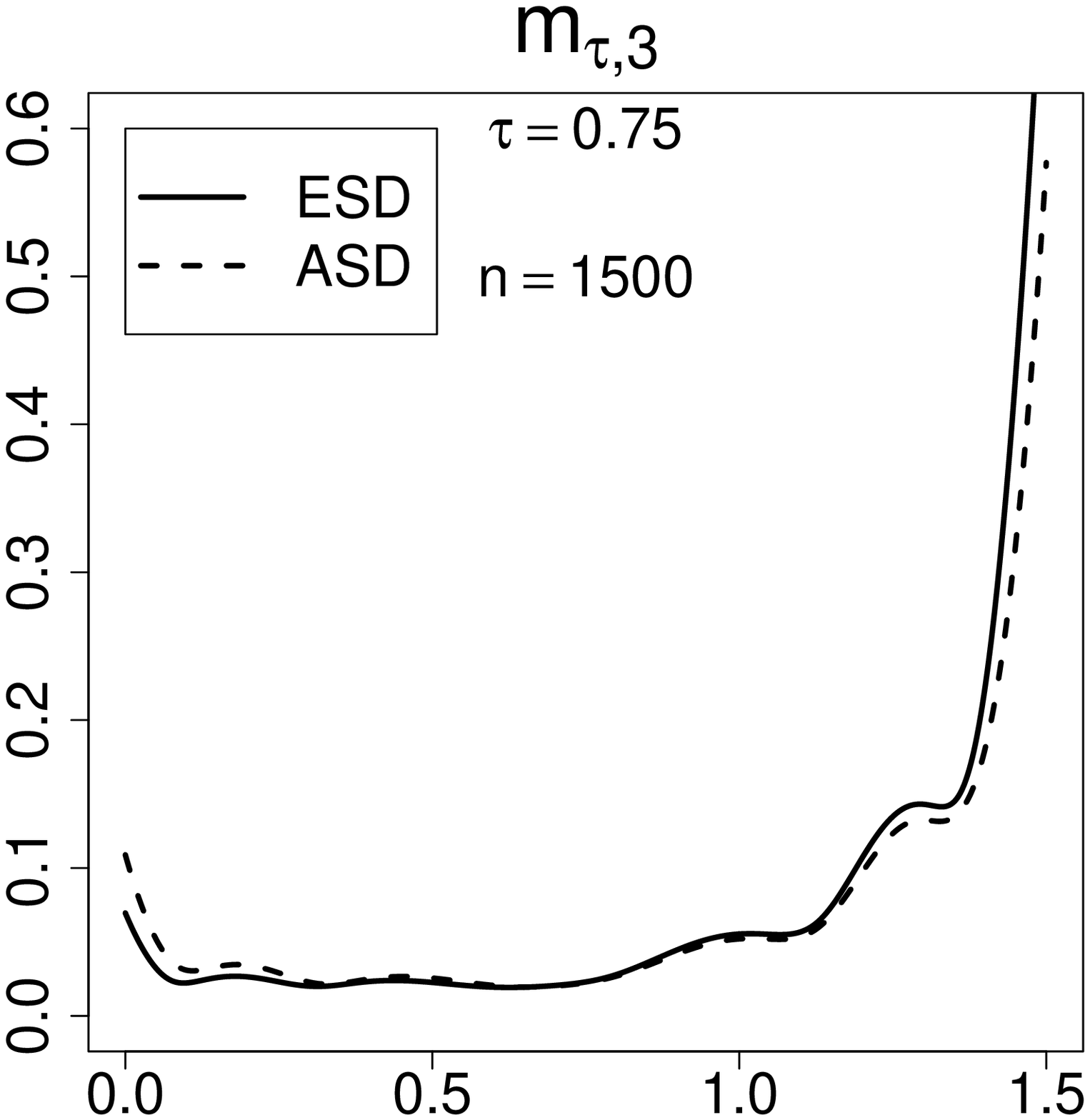}
\caption{The empirical standard deviation (ESD) and  estimated asymptotic standard deviation (ASD) for example 2 }
\label{figure2}
\end{figure}

\textsc{Example 2.} In this example, we specify the conditional quantile function $Q_{\tau}(Y_i |\bm X_i,\bm Z_i)$ to be
\[Q_{\tau}(Y_i |\bm X_i,\bm Z_i)= m_{\tau,1} ({\bm Z_i^T}{\bm \beta _{\tau,1}}){X_{i1}}+m_{\tau,2} ({\bm Z_i^T}{\bm \beta _{\tau,2}}){X_{i2}}+m_{\tau,3} ({\bm Z_i^T}{\bm \beta _{\tau,3}}){X_{i3}},\]
where ${\beta _{\tau ,1}} = \frac{{{{\left( {{\tau ^{1/2}},\tau ,2\tau } \right)}^T}}}{{\sqrt {5{\tau ^2} + \tau } }}$, ${\beta _{\tau ,2}}=\frac{{{{\left( {\tau ,{\tau ^{1/2}},2\tau } \right)}^T}}}{{\sqrt {5{\tau ^2} + \tau } }}$, ${\beta _{\tau ,3}}=\frac{{{{\left( {2\tau ,\tau ,{\tau ^{1/2}}} \right)}^T}}}{{\sqrt {5{\tau ^2} + \tau } }}$, ${m_{\tau ,1}}\left( {{u_1}} \right) = {\tau ^{1/2}}{u_1}$, ${m_{\tau ,2}}\left( {{u_2}} \right) = \tau \sin \left( {0.5\pi {u_2}} \right)$ and ${m_{\tau ,3}}\left( {{u_3}} \right) =  - 0.5\log \left( {1 - \tau } \right)u_3^2$. The covariate $X_{i1}=1$ and $\bm \left(X_{i2},X_{i3}\right)^T$ are generated from the independent standard normal distribution. The covariate $\bm Z_i=\left(Z_{i1},Z_{i2},Z_{i3}\right)^T$ are independently generated from Uniform [0,1]. Similar to \cite{FB16,MH16}, we generate $Y_i$ as
\[Y_i=m_{U_i,1} ({\bm Z_i^T}{\bm \beta _{U_i,1}}){X_{i1}}+m_{U_i,2} ({\bm Z_i^T}{\bm \beta _{U_i,2}}){X_{i2}}+m_{U_i,3} ({\bm Z_i^T}{\bm \beta _{U_i,3}}){X_{i3}}, \]
where $U_i$ follows a uniform distribution U(0,1). In this example, it is easy to see that the loading coefficients $\bm\beta_{\tau,l}$ and nonparametric functions $m_{\tau,l}$ for $l=1,2,3$ are functions of $\tau$, suggesting different covariate effects at different quantile levels. Thus, the VIC model structure is more sophisticated than that of example 1 and the mean regression method is no longer appropriate.

In this example, we consider estimation at the quartiles $\tau =0.5$ and $ \tau=0.75$, and simulate 500 data sets with $n = 500$ and $n=1500$. Tables \ref{table6} and \ref{table7} give the bias, ESD, ASD and MAD of $\bm\beta_{\tau,l}$, and RASE for $m_{\tau,l}$ for the proposed method, $l=1,2,3$. We may note that the true loading coefficients $\bm\beta_{\tau,l}$ and nonparametric functions $m_{\tau,l}$ are different at $\tau=0.5$ and $0.75$. The proposed estimation is also consistent with small biases, and the ESD, ASD, MAD and RASE become smaller with the increasing sample size. To evaluate the performance of nonparametric sandwich formula (\ref{eq8}), we define $\left\{0=t_1< t_2 <...< t_n=1.5\right\}$ as a grid set over the observed range of ${\bm Z}^T{\bm \beta}$. Figure \ref{figure2} depicts the ESD (calculated as the sample standard deviation of 500 estimates) and ASD (calculated by the nonparametric sandwich formula (\ref{eq8})) of $\hat m(t_j,\hat{\bm\beta})$ for $j=1,...,n$. It is evident that ESD and ASD are very similar and their difference decreases rapidly with the increasing sample size, indicating that the sandwich covariance formula (\ref{eq8}) performs well. This provides an assurance for the use of the nonparametric sandwich formula (\ref{eq8}) in practice.

\textsc{Example 3.} The main goal of this example is to investigate the finite sample performance of the proposed penalized estimation approach for identifying the linear components in quantile regression VICM. We generate random samples from model (\ref{eq13}) with $\sigma=0.2$, $d=4$, ${m_1}( u_1 )=0.2u_1^3$, ${m_2}( u_2 )=cos (0.5\pi {u_2})$, ${m_3}({u_3}) =0.5 u_3 $ and ${m_4}({u_4}) =-0.5 u_4 $. In this case we allow the last two nonparametric components to be linear functions. The true loading parameters are ${\bm\beta _1} = {\left( {\sqrt 2 /2,\sqrt 3 /3,\sqrt 6 /6,\bm 0_{p_n-3}} \right)^T}$, ${\bm\beta _2} = {\left( {\sqrt 3 /3,\sqrt 2 /2,\sqrt 6 /6,\bm 0_{p_n-3}} \right)^T}$, ${\bm\beta _3} = \frac{1}{{\sqrt {50} }}{\left( {3,4,5,\bm 0_{p_n-3}} \right)^T}$ and ${\bm\beta _4} = \frac{1}{{\sqrt {50} }}{\left( {4,3,5,\bm 0_{p_n-3}} \right)^T}$, where $\bm{0}_m$ denotes a $m$-vector of zeros. The dimension of $\bm\beta_l ( 1\leq l\leq 4)$ is set as $p_n=[n^{1/3}]$ for $n=500$ and $1500$. In this example, we focus on the quantile levels at $\tau=0.5$ and 0.75. To ensure $Q_{\tau}\left( Y | \bm X, \bm Z \right)=\sum\limits_{l = 1}^d {m_l\left(\bm Z_i^T\bm \beta_l\right)X_{il}}$ at $\tau=0.5$ and 0.75, we consider ${\epsilon _{\tau,i }} =\varsigma _{i}- {c_\tau }  $ and $c_\tau$ being the $\tau$th quantile of the random error $\varsigma _i$, resulting in $Q_{\tau}\left( \epsilon _{\tau,i }\right)=0$. Here $\{\varsigma _{i}\}$ is an i.i.d. random sample from SN, $t_3$, LA or MN. Other settings are the similar to that of example 1.

To evaluate the performance of variable selection and identification of linear components, we consider the following five criteria: (1) the average number of zero coefficients that are correctly estimated to be zero (C); (2) the average number of non zero coefficients that are incorrectly estimated to be zero (IC); (3) the average correctly fit percentage (CF) measures the accuracy of the variable selection procedure, where `` correctly fit" means that the procedure correctly select significant components from all $\bm\beta_l, l=1,2,3,4$; (4)  the proportion of $m_l$ being identified as the linear component for $l=1,2,3,4$ (${\rm{ILC}}_l$); (5) the proportion of correctly identification of linear components (CIL) among the four components. For the loading parameters, we compute the mean square error of the oracle estimators (O.MSE), the penalized estimators (P.MSE) and the unpenalized estimators (U.MSE). We also consider RASE of penalized estimators (P.RASE) and unpenalized estimators (U.RASE) that are used to measure the accuracy of nonparametric estimation. In each case, 500 data sets are generated. The simulation results are summarized in Tables \ref{table8}--\ref{table10}.

Eyeballing Tables \ref{table8} -- \ref{table10}, we can make several observations. Firstly, the values in the column labeled C are very close to the true number of zero loading parameters in Table \ref{table8}. The CF values steadily increase with the sample size $n$ and approach one quickly, which indicates that the proposed procedure is consistent in variable selection. Secondly, the proposed penalized estimator performs similarly as the oracle estimator in terms of estimation accuracy, and significantly reduces the MSE of the unpenalized estimator. Thirdly, we should realize that only the last two functions $m_3$ and $m_4$ are linear in this example. Thus, it is appealing to note that ${\rm{ILC}}_l$ is close to zero for $l=1,2$ and ${\rm{ILC}}_l$ approaches one for $l=3,4$ as the sample size increases. Table \ref{table9}  also shows that our penalized method can correctly distinguish linear components from nonparametric functions with a high probability. Fourthly, for the nonlinear functions ($m_1$ and $m_2$), there is a small difference for RASE between penalized and unpenalized estimators in Table \ref{table10}. However, our proposed penalized estimator is obviously more efficient leading to about 40$\%$-60$\%$ reduction in RASE for the linear components $m_3$ and $m_4$. The reason is that we apply a regularized estimation procedure to identify linear functions, namely, penalized method can discriminate the model structure. Therefore, comparing with the unpenalized estimators,  we see that the proposed penalized estimators improve the RASE. In summary, the proposed methods are satisfactory at different quantile levels in terms of variable selection and identification of linear components.

\tabcolsep=14pt
\begin{table}\scriptsize
\caption{Simulation results of RASE for $m_{\tau,1}$, $m_{\tau,2}$, $m_{\tau,3}$ with $\tau=0.5,0.75$ and $n = 500,1500$ in example 2.}
\label{table7}
\begin{tabular}{ccccccccccccccc} \noalign{\smallskip}\hline
\multirow{1}{*}{$n$}
&\multicolumn{3}{c}{$\tau=0.5$}
&&\multicolumn{3}{c}{$\tau=0.75$}
\\
\cline{2-4}
\cline{6-8}
&\multicolumn{1}{c}{$m_{\tau,1}$}&\multicolumn{1}{c}{$m_{\tau,2}$}&\multicolumn{1}{c}{$m_{\tau,3}$}
&&\multicolumn{1}{c}{$m_{\tau,1}$}&\multicolumn{1}{c}{$m_{\tau,2}$}&\multicolumn{1}{c}{$m_{\tau,3}$}
\\
\hline
500& 0.066&0.081&0.102&&0.098&0.189&0.210 \\
1500&0.041&0.043&0.059 &&0.070&0.170&0.196 \\
\hline
\end{tabular}
\end{table}

\tabcolsep=8pt
\begin{table}\scriptsize
\caption{{Simulations results of variable selection for $\bm\beta$ with $\tau=0.5$ and $\tau=0.75$ in example 3.}}
\label{table8}
\begin{tabular}{cccccccccc} \noalign{\smallskip}\hline
\multirow{1}{*}{$n$}&
\multirow{1}{*}{$p_n$}&
\multirow{1}{*}{$\tau$}&
\multirow{1}{*}{Error}&
\multicolumn{1}{c}{C}
&\multicolumn{1}{c}{IC}
&\multicolumn{1}{c}{CF}
&\multicolumn{1}{c}{O.MSE}
&\multicolumn{1}{c}{P.MSE}
&\multicolumn{1}{c}{U.MSE}\\
\hline
500&7&0.5&SN&15.87&  0.000&  0.928&0.890 &0.918& 2.370\\
&&&$t_3$&15.89 & 0.000&  0.944 &1.173& 1.214 & 3.263\\
&&&MN&15.91&  0.000&  0.938 &1.106& 1.149& 3.152\\
&&&LA&15.86  &0.000&  0.920& 1.088& 1.119& 3.081\\
&&0.75&SN&15.75&  0.008&  0.892& 0.938& 1.718& 3.418\\
&&&$t_3$&15.77&  0.002&  0.890& 1.514& 1.957& 4.479\\
&&&MN&15.81&  0.002 & 0.910&1.358& 1.564& 3.708\\
&&&LA& 15.74&  0.006&  0.876& 1.261& 1.806& 4.129\\
\hline
1500&11&0.5&SN&31.97 & 0.000&  0.982 &0.175& 0.176& 0.781\\
&&&$t_3$&31.99&  0.002&  0.988& 0.239 &0.240& 1.024\\
&&&MN&31.99& 0.000&  0.992& 0.222& 0.222& 1.033\\
&&&LA&31.97&  0.000  &0.980&  0.212 &0.213 &1.002\\
&&0.75&SN&31.95&  0.002&  0.966& 0.205 &0.205 & 0.936\\
&&&$t_3$&31.87&  0.006&  0.968 &0.302& 0.305& 1.437\\
&&&MN&31.97&  0.002&  0.986&  0.272 &0.272 &1.256\\
&&&LA&31.94&  0.006&  0.962&0.284 &0.284 &1.306\\
\hline
\end{tabular}
\\
Notation: the values of last three columns multiplied by $10^{-2}$ are true simulation results of O.MSE, P.MSE and U.MSE. In addition the number of zero coefficients is 16 for $n=500$ and 32 for $n=1500$.
\end{table}

\tabcolsep=5pt
\begin{table}\scriptsize
\caption{Simulations results of linear component identification for $m_l,l=1,2,3,4$ with $\tau=0.5$ and $\tau=0.75$ in example 3. }
\label{table9}
\begin{tabular}{ccccccccccccccc} \noalign{\smallskip}\hline
\multicolumn{1}{c}{$n$}&
\multirow{1}{*}{Error}
&\multicolumn{5}{c}{$\tau=0.5$}
&&\multicolumn{5}{c}{$\tau=0.75$}
\\
\cline{3-7}
\cline{9-13}
&&\multicolumn{1}{c}{${\rm{ILC}}_1$}&\multicolumn{1}{c}{${\rm{ILC}}_2$}
&\multicolumn{1}{c}{${\rm{ILC}}_3$}&\multicolumn{1}{c}{${\rm{ILC}}_4$}&\multicolumn{1}{c}{CIL}
&&\multicolumn{1}{c}{${\rm{ILC}}_1$}&\multicolumn{1}{c}{${\rm{ILC}}_1$}
&\multicolumn{1}{c}{${\rm{ILC}}_3$}&\multicolumn{1}{c}{${\rm{ILC}}_4$}&\multicolumn{1}{c}{CIL}
\\
\hline
500&SN&0.000& 0.000& 0.792& 0.812 &0.704&&0.000 &0.002& 0.702& 0.734& 0.620\\
&$t_3$&0.000& 0.000& 0.844& 0.842 &0.754&&0.000& 0.000& 0.774 &0.792 &0.690\\
&MN&0.000& 0.000& 0.794& 0.868& 0.734&&0.000& 0.000& 0.786& 0.804& 0.686\\
&LA&0.000& 0.000& 0.822 &0.846& 0.730&&0.000& 0.000 &0.790 &0.804& 0.696\\
\hline
1500&SN&0.000& 0.000& 0.922& 0.952& 0.892&&0.000 &0.002& 0.910 &0.924& 0.862\\
&$t_3$&0.000 &0.000& 0.962& 0.962& 0.938&&0.000& 0.000& 0.946& 0.942& 0.898\\
&MN&0.000& 0.000& 0.968 &0.986 &0.960&&0.000 &0.002& 0.942& 0.948& 0.916\\
&LA&0.000& 0.000& 0.960& 0.970& 0.938&&0.000& 0.000& 0.938 &0.952& 0.904\\
\hline
\end{tabular}
\end{table}

\tabcolsep=6pt
\begin{table}\scriptsize
\caption{Simulation results of RASE for $m_l,l=1,2,3,4$ with $\tau=0.5$ and $\tau=0.75$ in example 3.}
\label{table10}
\begin{tabular}{ccccccccccccccc} \noalign{\smallskip}\hline
\multicolumn{1}{c}{$n$}&
\multicolumn{1}{c}{$\tau$}&
\multirow{1}{*}{Error}&
\multicolumn{4}{c}{P.RASE}
&&\multicolumn{4}{c}{U.RASE}
\\
\cline{4-7}
\cline{9-12}
&&&\multicolumn{1}{c}{$m_1$}&\multicolumn{1}{c}{$m_2$}&\multicolumn{1}{c}{$m_3$}
&\multicolumn{1}{c}{$m_4$}
&&\multicolumn{1}{c}{$m_1$}&\multicolumn{1}{c}{$m_2$}&\multicolumn{1}{c}{$m_3$}
&\multicolumn{1}{c}{$m_4$}
\\
\hline
500&0.5&SN&0.125& 0.516 &0.104 &0.091&& 0.129& 0.530& 0.184 &0.162\\
&&$t_3$&0.145& 0.515 &0.105 &0.087 &&0.146 &0.533 & 0.195& 0.165\\
&&MN&0.132& 0.494 &0.115 &0.082&& 0.138& 0.507& 0.203 &0.163\\
&&LA&0.125& 0.490 &0.093 &0.080 &&0.130& 0.507 &0.179& 0.152\\
&0.75&SN&0.165& 0.525& 0.128& 0.107&& 0.165& 0.540 &0.196& 0.174\\
&&$t_3$&0.162& 0.498 &0.120& 0.109&& 0.161& 0.522& 0.205& 0.189\\
&&MN&0.170& 0.492& 0.108& 0.101&& 0.166 &0.515& 0.202 &0.176\\
&&LA&0.167& 0.493& 0.107 &0.095 &&0.164& 0.509& 0.187& 0.166\\
\hline
1500&0.5&SN&0.048& 0.217& 0.028 &0.021&& 0.051& 0.224& 0.064 &0.057\\
&&$t_3$&0.050& 0.210 &0.027 &0.022 &&0.051& 0.223 & 0.071& 0.061\\
&&MN&0.050& 0.210 &0.023 &0.019 &&0.050& 0.222& 0.069 &0.058\\
&&LA&0.047 &0.217 &0.023 &0.020&& 0.049& 0.228& 0.066& 0.060\\
&0.75&SN&0.064& 0.216 &0.030& 0.025&& 0.064& 0.226 & 0.074& 0.062\\
&&$t_3$&0.073& 0.207& 0.032 &0.029&& 0.071 &0.220& 0.078 &0.070\\
&&MN&0.073& 0.208 &0.030& 0.027 &&0.072& 0.223& 0.079 &0.068\\
&&LA&0.076& 0.205&0.029 &0.025&& 0.076& 0.218& 0.078 &0.068\\
\hline
\end{tabular}
\end{table}

\begin{figure}
\includegraphics[scale=0.25]{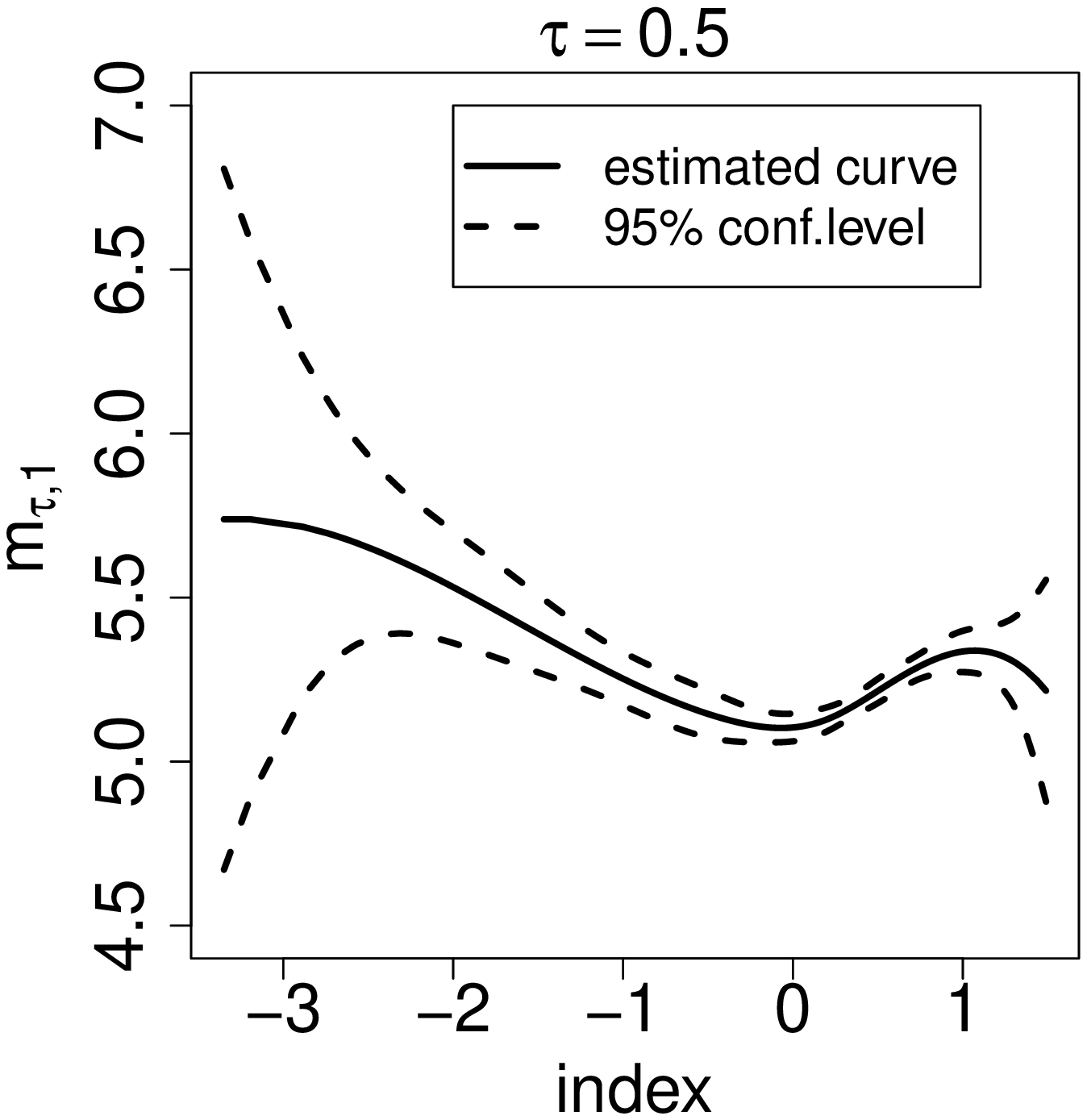}
\includegraphics[scale=0.25]{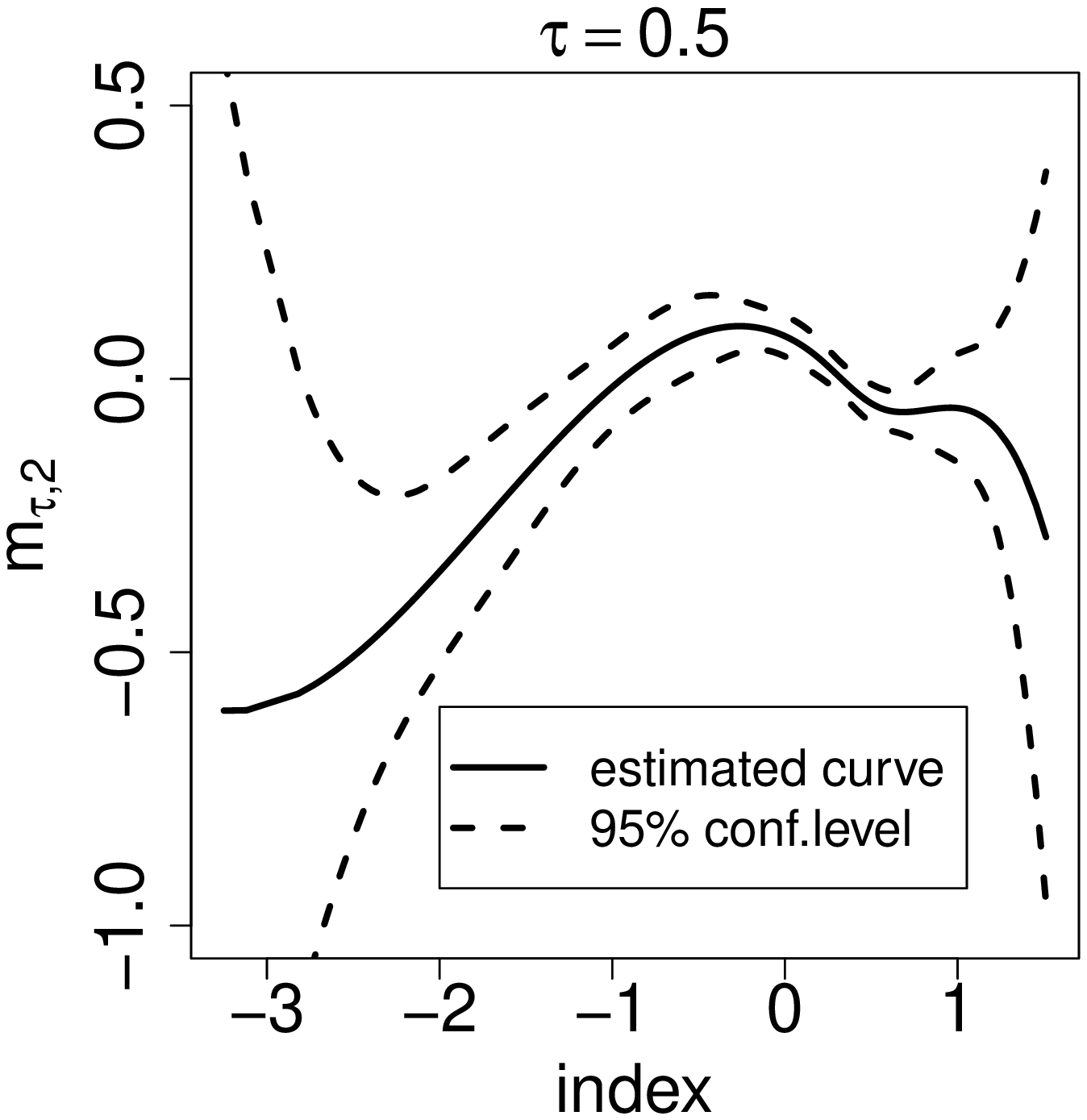}
\includegraphics[scale=0.25]{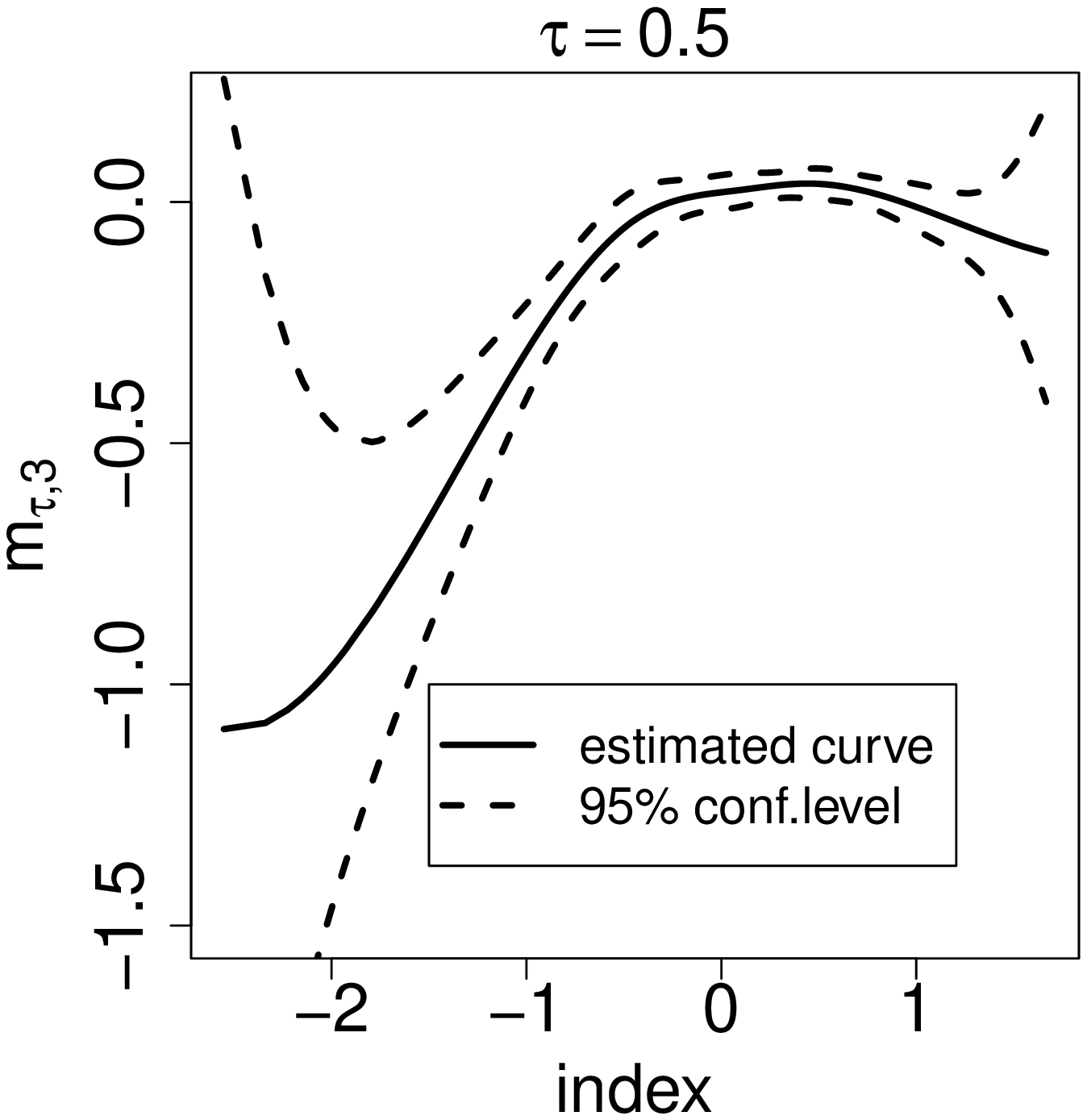}
\includegraphics[scale=0.25]{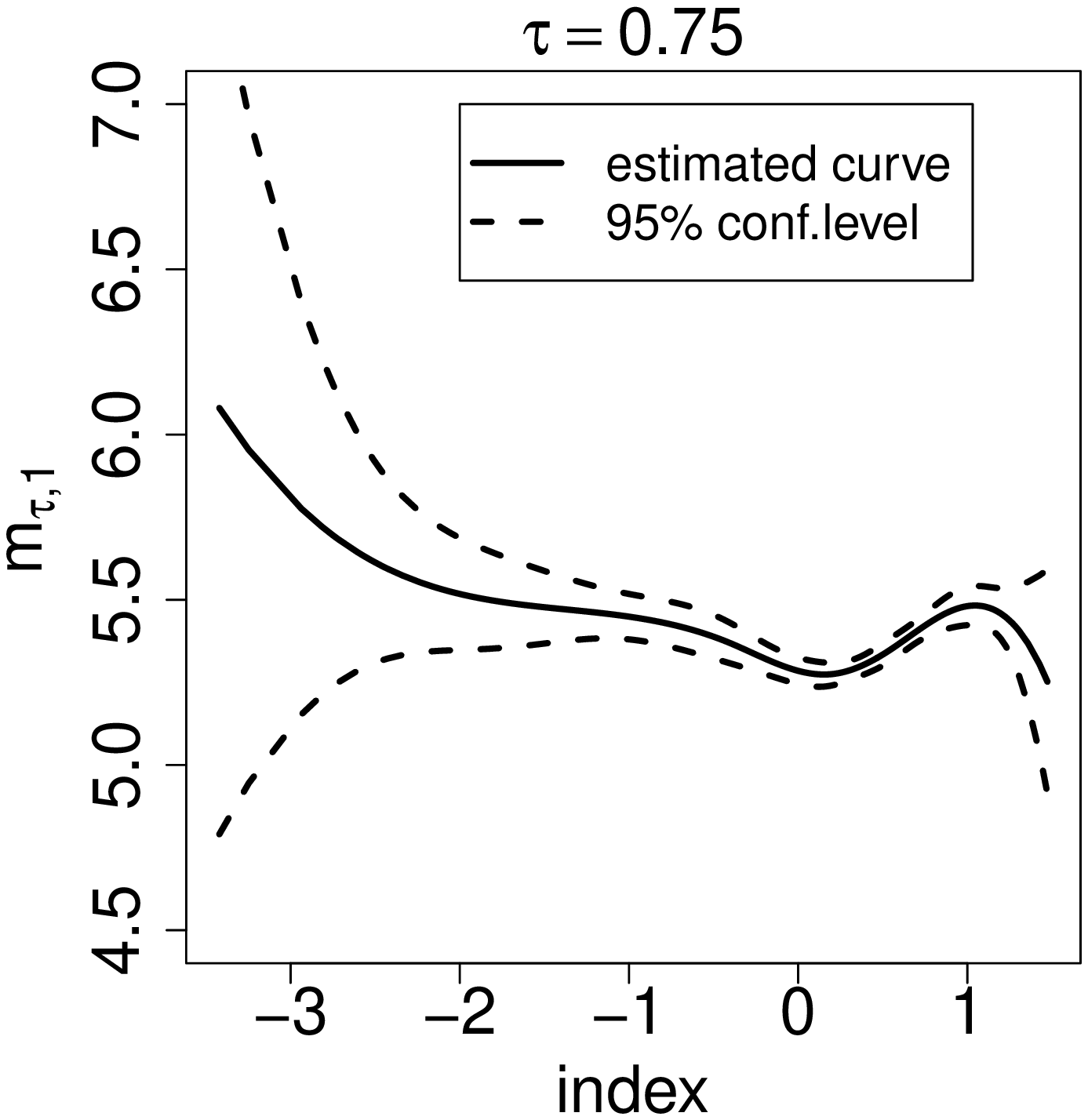}
\includegraphics[scale=0.25]{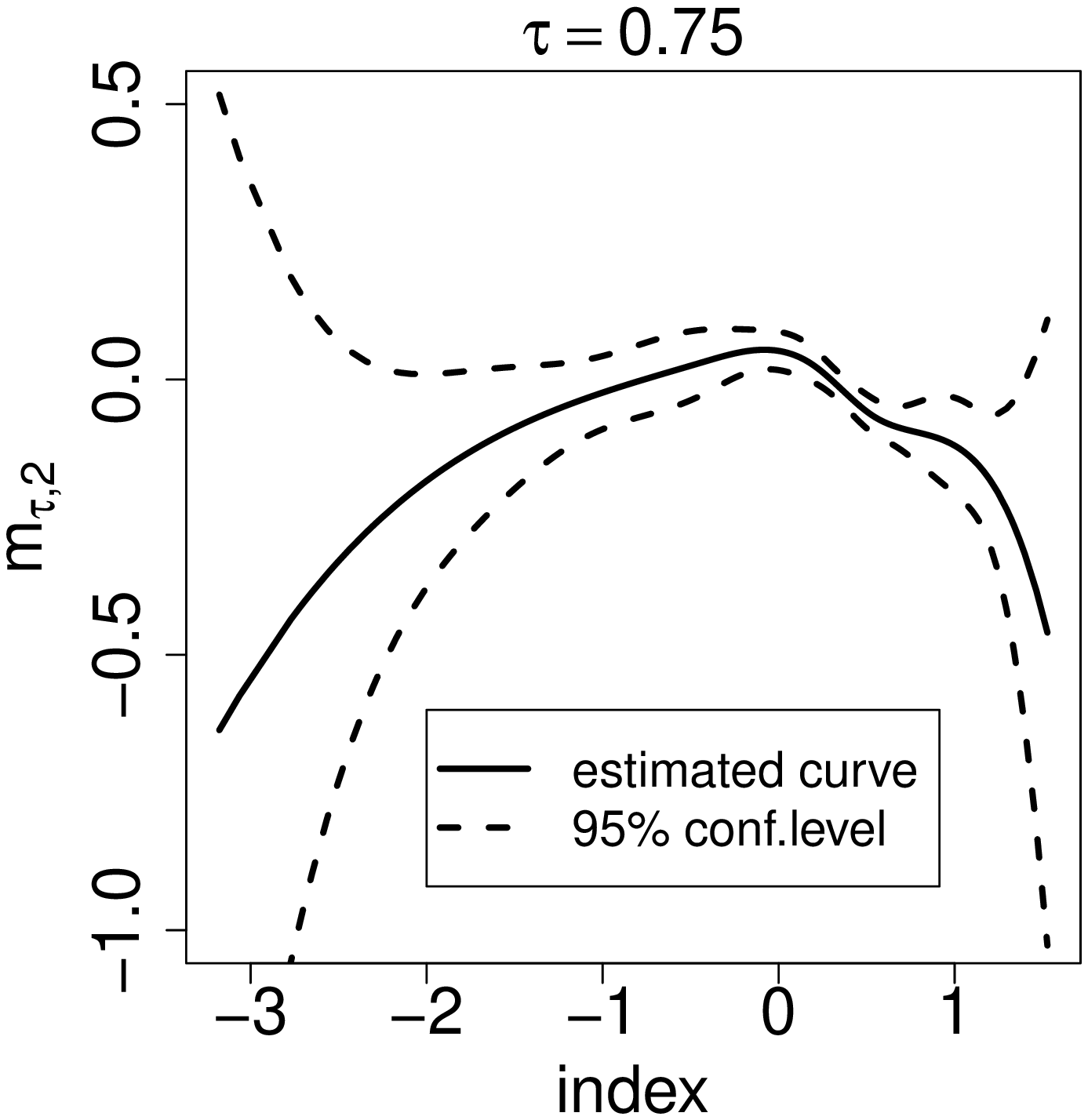}
\includegraphics[scale=0.25]{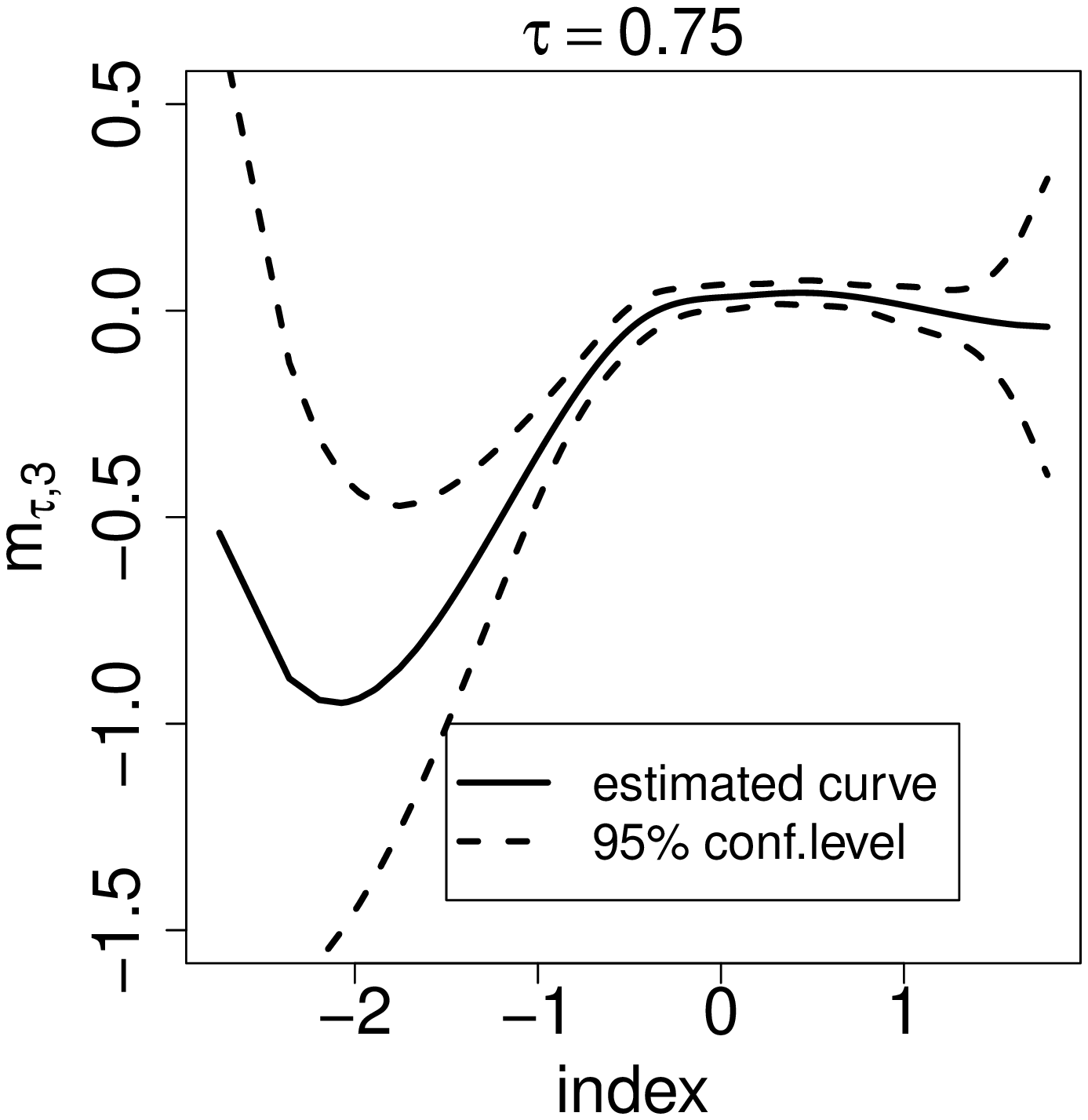}
\caption{ Estimated nonparametric curves for $\tau = 0.5$ and $ 0.75$ and their $95\%$ confidence intervals.}
\label{figure3}
\end{figure}

\subsection{Real data analysis}\label{subsect.5.3}

In this subsection, we will illustrate the proposed approaches by analyzing an environmental dataset in Hong Kong (\cite{FZ99,HZ13}). This dataset consists of a collection of daily measurements of pollutants and two environmental factors, with a total of  $n=1461$ observations. Three pollutants, nitrogen dioxide ($no_2$ ($\mu g/m^3$)), sulphur dioxide ($so_2$ ($\mu g/m^3$)) and ozone ($o_3$ ($\mu g/m^3$)) and two weather elements, temperature ($temp$ ($^{\circ} C$)) and relative humidity ($hum$ ($\%$)), are considered here. These five factors may give rise to circulatory and respiratory problems of humans. In this study, our goal is to explore whether these five variables may influence the number of daily total hospital admissions. Thus, we study the relationship between the number of daily hospital admissions ($Y$) and the following covariates: $no_2$ ($Z_1$), $so_2$ ($Z_2$), $o_3$ ($Z_3$), $tem$ ($X_2$) and $hum$ ($X_3$). Here we take $X_1=1$ as the intercept term. Specifically, we use model (\ref{eq1}) to fit the data, where $\bm{Z_i} = {\left( {Z_{i1},Z_{i2},Z_{i3}} \right)^T}$ is a covariate vector with length $p=3$, $m_l(\cdot)$ are the unknown smooth functions and $\bm{\beta }_l= {\left( {{\beta _{l1}},{\beta _{l2}},{\beta _{l3}}} \right)^T}$ are unknown loading parameters for $l=1,2,3$. Before implementing the estimation procedure, we normalize all predictor variables and take the logarithm of the response variable. The initial estimates of the parameters are obtained by the profile least squares method proposed in \cite{MS15}. In our analysis, tuning parameters (e.g., $\alpha_1$ and $\alpha_2$) are chosen based on the information criteria given in subsection \ref{subsect.5.1}, and the bandwidth is set as $h=cn^{-0.2}$ with $c=0.1,0.2,...,1$, the optimal bandwidth is selected by the 5-fold cross-validation. In this real data analysis, we consider the unpenalized estimators ($\hat{\bm\beta}_l$ and $\hat{m}_l$) and the penalized estimator ($\bar{\bm\beta}_l$ and $\bar{m}_l$) at two quantile levels $\tau = 0.5, 0.75$.

\tabcolsep=10pt
\begin{table}\scriptsize
\renewcommand\arraystretch{1}
\caption{The estimates (EST),  estimated asymptotic standard deviation (ASD) of $\bm\beta_l$, and $p$-values for testing significance of each component in $\bm\beta_l$ for $l=1,2,3$ in environmental data.}
\label{table11}
\begin{tabular}{cccccccccc} \noalign{\smallskip}\hline
\multicolumn{1}{c}{$\tau$}
&&&\multicolumn{3}{c}{unpenalized}
&&\multicolumn{2}{c}{penalized}
\\
\cline{4-6}
\cline{8-10}
&&&\multicolumn{1}{c}{EST}&\multicolumn{1}{c}{ASD}&\multicolumn{1}{c}{p-value}
&&\multicolumn{1}{c}{EST}&\multicolumn{1}{c}{ASD}&
\\
\hline
\multicolumn{10}{c}{$X_1$=intercept}\\
0.5&$\bm\beta_1$&$Z_1$&0.537& 0.039& $<10^{-3}$&&0.554 &0.040\\
&&$Z_2$&-0.843 &0.026& $<10^{-3}$&&-0.832 &0.026\\
&&$Z_3$&-0.041 &0.051& 0.427&&0&0\\

\multicolumn{10}{c}{$X_2$=tem}\\
&$\bm\beta_2$&$Z_1$&0.604& 0.043& $<10^{-3}$&&0.573& 0.040\\
&&$Z_2$& -0.788& 0.037&$<10^{-3}$&&-0.820& 0.028\\
&&$Z_3$&-0.123& 0.064 &0.056&&0&0\\

\multicolumn{10}{c}{$X_3$=hum}\\
&$\bm\beta_3$&$Z_1$&0.773& 0.042& $<10^{-3}$&&0.787 &0.041\\
&&$Z_2$&-0.601 &0.060& $<10^{-3}$&&-0.585 &0.063\\
&&$Z_3$&-0.202 &0.083 &0.015&&-0.197& 0.092\\
\hline
\multicolumn{10}{c}{$X_1$=intercept}\\
0.75&$\bm\beta_1$&$Z_1$&0.543& 0.036 &$<10^{-3}$&&0.541& 0.040\\
&&$Z_2$&-0.840& 0.024& $<10^{-3}$&&-0.841& 0.026
\\
&&$Z_3$&-0.020& 0.054& 0.712&&0&0\\

\multicolumn{10}{c}{$X_2$=tem}\\
&$\bm\beta_2$&$Z_1$&0.653 &0.030& $<10^{-3}$&&0.586 &0.043\\
&&$Z_2$&-0.711 &0.038 &$<10^{-3}$&&-0.811& 0.031\\
&&$Z_3$&-0.261 &0.069& $<10^{-3}$&&0&0\\

\multicolumn{10}{c}{$X_3$=hum}\\
&$\bm\beta_3$&$Z_2$&0.775& 0.027& $<10^{-3}$&&0.791 &0.032\\
&&$Z_2$&-0.558 &0.051 &$<10^{-3}$&&-0.547& 0.054\\
&&$Z_3$&-0.298& 0.066& $<10^{-3}$&&-0.275& 0.065\\
\hline
\end{tabular}
\end{table}

\begin{figure}\centering
\includegraphics[scale=0.32]{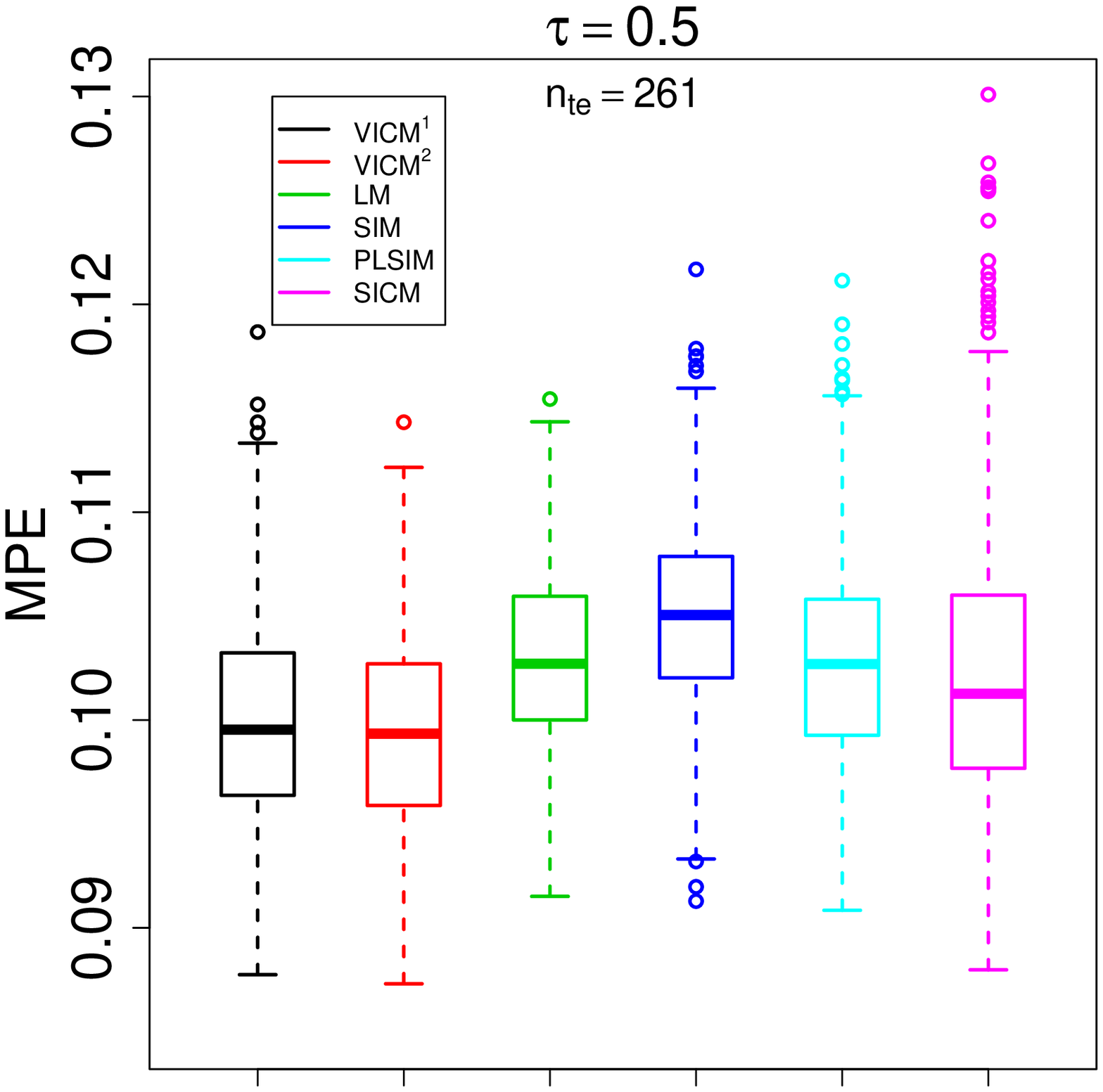}
\includegraphics[scale=0.32]{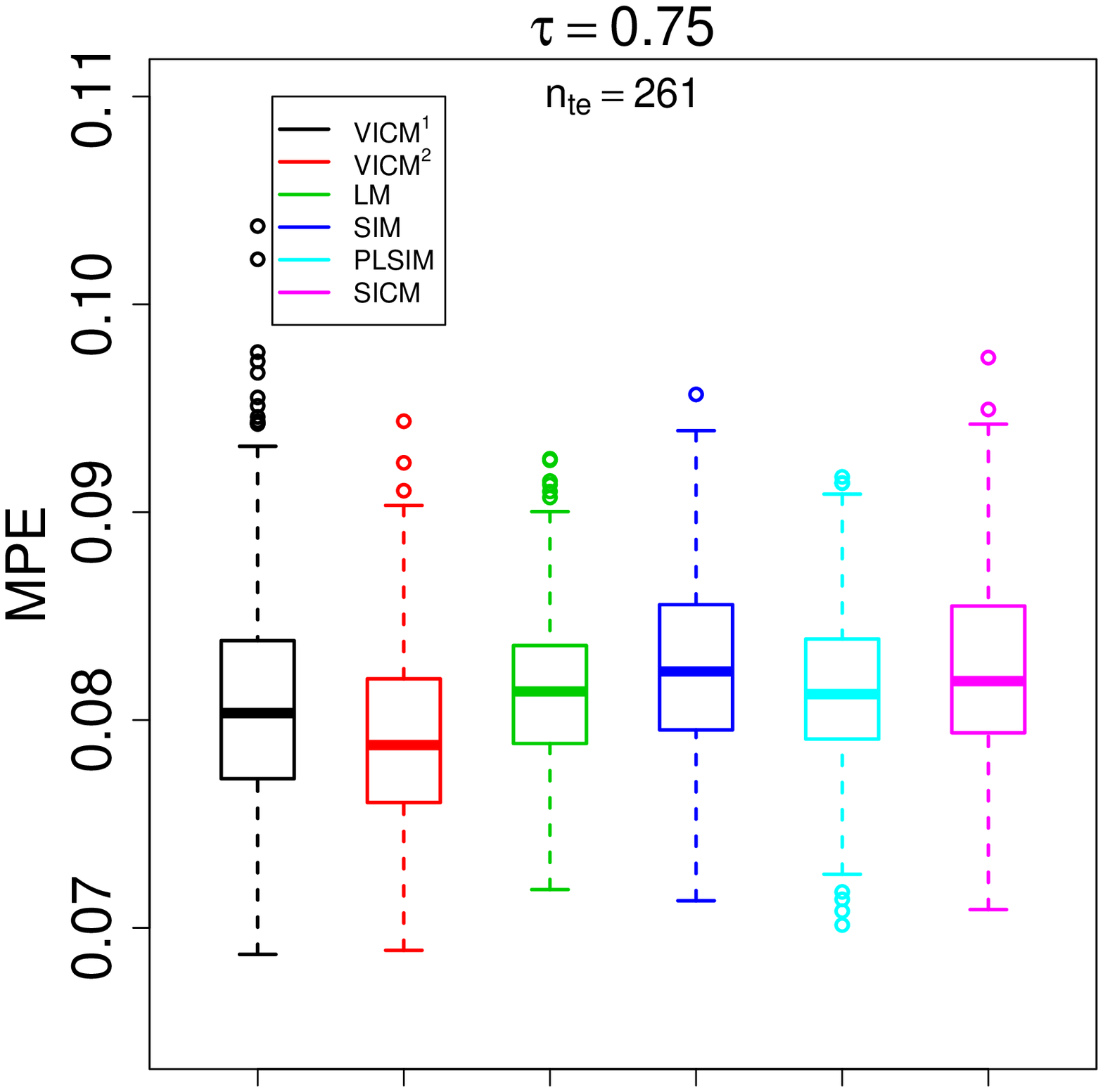}
\includegraphics[scale=0.32]{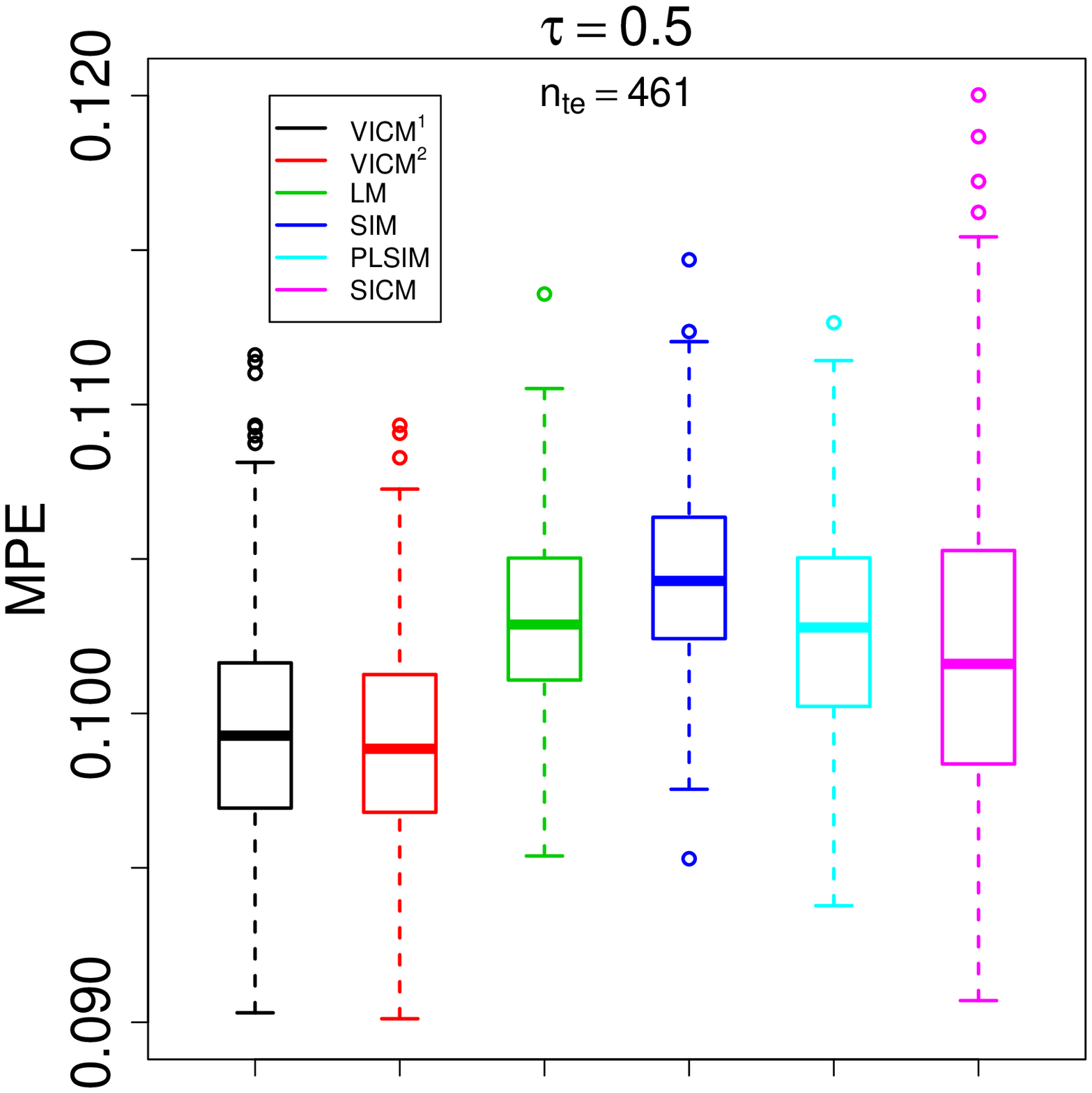}
\includegraphics[scale=0.32]{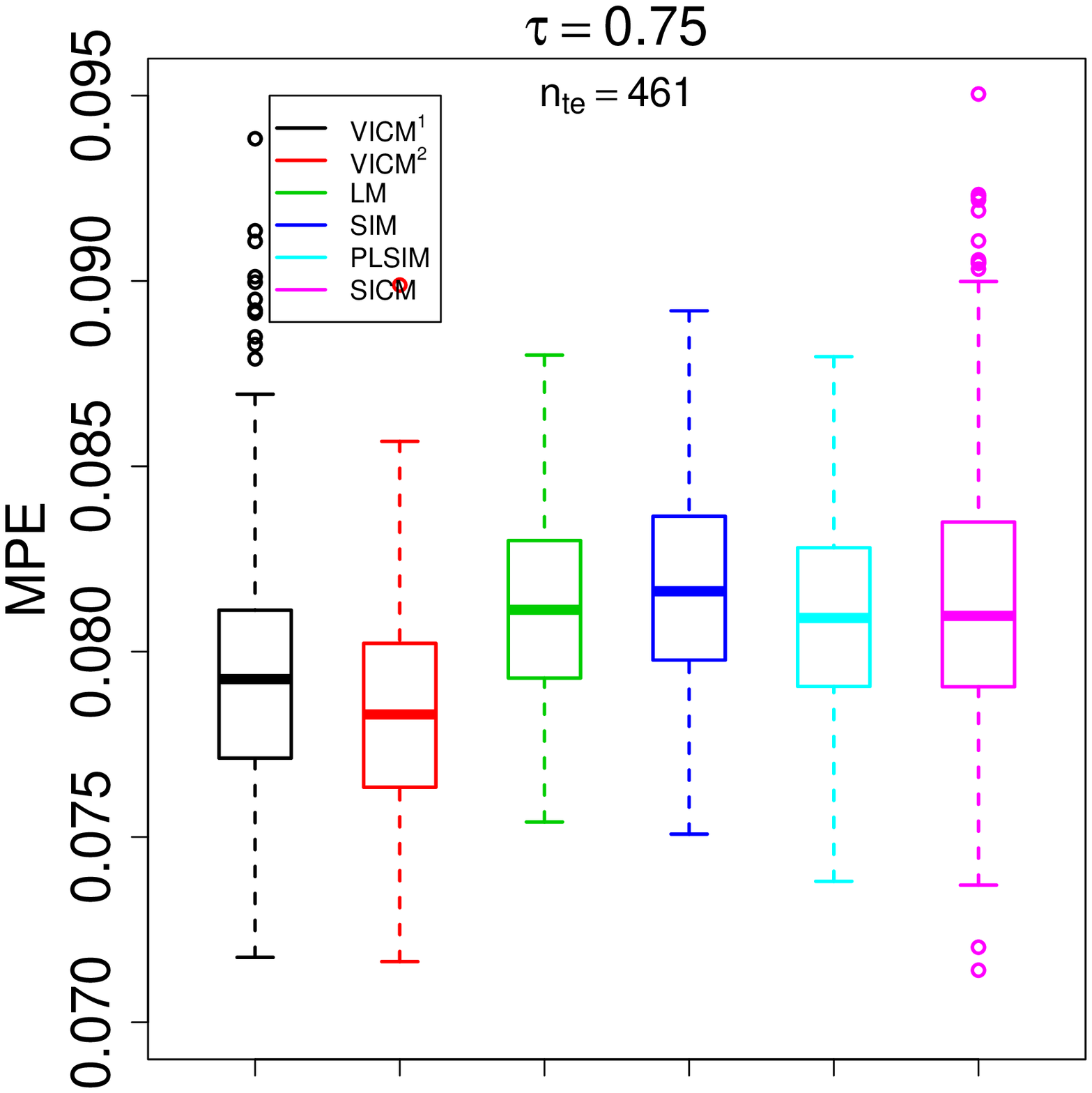}
\caption{ The mean prediction error (MPE) of VICM$^1$, VICM$^2$, LM, SIM, PLSIM and SICM at $\tau =0.5$ and $\tau=0.75$ for $n_{te}=261$ and 461.}
\label{figure4}
\end{figure}

Table \ref{table11} shows the estimated coefficients (EST), their estimated asymptotic standard deviation (ASD) calculated by the sandwich formula (\ref{eq11}), as well as the $p$-values for testing significance of each pollutant. We first notice that the loading parameters for $Z_1$ and $Z_2$ are highly significant for the intercept term $X_1$, suggesting  $no_2$ and $so_2$ are strong predictors for the daily hospital admission numbers at both quantile levels. For the temperature variable $X_2$, we observe that the two loading parameters corresponding to $Z_1$ and $Z_2$ are significantly different from zero at significance level 0.05. This implies that $no_2$ and $so_2$ have significant nonlinear interaction effects with temperature on the number of daily hospital admissions. However, for relative humidity $X_3$, all predictors $Z_1$, $Z_2$ and $Z_3$ are significantly different from zero at significance level 0.05, which indicates that $no_2$, $so_2$ and $o_3$ have significant nonlinear interaction effects with relative humidity. Table \ref{table11} also shows the estimated coefficients and their standard errors for the penalized estimators, the loadings for $Z_3$ are penalized to be zero for the intercept term $X_1$ and temperature $X_2$, but not for relative humidity $X_3$. By using the penalized estimate $\bar{\bm\lambda}_l$ (l=1,2,3) given in Sect. \ref{sect.4}, we can obtain $\| \bar{\bm\lambda}_1 \|_{\bm D}=19.34$, $\| \bar{\bm\lambda}_2 \|_{\bm D}=24.84$ and $\| \bar{\bm\lambda}_3 \|_{\bm D}=16.42$ for $\tau=0.5$ and $\| \bar{\bm\lambda}_1 \|_{\bm D}=28.11$, $\| \bar{\bm\lambda}_2 \|_{\bm D}=22.53$ and $\| \bar{\bm\lambda}_3 \|_{\bm D}=20.34$ for $\tau=0.75$, which are far away from zero. Thus, $m_1$, $m_2$ and $m_3$ are identified as nonlinear functions. Fig. \ref{figure3} displays the estimated curves by the penalized method and their 95$\%$ confidence bands obtained by nonparametric sandwich formula (\ref{eq8}). The plot for the intercept shows that the estimated function $\bar m_1(.)$ is a decreasing function of index $\bm Z^T\bm{\bar \beta}_{\alpha_1,1}$, which indicates that the combination of environmental factors has a negative effect on the daily hospital admission numbers.
The plots of temperature and relative humidity demonstrate that the effects of temperature and relative humidity are nonlinearly modified. These finding are consistent with the penalized estimation results.

Next we consider the prediction performance of the proposed method. To this end, the data is randomly divided into two parts. The first part is reserved as a training dataset including $n_{tr}$ observations while the second part is reserved as a test dataset including $n_{te}$ observations, where $n=n_{tr}+n_{te}$. Here we consider $n_{te}=261$ and 461. We compute the mean prediction error (MPE) to evaluate the prediction performance for different quantile regression models. The MPE is defined as ${\rm{MPE}} =\sum\nolimits_{i \in \mathcal{I}} {{\rho _\tau }\left( {{Y_i} - {{\hat Y}_i}} \right)} /\left| \mathcal{I} \right|$, where $\mathcal{I}$ stands for an index set of the testing sample. For MPE, we repeat the random splitting procedure for 500 times and report the average. We compare the following six models and use quantile regression method to fit every model with $\tau=0.5$ and 0.75: the varying index coefficient model (VICM$^1$), the linear model (LM), the single-index model (\cite{XTLZ02}; SIM), the partially linear single index model (\cite{WXZC10}; PLSIM), the single-index coefficient model (\cite{HZ13}; SICM) and the penalized varying index coefficient model (VICM$^2$). The MPEs of six models are displayed in Fig. \ref{figure4}, clearly showing a superior performance of our proposed VICM to predict the future response values. In addition, we apply the mean quantile residual (MQR) to evaluate the in-sample performance of different models, which is defined as ${\rm{MQR}} ={1461^{ - 1}}\sum\limits_{i = 1}^{1461} {{\rho _\tau }\left( {{Y_i} - {{\hat Y}_i}} \right)}$. The MQRs of VICM$^1$, VICM$^2$, LM, SIM, PLSIM, SICM are 0.0939, 0.0938, 0.1023, 0.1025, 0.1024, 0.0993 for $\tau=0.5$ and 0.0743, 0.0738, 0.0802, 0.0792, 0.0790, 0.0771 for $\tau=0.75$, indicating that our VICM achieves the smallest in-sample error.

\section{Concluding remarks}\label{sect.6}
In this paper, we apply the SCAD penalty to develop robust variable selection and linear components identification procedures for the VICM under the quantile regression framework, in which the nonparametric functions are approximated by B-spline basis functions. While other types of basis functions are applicable, B-spline functions are relatively easy to implement with practically stable performance. We establish the consistency and oracle property of the estimators in the situation of a slowly diverging number of loading parameters. In addition, we develop a novel penalization method to distinguish linear components automatically. To reduce the computational burden caused by the non-smoothing estimating equations, we utilize a kernel function to approximate the quantile score function, which results in smoothing estimating equations and facilitates a sandwich formula for variance estimation. Some useful criteria are proposed to choose the tuning parameters, and simulation studies and real data analysis have been conducted to illustrate the proposed method and confirm the asymptotic results. Finally, it is interesting to study high dimensional variable selection and model identification for the VICM with complex data including longitudinal data, multi-level data, censored survival data and others. Research in these aspects is ongoing.

\section*{Acknowledgements}
Jing Lv is partially supported by National Natural Science Foundation of China Grant 11801466 and the Basic and Frontier Research Program of Chongqing Grant cstc2017jcyjAX0182. Jialiang Li is partially supported by Academic Research Funds R-155-000-174-114, R-155-000-195-114 and Tier 2 Ministry of Education funds in Singapore MOE2017-T2-2-082: R-155-000-197-112 (Direct cost) and R-155-000-197-113 (IRC).


\begin{supplement}
\sname{Supplement }\label{supp}
\stitle{``High-dimensional varying index coefficient quantile regression model''}
\slink[url]{Supplemental materials.zip}
\sdescription{The supplementary materials list additional numerical results and regular conditions, and provide a number of technical lemmas and the proofs of lemmas and theorems.}
\end{supplement}


\begin{thebibliography}{00}
\bibitem{A2001} \textsc{Abrevaya, J.} (2001). The effect of demographics and maternal behavior on the distribution of birth outcomes. \textit{Empirical Economics} \textbf{26} 247--259.

 \bibitem[Brown and Wang, 2007]{BW07} \textsc{Brown, B.M.} and \textsc{Wang, Y.} (2007). Induced smoothing for rank regression with censored survival times. \textit{Statistics in Medicine} \textbf{26} 828--836
\bibitem[B\"{u}hlmann and van de Geer, 2011]{BV11} \textsc{B\"{u}hlmann, P.} and \textsc{van de Geer, S.} (2011). \textit{Statistics for high-dimensional data: methods, theory and applications}. Springer Series in Statistics. Heidelberg: Springer.
\bibitem[Chen and Chen, 2008]{CC08} \textsc{Chen, J.} and \textsc{Chen, Z.} (2008). Extended Bayesian information criteria for model selection with large model spaces. \textit{Biometrika} \textbf{95} 759--771.
 \bibitem[Chiou, Kang and Yan, 2015]{CKY15} \textsc{Chiou, S.} \textsc{Kang, S.} and \textsc{Yan, J.} (2015). Semiparametric Accelerated failure time modeling for clustered failure times from stratified sampling. \textit{Journal of the American Statistical Association} \textbf{110} 621--629.
\bibitem[Christou and Akritas, 2016]{CA16} \textsc{Christou, E.} and \textsc{Akritas, M. G.} (2016). Single index quantile regression for heteroscedastic data. \textit{Journal of Multivariate Analysis} \textbf{150} 169--182.
 \bibitem[Cui, H\"{a}rdle and Zhu, 2011]{CHZ11} \textsc{Cui, X.} \textsc{H\"{a}rdle, W. K.} and \textsc{Zhu, L.} (2011). The EFM approach for single-index models. \textit{The Annals of Statistics} \textbf{39} 1658--1688.
  \bibitem[de Boor, 2001]{De01} \textsc{de Boor, C.} (2001). \textit{A practical guide to splines}. Springer, New York.
\bibitem{EKJ2008} \textsc{Elsner, J. B.} \textsc{Kossin, J. P.} and \textsc{Jagger, T. H.} (2008). The increasing intensity of the strongest tropical cyclones. \textit{Nature} \textbf{455} 92--95.
 \bibitem[Fan and Li, 2006]{FL06} \textsc{Fan, J.} and \textsc{Li, R.} (2006). \textit{Statistical challenges with high dimensionality: feature selection in knowledge discovery}. Proceedings of the Madrid International Congress of Mathematicians, III: 595--622.
 \bibitem[Fan, Liu and Lu, 2017]{FLL17}\textsc{Fan, J.} \textsc{Liu, W.} and \textsc{ Lu, X.} (2017). Penalized empirical likelihood for semiparametric models with a diverging number of parameters. \textit{Journal of Statistical Planning and Inference} \textbf{186} 42-57.
 \bibitem[Fan and Peng, 2004]{FP04} \textsc{Fan, J.} and \textsc{Peng, H.} (2004). Nonconcave penalized likelihood with a diverging number of parameters. \textit{The Annals of Statistics} \textbf{32} 928-961.
 \bibitem[Fan and Zhang, 1999]{FZ99} \textsc{Fan, J.} and \textsc{Zhang, W.} (1999). Statistical estimation in varying coefficient models. \textit{The Annals of Statistics} \textbf{27} 1491--1518.
 \bibitem[Frumento and Bottai, 2016]{FB16} \textsc{Frumento P.} and \textsc{Bottai, M.} (2016). Parametric modeling of quantile regression coefficient functions. \textit{Biometrics} \textbf{72} 74--84.
 \bibitem[Giraud, 2015]{G15} \textsc{Giraud, C.} (2015). \textit{Introduction to high-dimensional statistics}. Chapman \& Hall/CRC Monographs on Statistics \& Applied Probability
\bibitem[Hastie, Tibshirani and Wainwright, 2015]{HTW15}  \textsc{Hastie, T.} \textsc{Tibshirani, R.} and  \textsc{Wainwright, M. J.} (2015). \textit{Statistical Learning with Sparsity: the Lasso and Generalizations}. Chapman \& Hall/CRC Press, Series in Statistics and Applied Probability.
 \bibitem[Huang and Zhang, 2013]{HZ13} \textsc{Huang, Z.} and \textsc{Zhang, R.}  (2013). Profile empirical-likelihood inferences for the single-index-coefficient regression model. \textit{Statistics and Computing} \textbf{23} 455--465.
 \bibitem[Hunter and Li, 2005]{HL05} \textsc{Hunter, D.} and \textsc{Li, R.} (2005). Variable selection using MM algorithms. \textit{The Annals of Statistics} \textbf{33} 1617--1642.
 \bibitem[Ji, Peng, Cheng and Lai]{JPCL12}\textsc{Ji, S.} \textsc{Peng, L.} \textsc{Cheng, Y.} and \textsc{Lai, H.} (2012). Quantile regression for doubly censored data. \textit{Biometrics} \textbf{68} 101--112.
 \bibitem[Jiang and Qian, 2016]{JQ16} \textsc{Jiang, R.} and \textsc{Qian, W.} (2016). Quantile regression for single-index-coefficient regression models. \textit{Statistics and Probability Letters} \textbf{110} 305-317.
 \bibitem[Jiang, Zhou and Qian, 2013]{JZQC13} \textsc{Jiang, R.} \textsc{Zhou, Z.} \textsc{Qian, W.} and \textsc{Chen, Y.} (2013). Two step composite quantile regression for single-index models. \textit{Computational Statistics and Data Analysis}  \textbf{64} 180--191.
 \bibitem[Jin, Lin, Wei and Ying, 2003]{JLWY03} \textsc{Jin, Z.} \textsc{Lin, D. Y.} \textsc{Wei, L. J.} and \textsc{Ying, Z.} (2003). Rank-Based Inference for the Accelerated Failure Time Model. \textit{Biometrika} \textbf{90} 341--353.
      \bibitem[Koenker, 2005]{K05} \textsc{Koenker, R.} (2005). \textit{Quantile regression}. Combridge University Press, New York.
 \bibitem[Koenker and Bassett, 1978 ]{KB78} \textsc{Koenker, R.} and \textsc{Bassett, G.} (1978). Regression quantiles. \textit{Econometrica} \textbf{ 46} 33--50.
 \bibitem[Kong and Xia, 2012]{KX12} \textsc{Kong, E.} and \textsc{Xia, Y.} (2012). A single-index quantile regression model and its estimation. \textit{Econometric Theory} \textbf{ 28} 730--768.
\bibitem[Li and Peng, 2015]{LP15} \textsc{Li, R.} and \textsc{Peng, L.} (2015). Quantile regression adjusting for dependent censoring from Semicompeting risks. \textit{Journal of the Royal Statistical Society Series B} \textbf{77} 107--130.
 \bibitem{L12} \textsc{Lian, H.} (2012). A note on the consistency of Schwarz's criterion in linear quantile regression with the SCAD penalty. \textit{Statistics and Probability Letters} \textbf{82} 1224--1228.
 \bibitem{LY12} \textsc{Lin, Z.} and \textsc{Yuan, Y.} (2012). Variable selection for generalized varying coefficient partially linear models with diverging number of parameters. \textit{Acta Mathematicae Applicatae Sinica} \textbf{28} 237--246.
 \bibitem{MH16} \textsc{Ma, S.} and \textsc{He, X.} (2016). Inference for single-index quantile regression models with profile optimization. \textit{The Annals of Statistics} \textbf{44} 1234--1268.
 \bibitem[Ma and Song, 2015]{MS15} \textsc{Ma, S.} and \textsc{Song, P. X.-K.} (2015). Varying index coefficient models. \textit{Journal of the American Statistical Association} \textbf{110} 341--356.
      \bibitem{MRT2009} \textsc{Marimoutou, V.} \textsc{Raggad, B.} \textsc{Trabelsi, A.} (2009). Extreme value theory and value at risk: application to oil market.  \textit{Energy Economics} \textbf{31} 519--530.
\bibitem{PXK14} \textsc{Peng, L.} \textsc{Xu, J.} and \textsc{Kutner, N.} (2014). Shrinkage estimation of varying covariate effects based on quantile regression. \textit{Statistics and Computing} \textbf{24} 853--869.
 \bibitem{SPMM16} \textsc{Sun, X.} \textsc{Peng, L.} \textsc{Manatunga, A.} \textsc{Marcus, M.} (2016). Quantile regression analysis of censored longitudinal data with irregular outcome-dependent follow-up. \textit{Biometrics} \textbf{72} 64--73.
      \bibitem{TWZ13} \textsc{Tang, Y.} \textsc{Wang, H.J.} and \textsc{Zhu, Z.} (2013). Variable selection in quantile varying coefficient models with longitudinal data. \textit{Computational Statistics and Data Analysis} \textbf{57} 435--449.

 \bibitem{WW15} \textsc{Wang, G.} and \textsc{Wang, L.} (2015). Spline estimation and variable selection for single-index prediction models with diverging number of index parameters. \textit{Journal of Statistical Planning and Inference} \textbf{162} 1--19.
 \bibitem{WLL09} \textsc{Wang, H.} \textsc{Li, B.} and \textsc{Leng, C.} (2009). Shrinkage tuning parameter selection with a diverging number of parameters. \textit{Journal of the Royal Statistical Society, Series B} \textbf{71} 671--683.
 \bibitem{WSZ12} \textsc{Wang, H.J.} \textsc{Stefanski, L.A.} and \textsc{Zhu, Z.} (2012). Corrected-loss estimation for quantile regression with covariate measurement error. \textit{Biometrika} \textbf{99} 405-421.
 \bibitem{WZL13} \textsc{Wang, H.J.} \textsc{Zhou, J.} and \textsc{Li, Y.} (2013). Variable selection for censored quantile regression. \textit{Statistica Sinica} \textbf{23} 145--167.
\bibitem{WZ11} \textsc{Wang, H.J.} and \textsc{Zhu, Z.} (2011). Empirical likelihood for quantile regression models with longitudinal data. \textit{Journal of Statistical Planning and Inference} \textbf{141} 1603--1615.

\bibitem{WXZC10} \textsc{Wang, J. L.} \textsc{Xue, L.} \textsc{Zhu, L.} and \textsc{Chong, Y. S.} (2010). Estimation for a partial linear single-index model. \textit{The Annals of Statistics} \textbf{38} 246--274.
 \bibitem{WZQ12} \textsc{Wang, L.} \textsc{Zhou, J.} and \textsc{Qu, A.} (2012). Penalized generalized estimating equations for high-dimensional longitudinal data analysis. \textit{Biometrics} \textbf{68} 353--360.
 \bibitem{W06} \textsc{Whang, Y. J.} (2006). Smoothed empirical likelihood methods for quantile regression models. \textit{Econometric Theory} \textbf{22} 173--205.
 \bibitem{WYY10} \textsc{Wu, T. Z.} \textsc{Yu, K.} and \textsc{Yu, Y.} (2010). Single-index quantile regression. \textit{Journal of Multivariate Analysis} \textbf{101} 1607--1621.
 \bibitem{XTLZ02} \textsc{Xia, Y.} \textsc{Tong, H.} \textsc{Li, W. K.} and \textsc{Zhu, L.} (2002). An adaptive estimation of dimension reduction space. \textit{Journal of the Royal Statistical Society Series B} \textbf{64} 363--410.
 \bibitem{XP13} \textsc{Xue, L.} and \textsc{Pang, Z.} (2013). Statistical inference for a single-index varying-coefficient model. \textit{Statistics and Computing} \textbf{23} 589--599.
 \bibitem{XQ12} \textsc{Xue, L.} and \textsc{Qu, A.} (2012). Variable selection in high-dimensional varying coefficient
models with global optimality. \textit{The Journal of Machine Learning Research} \textbf{ 13} 1973--1998.
 \bibitem{ZLL18} \textsc{Zhao, W.} \textsc{Li, J.} and \textsc{Lian, H.} (2018). Adaptive varying-coefficient linear quantile model: a profiled estimating equations approach. \textit{Annals of the Institute of Statistical Mathematics} \textbf{70} 553--582.
 \bibitem{ZL17} \textsc{Zhao, W.} and \textsc{Lian, H.} (2017). Quantile index coefficient model with variable selection. \textit{Journal of Multivariate Analysis} \textbf{154} 40--58.
 \bibitem{ZHL12} \textsc{Zhu, L.} \textsc{Huang, M.} and \textsc{Li, R.} (2012). Semiparametric quantile regression with high-dimensional covariates. \textit{Statistica Sinica} \textbf{22} 1379--1401.


\end{thebibliography}
\end{document}


\doublespacing

\begin{center}\textbf{SUPPLEMENTARY MATERIAL for ``HIGH-DIMENSIONAL VARYING INDEX COEFFICIENT
QUANTILE REGRESSION MODEL"}
\end{center}

\begin{aug}
\author{\fnms{JIALIANG} \snm{LI}\ead[label=e1]{stalj@nus.edu.sg}\thanksref{t1}}
\and
\author{\fnms{JING} \snm{LV}\corref{}\ead[label=e2]{lvjing@swu.edu.cn}\thanksref{t2}}

\runauthor{J. LI AND J. LV}
\affiliation{National University of Singapore\thanksmark{t1} and Southwest University\thanksmark{t2}}

\address{JIALIANG LI\\
Department of Statistics and Applied Probability\\ National University of Singapore\\
Singapore, 119077\\
\printead{e1}\\}

\address{JING LV\\
School of Mathematics and Statistics\\ Southwest University\\Chongqing, 400715, China\\
\printead{e2}\\}
\end{aug}

\begin{center}
\textbf{APPENDIX A: ADDITIONAL NUMERICAL RESULTS}
\end{center}

The tables of the Appendix A provide additional numerical results. Table S1 and S2 report the simulation results of estimators $\bm{\hat\beta}_l$ for two error distributions (MN and LA) in example 1.
\clearpage

\tabcolsep=7pt
\renewcommand\thetable{S1}
\begin{table}\scriptsize
\caption{Simulation results ($\times 10^{-2}$) of Bias, MAD, SD and SE for MN with $\tau=0.5$ in example 1. LS stands for the method by Ma and Song \cite{MS15} and QR is the proposed quantile regression.}
\label{table3}
\begin{tabular}{cccccccccccc} \noalign{\smallskip}\hline
\multicolumn{1}{c}{method}
&&&\multicolumn{2}{c}{$n=500$}
&&&&\multicolumn{2}{c}{$n=1500$}
\\
\cline{3-6}
\cline{8-11}
&&\multicolumn{1}{c}{Bias}&\multicolumn{1}{c}{MAD}&\multicolumn{1}{c}{SD}&\multicolumn{1}{c}{SE}
&&\multicolumn{1}{c}{Bias}&\multicolumn{1}{c}{MAD}&\multicolumn{1}{c}{SD}&\multicolumn{1}{c}{SE}

\\
\hline
LS&$\beta_{11}$&-0.352&  4.894& 6.439& 5.101&&-0.260&  2.747 &3.408 &2.845\\
&$\beta_{12}$&0.047&  6.774 &8.846& 7.111&&-0.292 & 3.644& 4.588 &3.901\\
&$\beta_{13}$&-0.678 & 3.680 &4.873 &3.597&&0.033 & 1.882& 2.333& 1.935\\

&$\beta_{21}$&-0.172&  2.706& 3.473& 2.964&&-0.096&  1.408 &1.776& 1.737\\
&$\beta_{22}$&-0.552&  4.603& 5.893& 4.643&&0.157 & 2.276 &2.874 &2.695\\
&$\beta_{23}$&0.211&  4.039 &5.327 &5.245&&-0.376 & 2.105 &2.682& 3.019\\

&$\beta_{31}$&-0.094&  1.760& 2.318& 1.911&&0.034  &0.987& 1.246& 1.085\\
&$\beta_{32}$&-0.048&  1.456 &1.856& 1.612&&-0.046 & 0.772& 0.983& 0.911\\
&$\beta_{33}$& 0.053 & 1.921 &2.466& 2.285&&-0.010 & 1.019 &1.294& 1.293\\
\hline
QR&$\beta_{11}$& -0.336 & 3.655 &4.694& 3.869&&0.050 & 1.991& 2.483& 2.019\\
&$\beta_{12}$&0.051 & 5.217& 6.706& 5.340&&-0.002 & 2.620& 3.331 &2.754\\
&$\beta_{13}$&-0.288&  2.674& 3.517& 2.619&&-0.159 & 1.384& 1.719& 1.357\\

&$\beta_{21}$&-0.124&  1.854& 2.298 &2.323&&-0.039 & 0.996 &1.250& 1.263\\
&$\beta_{22}$&-0.293&  3.024& 3.769 &3.617&&0.148 & 1.603& 2.019 &1.957\\
&$\beta_{23}$&0.383&  2.646& 3.327 &4.018&&-0.351 & 1.513& 1.859 &2.172\\

&$\beta_{31}$&-0.023&  1.255 &1.623 &1.453&&0.027 & 0.634 &0.817 &0.775\\
&$\beta_{32}$&-0.067 & 1.032& 1.327 &1.223&&-0.027 & 0.518& 0.667& 0.649\\
&$\beta_{33}$&0.105 & 1.426 &1.784& 1.735&&-0.008 & 0.740& 0.912 &0.914\\
\hline
\end{tabular}
\end{table}

\tabcolsep=7pt
\renewcommand\thetable{S2}
\begin{table}\scriptsize
\caption{Simulation results ($\times 10^{-2}$) of Bias, MAD, SD and SE for LA with $\tau=0.5$ in example 1. LS stands for the method by Ma and Song \cite{MS15} and QR is the proposed quantile regression.}
\label{table4}
\begin{tabular}{cccccccccccc} \noalign{\smallskip}\hline
\multicolumn{1}{c}{method}
&&&\multicolumn{2}{c}{$n=500$}
&&&&\multicolumn{2}{c}{$n=1500$}
\\
\cline{3-6}
\cline{8-11}
&&\multicolumn{1}{c}{Bias}&\multicolumn{1}{c}{MAD}&\multicolumn{1}{c}{SD}&\multicolumn{1}{c}{SE}
&&\multicolumn{1}{c}{Bias}&\multicolumn{1}{c}{MAD}&\multicolumn{1}{c}{SD}&\multicolumn{1}{c}{SE}

\\
\hline
LS&$\beta_{11}$&-0.141&  4.030& 5.052 &4.069&&0.088 & 2.181 &2.712& 2.252\\
&$\beta_{12}$&0.446&  5.318& 6.854 &5.677&&-0.085 & 2.996& 3.902& 3.119\\
&$\beta_{13}$&-0.595&  2.901& 3.693 &2.864&&-0.193&  1.490& 1.860 &1.555\\

&$\beta_{21}$&0.022&  1.937 &2.476 &2.317&&-0.194&  1.118 &1.433& 1.359\\
&$\beta_{22}$&-0.423&  3.303& 4.156& 3.623&&0.289 & 1.858 &2.321& 2.093\\
&$\beta_{23}$&0.089&  2.842 &3.662& 4.059&&-0.214 & 1.595 &2.025& 2.328\\

&$\beta_{31}$&0.022 & 1.255 &1.617& 1.496&&0.009 & 0.740& 0.929& 0.845\\
&$\beta_{32}$&-0.119&  1.014& 1.304& 1.270&&-0.049&  0.580& 0.725& 0.713\\
&$\beta_{33}$&0.166&  1.509& 1.856 &1.773&&0.083 & 0.820& 1.023& 0.999\\
\hline
QR&$\beta_{11}$&-0.066 & 3.606& 4.604& 3.785&&0.070&  1.813 &2.294& 1.876\\
&$\beta_{12}$&0.339 & 4.679 &6.112& 5.199&&0.000&  2.504 &3.173 &2.613\\
&$\beta_{13}$&-0.501&  2.533& 3.241& 2.633&&-0.158&  1.259& 1.585& 1.269\\

&$\beta_{21}$&-0.073 & 1.684& 2.132& 2.238&&-0.132&  0.868& 1.105 &1.151\\
&$\beta_{22}$&-0.167&  2.856& 3.583& 3.482&&0.225&  1.428 &1.783 &1.789\\
&$\beta_{23}$&0.035 & 2.512& 3.212& 3.897&&-0.183&  1.267& 1.564 &1.945\\

&$\beta_{31}$&0.029 & 1.172 &1.477& 1.399&&-0.014 & 0.615& 0.763 &0.713\\
&$\beta_{32}$&-0.138&  0.986& 1.225& 1.189&&-0.004 & 0.468& 0.598 &0.602\\
&$\beta_{33}$&0.230&  1.391& 1.726 &1.658&&0.010 & 0.666& 0.832 &0.840\\
\hline
\end{tabular}
\end{table}

\begin{center}
\textbf{APPENDIX B: REGULAR CONDITIONS AND PROOFS OF TECHNICAL LEMMAS AND THEOREMS}
\end{center}

In this document, we present the regular conditions and the proofs of technical lemmas and theorems. The technical lemmas are used to prove the theorems of the paper.
\\

We use $C$ to denote a generic positive constant that can vary from line to line. For any vector $\bm\xi  = {\left( {{\xi _1},...,{\xi _s}} \right)^T} \in {\mathbb{R}^s}$, denote ${\left\| \bm\xi  \right\|_\infty } = {\max _{1 \le l \le s}}\left| {{\xi _l}} \right|$. For any symmetric matrix $\bm A_{s\times s}$, denote its $L_r$ norm as ${\left\|\bm A \right\|_r} = {\max _{\bm\xi  \in {\mathbb{R}^s},\bm\xi  \ne 0}}{\left\| {\bm A\bm\xi } \right\|_r}\left\| \bm\xi  \right\|_r^{ - 1}$. For any matrix $\bm A = \left( {{A_{ij}}} \right)_{i = 1,j = 1}^{s,t}$, denote
 ${\left\| \bm A \right\|_\infty } = {\max _{1 \le i \le s}}\sum\nolimits_{j = 1}^t {\left| {{A_{ij}}} \right|} $. Denote the space of the $r$th order smooth function as
${C^{\left( r \right)}}\left[ {0,1} \right] = \left\{ {\varphi \left| {{\varphi ^{(r)}} \in C\left[ {0,1} \right]} \right.} \right\}$. Throughout the paper, $\rho_{min}$ and $\rho_{max}$ denote the minimum and maximum eigenvalues of a matrix and $tr\left( \cdot \right)$ is the trace operator of a matrix.
To prove the theoretical results of the proposed estimators, we next present the following technical conditions.

\noindent (C1) For $\bm \beta_l$ in the neighborhood of $\bm \beta_l^0$, the density function ${f_{{U_l}({\bm \beta _l})}}( \cdot )$ of random variable ${U_l}({\bm \beta _l}) = {\bm Z^T}{\bm \beta _l}$ is bounded away from zero on $[0, 1]$ for $1\leq l\leq d$ and satisfies the Lipschitz condition of order 1 on $[0, 1]$.

\noindent (C2) For every $1\leq l\leq d$ and $1\leq k\leq p_n$, $g_{l,k}^0 \in C^{(1)}[0,1] $ and ${m_l} \in C^{(r)}[0,1]$ for some integer $r\geq 2$. Meanwhile, the spline order $q$ satisfies $q\geq r$.

\noindent (C3) There are $N_0$ and constants $b$ and $B$ such that, for $n\geq N_0$,
\[0 < b \le {\rho _{min}}\left( {{n^{ - 1}}\sum\limits_{i = 1}^n {{\bm X_i}\bm X_i^T} } \right) \le {\rho _{max}}\left( {{n^{ - 1}}\sum\limits_{i = 1}^n {{\bm X_i}\bm X_i^T} } \right) \le B.\]

 \noindent (C4) ${\sup _i}{\left\| {{\bm Z_i}} \right\|_2} = O\left( {\sqrt {{p_n}} } \right)$.

 \noindent (C5) {For a fixed $\tau$, we assume $\varepsilon_{\tau,1},...,\varepsilon_{\tau,n}$ are independent identically distributed errors and the conditional density function $f_{\tau}(\cdot|\bm X,\bm Z)$ is bounded and bounded away from zero in a neighborhood of zero, uniformly over the support of $\bm X$ and $\bm Z$. The
derivative of $f_{\tau}(\cdot|\bm X,\bm Z)$ is uniformly bounded in a neighborhood of zero over the support of $\bm X$ and $\bm Z$. In the following proof, we omit the subscript $\tau$ from  $\varepsilon _{\tau,i}$ and $f_{\tau}\left(\cdot|\bm X, \bm Z \right)$ to simplify the notations.}

\noindent (C6) $K\left(\cdot\right)$ is bounded and compactly supported on $[-1,1]$. Furthermore, we assume that $K\left(\cdot\right)$ is a $\nu$th order kernel function. For example, $\int {{u^j}} K\left( u \right)du = 1$ if $j=0$; 0 if $1 \le j \le \nu - 1$; $C_K$ if $j=\nu$, where $\nu\geq 2$ is an integer and $C_K\neq 0$ is a constant.

\noindent (C7) The bandwidth $h$ satisfies $n{h^{2 \nu}} \to 0$.

\noindent (C8) $\lim {\inf _{n \to \infty }}\lim {\inf _{x \to {0_ + }}}{\dot{p}_{{\alpha _1}}}\left( x \right)/{\alpha _1} > 0$.

\noindent (C9) ${a_n} = {\max _{2 \le j \le p_n, 1 \le l \le d}}\{ {\dot{p}_{\alpha_1}}(\left| {\beta _{lj}^0} \right|),\beta _{lj}^0\ne 0\} = O\left( n^{-1/2} \right)$.

\noindent (C10) ${b_n} = {\max _{2 \le j \le {p_n}, 1 \le l \le d}}\{ {\ddot{p}_{\alpha_1}}(\left| {\beta _{lj}^0} \right|),\beta _{lj}^0 \ne 0\}\to 0 $ as $n \to \infty $.

\noindent (C11) There are constants $C_1$ and $C_2$ such that $\left| {{\ddot{p}_{{\alpha _1}}}\left( {{x_1}} \right) - {\ddot{p}_{{\alpha _1}}}\left( {{x_2}} \right)} \right| \le {C_2}\left| {{x_1} - {x_2}} \right|$ when $x_1,x_2>C_1\alpha_1$.

\noindent (C12) Assume $\left\{\beta_{l2}^{0},...,\beta_{ls_l}^{0}\right\}_{l=1}^d$ satisfy ${\min _{1 \le l \le d,2 \le j \le s_l}}{{\left| {\beta _{lj}^0} \right|} \mathord{\left/
 {\vphantom {{\left| {\beta _{lj}^0} \right|} {{\alpha _1}}}} \right.
 \kern-\nulldelimiterspace} {{\alpha _1}}} \to \infty$ as $\to \infty$.

Conditions (C1)--(C2) are standard conditions for VICM similar to (C1), (C2) and (C5) in Ma and Song \cite{MS15}. Condition (C3) is a common assumption on $M$-estimation with diverging dimension. Condition (C4) is also popularly adopted in the literature on high-dimensional regression for independent data. Condition (C5) on conditional density is commonly used in quantile regression. Conditions (C6)--(C7) are the necessary conditions on the kernel function and bandwidth $h$, which are also required in Whang \cite{W06}. Conditions (C8)--(C11) are the regularity conditions on the penalty given in Fan and Peng \cite{FP04} and condition (C12) is similar to the condition (H) of Fan and Peng \cite{FP04}, which is used to obtain the oracle property.

\begin{lemma}\label{lemma:1}
 Under conditions (C1)--(C5), $J_n\rightarrow \infty$ and $n^{-1}J_n^{3} =o(1)$, for $1\leq l\leq d$, we have
 \[(i)\left| {{{\hat m}_l}({u_l},{\bm \beta ^0}) - {m_l}({u_l})} \right| = O_p\left( {\sqrt {{{{J_n}} \mathord{\left/
 {\vphantom {{{J_n}} n}} \right.
 \kern-\nulldelimiterspace} n}}  + J_n^{ - r}} \right), \left| {{{\hat {\dot m}}_l}({u_l},{\bm \beta ^0}) - {{\dot m}_l}({u_l})} \right| = O_p\left( {\sqrt {{{J_n^3} \mathord{\left/
 {\vphantom {{J_n^3} n}} \right.
 \kern-\nulldelimiterspace} n}}  + J_n^{ - r + 1}} \right) \]
uniformly for any $u_l\in [0,1]$,

\noindent $(ii)$ under $nJ_n^{ - \left( {1 + 2r} \right)} = o\left( 1 \right)$, $\sigma _{nl}^{ - 1}\left( {{u_l}} \right) ( {{\hat m}_{l }}\left( {{u_l},\bm\beta ^0} \right)-{m_l^0}({u_l})  )\mathop  \to \limits^d N\left( {0,1} \right),$ where $\sigma _{nl}^2(u_l)$ is defined in subsection 2.2.
  \end{lemma}

\begin{proof}
 According to the result on page 149 of de Boor \cite{De01}, for ${m _l}(u_l)$ satisfying condition (C2), we have
\begin{gather}
\mathop {\sup }\limits_{u _l\in [0,1]} \left| {m _l}(u_l) - m_l^0(u_l) \right| =O({J_n^{ - r}}).
\tag{A.1}\label{A.1}
\end{gather}
Let ${\bm \lambda ^0} = {\left\{ {{{\left( {\bm \lambda _l^0} \right)}^T}:1 \le l \le d} \right\}^T}$, ${k_n} = {J_n}{n^{{{ - 1} \mathord{\left/
 {\vphantom {{ - 1} 2}} \right.
 \kern-\nulldelimiterspace} 2}}} + J_n^{ - r + {1 \mathord{\left/
 {\vphantom {1 2}} \right.
 \kern-\nulldelimiterspace} 2}}$ and a vector $\bm v$ satisfying $\left\| \bm v \right\|_2 = C$, where $C$ is a sufficiently large constant. Our aim is to show that for
 any given $\epsilon$ there is a large constant $C$ such that, for a large $n$, we have
\begin{gather}
P\left\{ {\mathop {\inf }\limits_{{{\left\| \bm v \right\|}_2} = C} {{\rm{{\cal L}}}_{\tau n}}\left( {{\bm \lambda ^0} + { k_n}\bm v,\bm \beta ^0} \right) \ge {{\rm{{\cal L}}}_{\tau n}}\left( {{\bm \lambda ^0},\bm \beta ^0} \right)} \right\} \ge 1 -\epsilon.
\tag{A.2}\label{A.2}
\end{gather}
This implies with probability of at least $ 1 -\epsilon$ that there exists a local minimizer in the ball
$\left\{ {{\bm \lambda ^0} + {k_n}\bm v:{{\left\| \bm v \right\|}_2} \le C} \right\}$. Then it follows that there exists a local minimizer $\bm{\hat \lambda}$ such that ${\left\| {\bm{\hat{  \lambda}}  - {\bm \lambda ^0}} \right\|_2} = {O_p}\left( {{k_n}} \right)$. By Knight's identity,
\[\left( {\left| {r - s} \right| - \left| r \right|} \right)/2 =  - s(\frac{1}{2} - I(r < 0)) + \int_0^s {[I(r \le t) - I(r \le 0)]} dt.\]
we have
\[{\rho _\tau }\left( {r - s} \right) - {\rho _\tau }\left( r \right) = s\left( {I\left( {r < 0} \right) - \tau } \right) + \int_0^s {\left[ {I\left( {r \leq t} \right) - I\left( {r \leq 0} \right)} \right]dt}. \]
Let ${\Delta _i} = \sum\limits_{l = 1}^d {{m_l}\left( {\bm Z_i^T\bm\beta _{l}^0} \right){X_{il}}}  - \sum\limits_{l = 1}^d {{\bm B}{{\left( {\bm Z_i^T\bm\beta _{l}^0} \right)}^T}{\bm \lambda_l^0}{X_{il}}} $, we have
\[\begin{array}{l}
  {{\rm{{\cal L}}}_{\tau n}}\left( {{\bm\lambda ^0},\bm\beta ^0} \right)- {{\rm{{\cal L}}}_{\tau n}}\left( {{\bm\lambda ^0} + {k_n}\bm v,\bm\beta ^0} \right) \\
  = {k_n}\sum\limits_{i = 1}^n {\left[ {\tau  - I\left\{ {{\varepsilon _i} + {\Delta _i} \le 0} \right\}} \right]} \left\{ {\sum\limits_{l = 1}^d {{\bm B}{{(\bm Z_i^T\bm\beta _{l}^0)}^T}{\bm v_l}{X_{il}}} } \right\} \\
 ~~~ - \sum\limits_{i = 1}^n {\int_0^{{k_n}\sum\limits_{l = 1}^d {{\bm B}{{(\bm Z_i^T\bm\beta _{l}^0)}^T}{\bm v_l}{X_{il}}} } {\left[ {I\left\{ {{\varepsilon _i} + {\Delta _i} \le t} \right\} - I\left\{ {{\varepsilon _i} + {\Delta _i} \le 0} \right\}} \right]dt} }  \\
  \buildrel \Delta \over = I + II. \\
 \end{array}\]
Meanwhile, the first term can be divided into the following two parts.
\begin{gather}
\begin{array}{l}
 I = {k_n}\sum\limits_{i = 1}^n {\left[ {\tau  - I\left\{ {{\varepsilon _i} \le 0} \right\}} \right]} \left\{ {\sum\limits_{l = 1}^d {\bm B{{(\bm Z_i^T\bm \beta _l^0)}^T}{\bm v_l}{X_{il}}} } \right\} \\
  ~~~~~+ {k_n}\sum\limits_{i = 1}^n {\left[ {I\left\{ {{\varepsilon _i} \le 0} \right\} - I\left\{ {{\varepsilon _i} + {\Delta _i} \le 0} \right\}} \right]} \left\{ {\sum\limits_{l = 1}^d {\bm B{{(\bm Z_i^T\bm\beta _l^0)}^T}{\bm v_l}{X_{il}}} } \right\} \\
 ~~ = {I_1} + {I_2} .\\
 \end{array}
 \tag{A.3}\label{A.3}
\end{gather}
Thus $E\left( I_1 \right) =0$. By Theorem 5.4.2 of DeVore and Lorentz \cite{DL93}, for sufficiently large $n$,
we can prove that under conditions (C1) and (C3), there are constants $0 < c_1 \leq C_1 < \infty$, such that
\begin{gather}
{c_1}\left\| \bm v \right\|_2^2J_n^{ - 1} \le E{\left\{ {\sum\limits_{l = 1}^d {{\bm B}{{(\bm Z_i^T\bm \beta _{l}^0)}^T}{\bm v_l}{X_{il}}} } \right\}^2} \le {C_1}\left\| \bm v \right\|_2^2J_n^{ - 1}
\tag{A.4}\label{A.4}
\end{gather}
and by Bernstein's inequality in Bosq \cite{B98}, we have with probability 1,
\begin{gather}
{c_1}\left\| \bm v \right\|_2^2J_n^{ - 1} \le {n^{ - 1}}\sum\limits_{i = 1}^n {{{\left\{ {\sum\limits_{l = 1}^d {{\bm B}{{(\bm Z_i^T\bm \beta _{l}^0)}^T}{\bm v_l}{X_{il}}} } \right\}}^2}}  \le {C_1}\left\| \bm v \right\|_2^2J_n^{ - 1}.
\tag{A.5}\label{A.5}
\end{gather}
Therefore, by (\ref{A.4}), there is a constant $0 < {C_2} < \infty $ such that for sufficiently large $n$,
\[Var\left( I_1 \right) \le {C_2}nk_n^2E{\left\{ {\sum\limits_{l = 1}^d {{\bm B}{{(\bm Z_i^T\bm \beta _{l}^0)}^T}{\bm v_l}{ X_{il}}} } \right\}^2} \le nk_n^2{C_1}{C_2}{C^2}J_n^{ - 1}.\]
By weak law of large numbers, we have with probability approaching 1, for sufficiently large
$n$, $\left| I_1 \right| \le {n^{{1 \mathord{\left/
 {\vphantom {1 2}} \right.
 \kern-\nulldelimiterspace} 2}}}{k_n}{C_3}CJ_n^{{-1 \mathord{\left/
 {\vphantom {-1 2}} \right.
 \kern-\nulldelimiterspace} 2}}$ for some constant  $0 < {C_3} < \infty $. In addition,
\begin{gather}
\begin{array}{l}
 {I_2} = {k_n}\sum\limits_{i = 1}^n {\left[ {F\left( 0\left| {\bm X_i, \bm Z_i} \right. \right) - F\left( { - {\Delta _i}} \left| {\bm X_i, \bm Z_i} \right.\right)} \right]} \left\{ {\sum\limits_{l = 1}^d {\bm B{{(\bm Z_i^T\bm \beta _l^0)}^T}{\bm v_l}{X_{il}}} } \right\} \\
 ~~~~~~ + {k_n}\sum\limits_{i = 1}^n
 \left\{ {\left[ {I\left\{ {{\varepsilon _i} \le 0} \right\} - I\left\{ {{\varepsilon _i} + {\Delta _i} \le 0} \right\}} \right]} \right. \\
 ~~~~~~-\left. {\left[ {F\left( 0\left| {\bm X_i, \bm Z_i} \right.  \right) - F\left( { - {\Delta _i}} \left| {\bm X_i, \bm Z_i} \right. \right)} \right]} \right\}\left\{ {\sum\limits_{l = 1}^d {\bm B{{(\bm Z_i^T\bm\beta _l^0)}^T}{\bm v_l}{X_{il}}} } \right\} \\
 ~~ = {I_{21}} + {I_{22}}.\\
 \end{array}
 \tag{A.6}\label{A.6}
\end{gather}
Moreover, by (\ref{A.1}) and (\ref{A.5}), taking Taylor's explanation for $F\left( { - {\Delta _i}}\left| {\bm X_i, \bm Z_i} \right.  \right)$ at 0 gives
\begin{gather}
\begin{array}{l}
 {I_{21}} = {k_n}\sum\limits_{i = 1}^n {\left[ {f\left( 0\left| {\bm X_i, \bm Z_i} \right.  \right){\Delta _i}\left( {1 + o\left( 1 \right)} \right)} \right]} \left\{ {\sum\limits_{l = 1}^d {\bm B{{(\bm Z_i^T\bm \beta _l^0)}^T}{\bm v_l}{X_{il}}} } \right\} \\
 ~~~ \le n{k_n}C_1^{{1 \mathord{\left/
 {\vphantom {1 2}} \right.
 \kern-\nulldelimiterspace} 2}}CJ_n^{{{ - 1} \mathord{\left/
 {\vphantom {{ - 1} 2}} \right.
 \kern-\nulldelimiterspace} 2}}O\left( {J_n^{ - r}} \right) \\
  ~~~\le {C_4}C n{k_n}J_n^{{{ - 1} \mathord{\left/
 {\vphantom {{ - 1} {2 - r}}} \right.
 \kern-\nulldelimiterspace} {2 - r}}} \\
 \end{array}
 \tag{A.7}\label{A.7}
\end{gather}
for some constant $0<C_4<\infty$. By direct calculation of the mean and variance, we can show that $I_{22}={o_p}\left( {n{k_n}J_n^{{{ - 1} \mathord{\left/
 {\vphantom {{ - 1} {2 - r}}} \right.
 \kern-\nulldelimiterspace} {2 - r}}}} \right)$.
 This combined with (\ref{A.3}), (\ref{A.6}), (\ref{A.7}) and the above results leads to
\begin{gather}
\begin{array}{l}
 I \le {C_3}C{n^{{1 \mathord{\left/
 {\vphantom {1 2}} \right.
 \kern-\nulldelimiterspace} 2}}}{k_n}J_n^{{{ - 1} \mathord{\left/
 {\vphantom {{ - 1} 2}} \right.
 \kern-\nulldelimiterspace} 2}} + {C_4}Cn{k_n}J_n^{{{ - 1} \mathord{\left/
 {\vphantom {{ - 1} {2 - r}}} \right.
 \kern-\nulldelimiterspace} {2 - r}}} \\
  ~~= {C_3}C{n^{{1 \mathord{\left/
 {\vphantom {1 2}} \right.
 \kern-\nulldelimiterspace} 2}}}\left( {{{{J_n}} \mathord{\left/
 {\vphantom {{{J_n}} {{n^{{1 \mathord{\left/
 {\vphantom {1 2}} \right.
 \kern-\nulldelimiterspace} 2}}} + J_n^{{{ - r + 1} \mathord{\left/
 {\vphantom {{ - r + 1} 2}} \right.
 \kern-\nulldelimiterspace} 2}}}}} \right.
 \kern-\nulldelimiterspace} {{n^{{1 \mathord{\left/
 {\vphantom {1 2}} \right.
 \kern-\nulldelimiterspace} 2}}} + J_n^{{{ - r + 1} \mathord{\left/
 {\vphantom {{ - r + 1} 2}} \right.
 \kern-\nulldelimiterspace} 2}}}}} \right)J_n^{{{ - 1} \mathord{\left/
 {\vphantom {{ - 1} 2}} \right.
 \kern-\nulldelimiterspace} 2}}  {C_4}Cn\left( {{{{J_n}} \mathord{\left/
 {\vphantom {{{J_n}} {{n^{{1 \mathord{\left/
 {\vphantom {1 2}} \right.
 \kern-\nulldelimiterspace} 2}}} + J_n^{{{ - r + 1} \mathord{\left/
 {\vphantom {{ - r + 1} 2}} \right.
 \kern-\nulldelimiterspace} 2}}}}} \right.
 \kern-\nulldelimiterspace} {{n^{{1 \mathord{\left/
 {\vphantom {1 2}} \right.
 \kern-\nulldelimiterspace} 2}}} + J_n^{{{ - r + 1} \mathord{\left/
 {\vphantom {{ - r + 1} 2}} \right.
 \kern-\nulldelimiterspace} 2}}}}} \right)J_n^{{{ - 1} \mathord{\left/
 {\vphantom {{ - 1} {2 - r}}} \right.
 \kern-\nulldelimiterspace} {2 - r}}} \\
   ~~= {C_3}CJ_n^{{1 \mathord{\left/
 {\vphantom {1 2}} \right.
 \kern-\nulldelimiterspace} 2}} + {C_3}C{n^{{1 \mathord{\left/
 {\vphantom {1 2}} \right.
 \kern-\nulldelimiterspace} 2}}}J_n^{ - r} + {C_4}C{n^{{1 \mathord{\left/
 {\vphantom {1 2}} \right.
 \kern-\nulldelimiterspace} 2}}}J_n^{{{ - r + 1} \mathord{\left/
 {\vphantom {{ - r + 1} 2}} \right.
 \kern-\nulldelimiterspace} 2}} + {C_4}CnJ_n^{ - 2r}. \\
 \end{array}
  \tag{A.8}\label{A.8}
\end{gather}
Let $\bm{\mathcal{X }}= {\left( {\bm X_i^T,1 \le i \le n} \right)^T},\bm{\mathcal{Z} }= {\left( {\bm Z_i^T,1 \le i \le n} \right)^T}$ and ${\Upsilon _i} = {k_n}\sum\limits_{l = 1}^d {{\bm B}{{(\bm Z_i^T\bm \beta _{l}^0)}^T}{\bm v_l}{X_{il}}} $. By (\ref{A.1}) and (\ref{A.5}), we obtain
\begin{gather}
\begin{array}{l}
 E\left[ {II\left| {\bm {\mathcal{X }},\bm {\mathcal{Z }}} \right.} \right] \\
=  - E\left[ {\sum\limits_{i = 1}^n {\int_0^{{\Upsilon _i} } {\left[ {I\left\{ {{\varepsilon _i} + {\Delta _i} \le t} \right\} - I\left\{ {{\varepsilon _i} + {\Delta _i} \le 0} \right\}} \right]dt} } \left| {\bm {\mathcal{X }},\bm {\mathcal{Z }}} \right.} \right] \\
=  - \sum\limits_{i = 1}^n {\int_0^{{\Upsilon _i} } {\left[ {F\left( {(t - {\Delta _i})}\left| {\bm X_i, \bm Z_i} \right.  \right) - F\left( { - {\Delta _i}}\left| {\bm X_i, \bm Z_i} \right.  \right)} \right]dt} }  \\
=  - \sum\limits_{i = 1}^n {\int_0^{{\Upsilon _i} } {\left[ {f\left( { - {\Delta _i}} \left| {\bm X_i, \bm Z_i} \right. \right)t + o\left( t \right)} \right]dt} }  \\
=  - \frac{1}{2}k_n^2\sum\limits_{i = 1}^n {{{f\left( 0 | {\bm X_i, \bm Z_i} \right)\left\{ {\sum\limits_{l = 1}^d {{\bm B}{{(\bm Z_i^T\bm \beta _{l}^0)}^T}{\bm v_l}{X_{il}}} } \right\}}^2}} \left\{ {1 + o\left( 1 \right)} \right\} \\
 \le  - n{c_1}C_5k_n^2J_n^{ - 1}{C^2}\\
 \end{array}
 \tag{A.9}\label{A.9}
\end{gather}
for some constant $0<C_5<\infty$ and
\begin{gather}
\begin{array}{l}
 Var\left( {II\left| {\bm {\mathcal{X }},\bm {\mathcal{Z} }} \right.} \right)\\
  = Var\left[ { - \sum\limits_{i = 1}^n {\int_0^{{\Upsilon _i}} {\left[ {I\left\{ {{\varepsilon _i} + {\Delta _i} \le t} \right\} - I\left\{ {{\varepsilon _i} + {\Delta _i} \le 0} \right\}} \right]dt} } \left| {\bm {\mathcal{X} },\bm {\mathcal{Z }}} \right.} \right] \\
 \le \sum\limits_{i = 1}^n {E\left[ {{{\left( {\int_0^{{\Upsilon _i}} {\left[ {I\left\{ {{\varepsilon _i} + {\Delta _i} \le t} \right\} - I\left\{ {{\varepsilon _i} + {\Delta _i} \le 0} \right\}} \right]dt} } \right)}^2}\left| {\bm {\mathcal{X }},\bm {\mathcal{Z} }} \right.} \right]}  \\
  \le \sum\limits_{i = 1}^n {\int_0^{\left| {{\Upsilon _i}} \right|} {\int_0^{\left| {{\Upsilon _i}} \right|} {\left[ {F \left\{ \left( {\left| {{\Upsilon _i}} \right| - {\Delta _i}} \right) \left| {\bm X_i, \bm Z_i} \right. \right\}- F\left( { - {\Delta _i}}\left| {\bm X_i, \bm Z_i} \right.  \right)} \right]} } } d{t_1}d{t_2} \\
 \le o\left( {\sum\limits_{i = 1}^n {k_n^2} {{\left( {\sum\limits_{l = 1}^d {{\bm B}{{(\bm Z_i^T\bm \beta _{l}^0)}^T}{\bm v_l}{X_{il}}} } \right)}^2}} \right) \\
 \le o\left( {n{C_1}k_n^2J_n^{ - 1}} C^2\right).\\
 \end{array}
 \tag{A.10}\label{A.10}
\end{gather}
From (\ref{A.9}) and (\ref{A.10}), we have
\[II \le  - n{c_1}{C_5}k_n^2J_n^{ - 1}{C^2} =  - n{c_1}{C_5}{C^2}\left( {{{J_n^2} \mathord{\left/
 {\vphantom {{J_n^2} {n + J_n^{ - 2r + 1}}}} \right.
 \kern-\nulldelimiterspace} {n + J_n^{ - 2r + 1}}}} \right)J_n^{ - 1} =  - {c_1}{C_5}{C^2}{J_n} - {c_1}{C_5}{C^2}nJ_n^{ - 2r}.\]
Combined with (\ref{A.8}), we have
\[\begin{array}{l}
 {{\rm{{\cal L}}}_{\tau n}}\left( {{\bm \lambda ^0},\bm \beta ^0} \right)- {{\rm{{\cal L}}}_{\tau n}}\left( {{\bm \lambda ^0} + {k_n}\bm v,\bm \beta ^0} \right)  \\
\le{C_3}CJ_n^{{1 \mathord{\left/
 {\vphantom {1 2}} \right.
 \kern-\nulldelimiterspace} 2}} + {C_3}C{n^{{1 \mathord{\left/
 {\vphantom {1 2}} \right.
 \kern-\nulldelimiterspace} 2}}}J_n^{ - r} + {C_4}C{n^{{1 \mathord{\left/
 {\vphantom {1 2}} \right.
 \kern-\nulldelimiterspace} 2}}}J_n^{{{ - r + 1} \mathord{\left/
 {\vphantom {{ - r + 1} 2}} \right.
 \kern-\nulldelimiterspace} 2}} + {C_4}CnJ_n^{ - 2r}  {c_1}{C_5}{C^2}{J_n} - {c_1}{C_5}{C^2}nJ_n^{ - 2r}\\
 \le 0.
 \end{array}\]
Thus (\ref{A.2}) holds. That is, ${\left\| {\hat {\bm \lambda}({\bm\beta ^0})  - {\bm \lambda ^0}} \right\|_2} = {O_p}\left( {{k_n}} \right)$. By Taylor expansion, for $1\leq l \leq d$,
\begin{gather}
\begin{array}{l}
 {{\hat m}_l}({u_l},{\bm \beta ^0}) - m_l^0({u_l},{\bm\beta ^0}) \\
 = \bm e_l^T \mathbb{B}(\bm u)\left\{ {\bm{\hat \lambda} ({\bm\beta ^0}) - {\bm\lambda ^0}} \right\} \\
  =   -\bm e_l^T \mathbb{B}(\bm u){\left\{ {{{{\partial ^2}{{\cal L}_{\tau n}}\left( {{\bm \lambda ^0},\bm \beta ^0} \right)} \mathord{\left/
 {\vphantom {{{\partial ^2}{\rm{{\cal L}}_{\tau n}}\left( {{\bm \lambda ^0},\bm \beta _n^0} \right)} {\partial \lambda \partial {\bm \lambda ^T}}}} \right.
 \kern-\nulldelimiterspace} {\partial \bm \lambda \partial {\bm \lambda ^T}}}} \right\}^{ - 1}}\left\{ {{{\partial {{\cal L}_{\tau n}}\left( {{\bm \lambda ^0},\bm \beta ^0} \right)} \mathord{\left/
 {\vphantom {{\partial {\rm{{\cal L}}_{\tau n}}\left( {{\bm \lambda ^0},\bm \beta ^0} \right)} {\partial \bm \lambda }}} \right.
 \kern-\nulldelimiterspace} {\partial \bm \lambda }}} \right\}\left\{ {1 + {o_p}(1)} \right\} \\
  = \left\{ {{{\hat m}_{l,\varepsilon }}({u_l},{\bm \beta ^0}) + {{\hat m}_{l,m}}({u_l},{\bm \beta ^0})} \right\}\left\{ {1 + {o_p}(1)} \right\},
 \end{array}
  \tag{A.11}\label{A.11}
\end{gather}
where
 \[{{\hat m}_{l,\varepsilon }}({u_l},{\bm \beta ^0}) = \bm e_l^T\mathbb{B}(\bm u){\bm \Omega _n}{({\bm \beta ^0})^{ - 1}} n^{ - 1}{\sum\limits_{i = 1}^n  \psi_\tau (\varepsilon_i) {{\bm D_i}({\bm \beta ^0})} } , \]
  \[{{\hat m}_{l,m}}({u_l},{\bm \beta ^0}) = \bm e_l^T\mathbb{B}(\bm u){\bm \Omega _n}{({\bm \beta ^0})^{ - 1}}{n^{ - 1}} \sum\limits_{i = 1}^n{f\left( 0\left| {\bm X_i, \bm Z_i} \right. \right)} \Delta_i {\bm D_i}({\bm \beta ^0}),\]
and ${\bm \Omega _n}({\bm \beta ^0})$, $\mathbb{B}(\bm u)$ and $\bm e_l$ are given in subsection 2.2.

Then $E\left\{ {{\hat m}_{l,\varepsilon}}\left( {{u_l},\bm \beta ^0} \right)\right\} = 0$ and
$Var\left\{ {{{\hat m}_{l,\varepsilon }}\left( {{u_l},\bm \beta ^0} \right)\left| {\bm{\mathcal{X}},\bm{\mathcal{ Z}}} \right.} \right\} = \sigma _{nl}^2\left( {{u_l}} \right)$. It can be proved by Theorem 5.4.2 of DeVore and
Lorentz \cite{DL93} and Berstein's inequality in Bosq \cite{B98} that for large enough $n$ and with
probability 1, there are constants $0 < c_R \leq C_R < \infty$ such that for any $\bm \theta  \in {\mathbb{R}^{d{J_n}}}$,
\begin{gather}
{c_R}{J_n}\left\| \bm \theta  \right\|_2^2 \le {\bm \theta ^T}{\bm \Omega _n}{(\bm \beta ^0)^{ - 1}}\bm \theta  \le {C_R}{J_n}\left\| \bm \theta  \right\|_2^2,
  \tag{A.12}\label{A.12}
\end{gather}
and by the above result and Demko \cite{D86}, we obtain
\[{\left\| {{\bm \Omega _n}{{(\bm \beta ^0)}^{ - 1}}} \right\|_\infty } = {O_p}\left( {{J_n}} \right).\]
Thus, by (\ref{A.12}), for large enough $n$ and with probability 1,
\[\begin{array}{l}
n^{-1}{{{c'_R}{J_n}}} \le \mathop {\inf }\limits_{{u_l} \in [0,1]} Var\left\{ {{{\hat m}_{l,\varepsilon }}\left( {{u_l},\bm \beta ^0} \right)\left| {\bm{\mathcal{X}},\bm{\mathcal{ Z}}} \right.} \right\}\le \mathop {\sup }\limits_{{u_l} \in [0,1]} Var\left\{ {{{\hat m}_{l,\varepsilon }}\left( {{u_l},\bm \beta ^0} \right)\left| {\bm{\mathcal{X}},\bm{\mathcal{ Z}}} \right.} \right\} \le n^{-1}{{{C'_R}{J_n}}}
  \end{array}\]
for some constants $0 < c'_R\leq C'_R< \infty$. By weak law of large numbers, we have for every $1\leq l \leq d$,
\begin{gather}
{{\hat m}_{l,\varepsilon }}\left( {{u_l},\bm \beta ^0} \right) = {O_p}\left( {\sqrt {{{{J_n}} \mathord{\left/
 {\vphantom {{{J_n}} n}} \right.
 \kern-\nulldelimiterspace} n}} } \right)
  \tag{A.13}\label{A.13}
\end{gather}
uniformly for any $u_l\in [0,1]$. It can be proved by Lindeberg-Feller central limit theorem that
\begin{gather}
\sigma _{nl}^{ - 1}\left( {{u_l}} \right){{\hat m}_{l,\varepsilon }}\left( {{u_l},\bm \beta ^0} \right)\mathop  \to \limits^d N\left( {0,1} \right).
  \tag{A.14}\label{A.14}
\end{gather}
Since $\mathop {\sup }\limits_{1 \le s \le {J_n},1 \le l \le d} E\left| {{D_{i,sl}}(\bm \beta _{l}^0)} \right| = O\left( {J_n^{ - 1}} \right)$ and $E\left\{ {{D_{i,sl}}{{(\bm \beta _{l}^0)}^2}} \right\} =O\left( {J_n^{ - 1}} \right) $. It can be proved by Berstein's inequality in Bosq \cite{B98}
that
\[\mathop {\sup }\limits_{1 \le s \le {J_n},1 \le l \le d} \sum\limits_{i = 1}^n {\left| {{D_{i,sl}}(\bm \beta _{l}^0)} \right|}  = O\left( {nJ_n^{ - 1}} \right).\]
Thus, by (\ref{A.1}), we have
\[{\left\| {\sum\limits_{i = 1}^n {f\left(0  \left| {\bm X_i, \bm Z_i} \right.\right)}  \Delta_i {\bm D_i}(\bm \beta ^0)} \right\|_\infty } = {O_p}\left( {nJ_n^{ - 1 - r}} \right).\]
Therefore, for every $1\leq l \leq d$,
\begin{gather}
\begin{array}{l}
 \mathop {\sup }\limits_{{u_l} \in [0,1]} \left| {{{\hat m}_{l,m}}\left( {{u_l},\bm \beta ^0} \right)} \right| \\
  \le \mathop {\sup }\limits_{{u_l} \in [0,1]} {\left\| {\bm e_l^T\mathbb{ B}(\bm u)} \right\|_\infty }{\left\| {{n^{ - 1}}{\bm \Omega _n}{{(\bm \beta ^0)}^{ - 1}}} \right\|_\infty }  {\left\| {\sum\limits_{i = 1}^n {f\left(0\left| {\bm X_i, \bm Z_i} \right.\right)}  \Delta_i {\bm D_i}(\bm \beta ^0)} \right\|_\infty } \\
  = {O_p}\left( {{{{J_n}} \mathord{\left/
 {\vphantom {{{J_n}} n}} \right.
 \kern-\nulldelimiterspace} n}} \right){O_p}\left( {nJ_n^{ - 1 - r}} \right)\\
  = {O_p}\left( {J_n^{ - r}} \right).
 \end{array}
  \tag{A.15}\label{A.15}
\end{gather}
By (\ref{A.1}), (\ref{A.11}), (\ref{A.13}) and (\ref{A.15}), we have
\[\left| {{{\hat m}_l}({u_l},{\bm \beta ^0}) - {m_l}({u_l})} \right| = O\left( {\sqrt {{{{J_n}} \mathord{\left/
 {\vphantom {{{J_n}} n}} \right.
 \kern-\nulldelimiterspace} n}}  + J_n^{ - r}} \right) uniformly for any u_l\in [0,1].\]
According to de Boor (\cite{De01}, page 116), ${{\hat {\dot m}}_l}({u_l},\bm \beta ^0)=\bm B_{q-1}(u_l)^T \bm \Xi \bm{\hat \lambda} _l(\bm \beta^0)$, where $\bm B_{q-1}(u_l)=\{B_{q-1,s}(u_l): 2\leq s \leq J_n\}^T$ is the $(q-1)$-th order B-spline basis and
\[\bm \Xi = (q - 1){\left[ \begin{array}{l}
 \frac{{ - 1}}{{{\xi _{q + 1}} - {\xi _2}}}{\rm{  }}~\frac{1}{{{\xi _{q + 1}} - {\xi _2}}}{\rm{      ~~ ~~0              }} ~~~~~~~~~~~\cdots {\rm{                       ~~~~~~~~~~~~0}} \\
 {\rm{   ~~~~  0         }}~~~~~\frac{{ - 1}}{{{\xi _{q + 2}} - {\xi _3}}}{\rm{  }}~\frac{1}{{{\xi _{q + 2}} - {\xi _3}}}{\rm{      }} ~~~~~~\cdots {\rm{                        ~~~~~~~~~~~~0}} \\
 {\rm{     }} ~~~~\vdots {\rm{                  }} ~~~~~~~~\vdots {\rm{            }} ~~~~~~~~~\ddots {\rm{          }} ~~~~~~~~\ddots {\rm{                       }} ~~~~~~~~~~~~\vdots  \\
 {\rm{   ~~~~  0           ~~~~~   ~~~   0             }} ~~~~~~~~\cdots {\rm{     }}~\frac{{ - 1}}{{{\xi _{N_n + 2q - 1}} - {\xi _{N_n + q}}}}{\rm{  }}~\frac{1}{{{\xi _{N_n + 2q - 1}} - {\xi _{N_n + q}}}} \\
 \end{array} \right]_{({J_n} - 1) \times {J_n}}}.\]
 It is easy to prove that ${\left\| {{\bm \Xi}} \right\|_\infty } = O({J_n})$. Following the similar
reasoning as the proof for ${{\hat {m}}_l}({u_l},{\bm \beta ^0})$, we have $\left| {{{\hat {\dot m}}_l}({u_l},{\bm \beta ^0}) - {{\dot m}_l}({u_l})} \right| = O_p\left( {\sqrt {{{J_n^3} \mathord{\left/
 {\vphantom {{J_n^3} n}} \right.
 \kern-\nulldelimiterspace} n}}  + J_n^{ - r + 1}} \right)$ uniformly for any $u_l\in [0,1]$ by Slutsky's theorem. Result (ii) follows from (\ref{A.11}), (\ref{A.13}), (\ref{A.15}) and condition $nJ_n^{ - \left( {1 + 2r} \right)} = o\left( 1 \right)$. This completes the proof. $\square$\\
\end{proof}

\begin{lemma}\label{lemma:2}
 Under conditions (C1)--(C7), for $r\geq 2$, ${n^{{1 \mathord{\left/
 {\vphantom {1 {(2r + 2)}}} \right.
 \kern-\nulldelimiterspace} {(2r + 2)}}}}  \ll {N_n}  \ll {n^{{1 \mathord{\left/
 {\vphantom {1 4}} \right.
 \kern-\nulldelimiterspace} 4}}}$, we have
\[{{\rm{{\cal R}}}_{\tau nh}}\left( {{\bm \beta _{ - 1}^0}} \right) =  \sum\limits_{i = 1}^n {{\psi _\tau }\left\{ {{\varepsilon _i}} \right\}} \left[ {{{\dot m}_l}\left( {\bm Z_i^T {\bm \beta _{l}^0} } \right){X_{il}}\bm J_l^{0T}{{\bm {\tilde Z}}_i}} \right]_{l = 1}^d + {o_p}\left( {\sqrt {n{p_n}} } \right).\]
\end{lemma}

\begin{proof}  Using conditions (C6) and (C7), similar to Lemma 3(k) of Horowitz \cite{H98}, we obtain
\[\begin{array}{l}
 {{\rm{{\cal R}}}_{\tau nh}}\left( {{\bm \beta _{ - 1}}} \right) = \sum\limits_{i = 1}^n {{\psi _\tau }\left\{ {{Y_i} - \sum\limits_{l = 1}^d {{\bm B}{{\left( {\bm Z_i^T {{\bm \beta _{l}}} } \right)}^T}{{\hat {\bm \lambda} }_l}\left( {{\bm \beta }} \right){X_{il}}} } \right\}}  \\
~~~~~~~~~~~~~~~\times \left[ {\begin{array}{*{20}{c}}
   {\left\{ {{{\hat {\dot m}}_1}\left( {\bm Z_i^T{\bm \beta _{1}},{\bm \beta }} \right){ X_{i1}}\bm J_1^T{\bm Z_i} + \left( {{{\partial \hat {\bm \lambda }{{({\bm \beta })}^T}} \mathord{\left/
 {\vphantom {{\partial \hat {\bm \lambda} {{({\bm \beta })}^T}} {\partial {\bm \beta _{1, - 1}}}}} \right.
 \kern-\nulldelimiterspace} {\partial {\bm \beta _{1, - 1}}}}} \right){\bm D_i}({\bm \beta })} \right\}}  \\
   { \vdots }  \\
   {\left\{ {{{\hat {\dot m}}_d}\left( {\bm Z_i^T{\bm \beta _{d}},{\bm \beta }} \right){X_{id}}\bm J_d^T{\bm Z_i} + \left( {{{\partial \hat {\bm \lambda} {{({\bm \beta })}^T}} \mathord{\left/
 {\vphantom {{\partial \hat {\bm \lambda} {{({\bm \beta})}^T}} {\partial {\bm \beta _{d, - 1}}}}} \right.
 \kern-\nulldelimiterspace} {\partial {\bm \beta _{d, - 1}}}}} \right){\bm D_i}({\bm \beta })} \right\}}  \\
\end{array}} \right] + {o_p}\left( {\sqrt {n{p_n}} } \right). \\
 \end{array}\]
 By Taylor expansion and (\ref{A.11}), we have
\begin{gather}
{\bm D_i}{\left( {\bm \beta ^0} \right)^T}\left\{ {\hat {\bm \lambda} \left( {\bm \beta ^0} \right) - {\bm \lambda ^0}} \right\} = {\bm D_i}{\left( {\bm \beta ^0} \right)^T}\left\{ {{{\hat {\bm \lambda} }_\varepsilon }\left( {\bm \beta ^0} \right) + {{\hat {\bm \lambda } }_m}\left( {\bm \beta ^0} \right)} \right\}\left\{1+o_p(1) \right\},
  \tag{A.16}\label{A.16}
\end{gather}
where
\[{{\hat {\bm\lambda} }_\varepsilon }\left( {\bm\beta ^0} \right) =n^{ - 1} {\bm\Omega _n}{(\bm\beta ^0)^{ - 1}} {\sum\limits_{i = 1}^n {{\psi _\tau }\left\{ {{\varepsilon _i}} \right\}{\bm D_i}\left( {\bm\beta ^0} \right)} }, {{\hat{\bm \lambda} }_m}\left( {\bm\beta ^0} \right) = n^{-1}{\bm\Omega _n}{(\bm\beta ^0)^{ - 1}}
{\sum\limits_{i = 1}^n {f\left( 0\left| {\bm X_i, \bm Z_i} \right.\right)} \Delta_i {\bm D_i}(\bm \beta ^0)},\]
and $\Delta_i$ is given in the proof of Lemma \ref{lemma:1}. Define
$\hat {\bm \delta}_k = \mathop {\arg \min }\limits_{\bm \delta  \in {R^{d{J_n} }}} {\sum\limits_{i = 1}^n f(0\left| {\bm X_i, \bm Z_i} \right.){\left\{ {{Z_{ik}} - {\bm D_i}\left( {\bm \beta ^0} \right)^T{{\bm \delta} }} \right\}} ^2},$
\[\begin{array}{l}
 {\mathbb{P}_n}\left( {{Z_{ik}}} \right) = {\bm D_i}\left( {\bm \beta ^0} \right)^T\hat {\bm \delta}_k \\
 ~~~~~~~~~~={\bm D_i}{\left( {{\bm\beta ^0}} \right)^T}{\left( \sum\limits_{i = 1}^n f(0\left| {\bm X_i, \bm Z_i} \right.){{\bm D_i}({\bm \beta ^0})} {\bm D_i}{({\bm \beta ^0})^T} \right)^{ - 1}}\sum\limits_{i = 1}^n f(0\left| {\bm X_i, \bm Z_i} \right.){{\bm D_i}\left( {{\bm\beta ^0}} \right)Z_{ik}}, \\
  \end{array}\]
${\mathbb{P}_n}\left( {{\bm Z_i }} \right)=\left\{{\mathbb{P}_n}\left( {{Z_{ik}}} \right):1\leq k\leq p_n\right\}^T$ and $\bm D \left( {\bm \beta ^0} \right)=\left( \bm D_1\left( {\bm \beta ^0} \right),...,\bm D_n\left( {\bm \beta ^0} \right) \right)^T$. By (\ref{A.12}) and the weak law of large numbers, it can be proved that
\begin{gather}
\begin{array}{l}
 {\bm D_i}{\left( {\bm \beta ^0} \right)^T}\left\{ {{{\partial {{\hat {\bm \lambda } }_\varepsilon }\left( {\bm \beta ^0} \right)} \mathord{\left/
 {\vphantom {{\partial {{\hat {\bm \lambda} }_\varepsilon }\left( {\bm \beta ^0} \right)} {\partial {\bm \beta} _{l, - 1}^T}}} \right.
 \kern-\nulldelimiterspace} {\partial {\bm \beta }_{l, - 1}^T}}} \right\}\\
  =  - {\bm D_i}{\left( {{\bm \beta ^0}} \right)^T}{\left(\sum\limits_{i = 1}^n f(0\left| {\bm X_i, \bm Z_i} \right.){{\bm D_i}({\bm \beta ^0})} {\bm D_i}{({\bm \beta ^0})^T}\right)^{ - 1}}\sum\limits_{i = 1}^n f(0\left| {\bm X_i, \bm Z_i} \right.){{\bm D_i}\left( {{\bm\beta ^0}} \right){{\dot m}_l}\left( {\bm Z_i^T\bm\beta _l^0} \right){X_{il}}\bm Z_i^T{\bm J_l^0}} \\
  ~~~ + {O_p}\left( {\sqrt {{J_n}{p_n}/n} } \right) \\
  =  - {\mathbb{P}_n}\left( {{{\dot m}_l}\left( {\bm Z_i^T\bm\beta _l^0} \right){X_{il}}\bm Z_i^T{\bm J_l^0}} \right) + {O_p}\left( {\sqrt {{J_n}{p_n}/n} } \right). \\
 \end{array}
   \tag{A.17}\label{A.17}
\end{gather}
By B--spline properties, conditions (C2) and (C4), together with the definition of $\mathbb{P}(\cdot)$ given in subsection 2.2, it can be proved that
\[\begin{array}{l}{\left\| {{\mathbb{P}_n}\left(  {{{\dot m}_l}\left( {\bm Z_i^T\bm\beta _l^0} \right){X_{il}}\bm Z_i^T{\bm J_l^0}} \right) - \mathbb{P}\left(  {{{\dot m}_l}\left( {\bm Z_i^T\bm\beta _l^0} \right){X_{il}}\bm Z_i^T{\bm J_l^0}} \right)} \right\|_\infty }= {O_p}\left( {J_n^{{1 \mathord{\left/
 {\vphantom {1 2}} \right.
 \kern-\nulldelimiterspace} 2}}\sqrt {{{{p_n}} \mathord{\left/
 {\vphantom {{{p_n}} n}} \right.
 \kern-\nulldelimiterspace} n}}  + \sqrt {{p_n}} J_n^{ - 1}} \right).
 \end{array}\]
Thus
\begin{gather}
\begin{array}{l}
 {\bm D_i}{\left( {\bm \beta ^0} \right)^T}\left\{ {{{\partial {{\hat {\bm \lambda } }_\varepsilon }\left( {\bm \beta ^0} \right)} \mathord{\left/
 {\vphantom {{\partial {{\hat {\bm \lambda} }_\varepsilon }\left( {\bm \beta ^0} \right)} {\partial {\bm \beta} _{l, - 1}^T}}} \right.
 \kern-\nulldelimiterspace} {\partial {\bm \beta }_{l, - 1}^T}}} \right\}\\
  =  - {\mathbb{P}}\left( {{{\dot m}_l}\left( {\bm Z_i^T\bm\beta _l^0} \right){X_{il}}\bm Z_i^T{\bm J_l^0}} \right) + {O_p}\left( {J_n^{{1 \mathord{\left/
 {\vphantom {1 2}} \right.
 \kern-\nulldelimiterspace} 2}}\sqrt {{{{p_n}} \mathord{\left/
 {\vphantom {{{p_n}} n}} \right.
 \kern-\nulldelimiterspace} n}}  + \sqrt {{p_n}} J_n^{ - 1}} \right)\\
 =- {{\dot m}_l}\left( {\bm Z_i^T\bm \beta _l^0} \right){X_{il}}{\mathbb{P}}\left( {\bm Z_i^T} \right){\bm J_l^0} + {O_p}\left( {J_n^{{1 \mathord{\left/
 {\vphantom {1 2}} \right.
 \kern-\nulldelimiterspace} 2}}\sqrt {{{{p_n}} \mathord{\left/
 {\vphantom {{{p_n}} n}} \right.
 \kern-\nulldelimiterspace} n}}  + \sqrt {{p_n}} J_n^{ - 1}} \right). \\
 \end{array}
   \tag{A.18}\label{A.18}
\end{gather}
For ${\bm{\hat \lambda }_m}\left( {\bm\beta ^0} \right)$, by (\ref{A.12}), with probability approaching 1,
\begin{gather}
\begin{array}{l}
 {\left\| {{{\hat {\bm \lambda} }_m}\left( {\bm \beta^0} \right)} \right\|_\infty } = {\left\| {{{\left\{ {\bm D{{\left( {\bm \beta ^0} \right)}^T}\bm D\left( {\bm \beta ^0} \right)} \right\}}^{ - 1}}\sum\limits_{i = 1}^n  \Delta_i {\bm D_i}(\bm \beta ^0)} \right\|_\infty } \\
 ~~~~~~~~~~~~~~~~ \le n^{ - 1} {C_R}{J_n}{\left\| \sum\limits_{i = 1}^n  \Delta_i {\bm D_i}(\bm \beta ^0) \right\|_\infty } \\
  ~~~~~~~~~~~~~~~~ \le {C_R}{J_n}\mathop {\sup }\limits_{s,l} {n^{ - 1}}\sum\limits_{i = 1}^n {\left| {{B_{s}}\left( {\bm Z_i^T\bm \beta _{l}^0} \right){X_{il}}} \right|} O\left( {J_n^{ - r}} \right) \\
 ~~~~~~~~~~~~~~~~  = O\left( {J_n^{ - 1}} \right)O\left( {J_n^{ - r + 1}} \right) = O\left( {J_n^{ - r}} \right).
 \end{array}
   \tag{A.19}\label{A.19}
\end{gather}
By (\ref{A.16}), (\ref{A.18}) and (\ref{A.19}), we have
\begin{gather}
\begin{array}{l}
\left\{ {{{\partial \hat {\bm \lambda} {{\left( {\bm\beta ^0} \right)}^T}} \mathord{\left/
 {\vphantom {{\partial \hat {\bm\lambda} {{\left( {\bm\beta ^0} \right)}^T}} {\partial {\bm\beta _{l, - 1}}}}} \right.
 \kern-\nulldelimiterspace} {\partial {\bm\beta _{l, - 1}}}}} \right\}{\bm D_i}\left( {\bm\beta ^0} \right) =  - {{\dot m}_l}\left( {\bm Z_i^T {\bm\beta _{l}^0} } \right){X_{il}}\bm J_l^{0T}{\mathbb{P}}\left( {{\bm Z_i}} \right)  +  {O_p}\left( {J_n^{{1 \mathord{\left/
 {\vphantom {1 2}} \right.
 \kern-\nulldelimiterspace} 2}}\sqrt {{{{p_n}} \mathord{\left/
 {\vphantom {{{p_n}} n}} \right.
 \kern-\nulldelimiterspace} n}}  + \sqrt {{p_n}} J_n^{ - 1}} \right).
  \end{array}
 \tag{A.20}\label{A.20}
\end{gather}
Based on the result of Lemma \ref{lemma:1}, we can obtain
\begin{gather}
{{\hat{ \dot m}}_l}\left( {\bm Z_i^T {\bm\beta _{l}^0}} \right){X_{il}}\bm J_l^{0T}{\bm Z_i} = {{\dot m}_l}\left( {\bm Z_i^T {\bm\beta _{l}^0}} \right){X_{il}}\bm J_l^{0T}{\bm Z_i} + {O_p}\left( {J_n^{{3 \mathord{\left/
 {\vphantom {3 2}} \right.
 \kern-\nulldelimiterspace} 2}}\sqrt {{{{p_n}} \mathord{\left/
 {\vphantom {{{p_n}} n}} \right.
 \kern-\nulldelimiterspace} n}}  + \sqrt {{p_n}} J_n^{ - r + 1}} \right).
  \tag{A.21}\label{A.21}
\end{gather}
 Thus, combine with (\ref{A.20}) and (\ref{A.21}), we have
\begin{gather}
\begin{array}{l}
 {{\hat {\dot m}}_l}\left( {\bm Z_i^T {\bm\beta _{l}^0}} \right){X_{il}}\bm J_l^{0T}{\bm Z_i} + \left\{ {{{\partial \hat {\bm \lambda} {{\left( {\bm \beta ^0} \right)}^T}} \mathord{\left/
 {\vphantom {{\partial \hat {\bm \lambda }{{\left( {\bm \beta ^0} \right)}^T}} {\partial {\bm \beta _{l, - 1}}}}} \right.
 \kern-\nulldelimiterspace} {\partial {\bm \beta _{l, - 1}}}}} \right\}{\bm D_i}\left( {\bm \beta ^0} \right) \\
  = {{\dot m}_l}\left( {\bm Z_i^T {\bm\beta _{l}^0}} \right){X_{il}}\bm J_l^{0T}{{\tilde{\bm  Z}}_i} + {O_p}\left( {J_n^{{3 \mathord{\left/
 {\vphantom {3 2}} \right.
 \kern-\nulldelimiterspace} 2}}\sqrt {{{{p_n}} \mathord{\left/
 {\vphantom {{{p_n}} n}} \right.
 \kern-\nulldelimiterspace} n}}  + \sqrt {{p_n}} J_n^{ - r + 1}} + \sqrt {{p_n}} J_n^{ -  1}\right).
 \end{array}
    \tag{A.22}\label{A.22}
\end{gather}
By (\ref{A.22}) and Lemma \ref{lemma:1}, for $r\geq 2$, ${n^{{1 \mathord{\left/
 {\vphantom {1 {(2r + 2)}}} \right.
 \kern-\nulldelimiterspace} {(2r + 2)}}}}  \ll {J_n}  \ll {n^{{1 \mathord{\left/
 {\vphantom {1 4}} \right.
 \kern-\nulldelimiterspace} 4}}}$,
\[\begin{array}{l}
 {{\rm{{\cal R}}}_{\tau nh}}\left( {{\bm \beta _{ - 1}^0}} \right)\\
  =  \sum\limits_{i = 1}^n {\left\{ {{\psi _\tau }\left\{ {{\varepsilon _i}} \right\} + {O_p}\left( {\sqrt {{{{J_n}} \mathord{\left/
 {\vphantom {{{J_n}} n}} \right.
 \kern-\nulldelimiterspace} n}}  + J_n^{ - r}} \right)} \right\}}  \\
  ~~~~\times \left\{ {\left[ {{{\dot m}_l}\left( {\bm Z_i^T{\bm \beta _{l}^0} } \right){X_{il}}\bm J_l^{0T}{{\tilde{\bm  Z}}_i}} \right]_{l = 1}^d + {O_p}\left( {J_n^{{3 \mathord{\left/
 {\vphantom {3 2}} \right.
 \kern-\nulldelimiterspace} 2}}\sqrt {{{{p_n}} \mathord{\left/
 {\vphantom {{{p_n}} n}} \right.
 \kern-\nulldelimiterspace} n}}  + \sqrt {{p_n}} J_n^{ - r + 1}} +\sqrt {{p_n}} J_n^{ -  1}\right)} \right\} \\
 =  \sum\limits_{i = 1}^n {{\psi _\tau }\left\{ {{\varepsilon _i}} \right\}\left[ {{{\dot m}_l}\left( {\bm Z_i^T{\bm \beta _{l}^0}} \right){X_{il}}\bm J_l^{0T}{{\tilde {\bm Z}}_i}} \right]_{l = 1}^d}  \\
~~ + {O_p}\left( {J_n^{{3 \mathord{\left/
 {\vphantom {3 2}} \right.
 \kern-\nulldelimiterspace} 2}}\sqrt {{{{p_n}} \mathord{\left/
 {\vphantom {{{p_n}} n}} \right.
 \kern-\nulldelimiterspace} n}}  + \sqrt {{p_n}} J_n^{ - r + 1}}+ \sqrt {{p_n}} J_n^{ -  1} \right)  \left\{ {O_p}\left( n^{1/2}\right)+  {O_p}\left( {n\sqrt {{{{J_n}} \mathord{\left/
 {\vphantom {{{J_n}} n}} \right.
 \kern-\nulldelimiterspace} n}}  + nJ_n^{ - r}} \right) \right\} \\
 =   \sum\limits_{i = 1}^n {{\psi _\tau }\left\{ {{\varepsilon _i}} \right\}\left[ {{{\dot m}_l}\left( {\bm Z_i^T{\bm \beta _{l}^0}} \right){X_{il}}\bm J_l^{0T}{{\tilde {\bm Z}}_i}} \right]_{l = 1}^d}  + {o_p}\left( {\sqrt {n{p_n}} } \right). \\
 \end{array}\]
\end{proof}

\noindent \textbf{Proof of  Theorem~2.1} To prove (i), it suffices to show that $\forall \delta  > 0$, there exists a $\triangle>0$ such that for all $n$ sufficiently large,
\begin{gather}
P\left\{ {\mathop {\sup }\limits_{\left\| {{\bm \beta _{ - 1}} - \bm \beta _{ - 1}^0} \right\|_2 = \Delta\sqrt {{{{p_n}} \mathord{\left/
 {\vphantom {{{p_n}} n}} \right.
 \kern-\nulldelimiterspace} n}} } {{\left( {{\bm \beta _{ - 1}} - \bm \beta _{- 1}^0} \right)}^T}{{\rm{{\cal R}}}_{\tau nh}}\left( {{\bm \beta _{ - 1}}} \right) < 0} \right\} \ge 1 - \delta .
   \tag{A.23}\label{A.23}
\end{gather}
 This requires us to calculate the sign of $ {{\left( {{\bm \beta _{ - 1}} - \bm \beta _{- 1}^0} \right)}^T} {{\rm{{\cal R}}}_{\tau nh}}\left( {{\bm \beta _{ - 1}}} \right)$ on $\left\{ \bm \beta_{-1}: \left\| {{\bm \beta _{ - 1}} - \bm \beta _{ - 1}^0} \right\|_2 = \Delta\sqrt {{{{p_n}} \mathord{\left/
 {\vphantom {{{p_n}} n}} \right.
 \kern-\nulldelimiterspace} n}} \right\}$. Note that
\[\begin{array}{l}
 {\left( {{\bm \beta _{- 1}} - \bm \beta _{ - 1}^0} \right)^T}{{\rm{{\cal R}}}_{\tau nh}}\left( {{\bm \beta _{- 1}}} \right) \\
  = {\left( {{\bm \beta _{- 1}} - \bm \beta _{ - 1}^0} \right)^T}{{\rm{{\cal R}}}_{\tau nh}}\left( {\bm \beta _{ - 1}^0} \right) + {\left( {{\bm \beta _{ - 1}} - \bm \beta _{ - 1}^0} \right)^T}\frac{\partial }{{\partial {\bm \beta} _{ - 1}^T}}{{\rm{{\cal R}}}_{\tau nh}}\left( {\bm \beta _{ - 1}^*} \right)\left( {{\bm \beta _{ - 1}} - \bm \beta _{ - 1}^0} \right) \\
  = {I_{n1}} + {I_{n2}}
 \end{array}\]
where $\bm \beta _{ - 1}^*$ lies between $\bm \beta _{- 1}^0$ and $\bm \beta _{ - 1}$. We first consider $I_{n1}$. By Lemma \ref{lemma:2}, we have
 $ {{\rm{{\cal R}}}_{\tau nh}}\left( {{\bm \beta _{ - 1}^0}} \right)={\bm S_n}\left( {{\bm \beta _{ - 1}^0}} \right)+{o_p}\left( {\sqrt {n{p_n}} } \right)$, where
 \[{\bm S_n}\left( {\bm\beta _{ - 1}^0} \right) = \sum\limits_{i = 1}^n {{\psi _\tau }\left\{ {{\varepsilon _i}} \right\}} \left[ {{{\dot m}_l}\left( {\bm Z_i^T\bm\beta _l^0} \right){X_{il}}\bm J_l^{0T}{{ \bm{\tilde{Z}}}_i}} \right]_{l = 1}^d.\]
 For any $\bm \beta _{ - 1}$ such that $\left\| {{\bm \beta _{ - 1}} - \bm \beta _{ - 1}^0} \right\| = \Delta \sqrt {{{{p_n}} \mathord{\left/
 {\vphantom {{{p_n}} n}} \right.
 \kern-\nulldelimiterspace} n}}  $, we have $\left| {{I_{n1}}} \right| \le \Delta \sqrt {{{{p_n}} \mathord{\left/
 {\vphantom {{{p_n}} n}} \right.
 \kern-\nulldelimiterspace} n}} \left\| {{\bm S_n}\left( {\bm \beta _{- 1}^0} \right)} \right\|_2+o_p(p_n)$. Note that $E\{ {\bm S_n}\left( {{\bm \beta _{ - 1}^0}} \right) \}=\bm 0$, and
\[\begin{array}{l}
 E\left[ {{{\left\| {{\bm S_n}\left( {\bm \beta _{ - 1}^0} \right)} \right\|_2^2}}} \right] = E\left[ {\sum\limits_{i = 1}^n {\psi _{\tau }^2\left\{ {{\varepsilon _i}} \right\}} tr \left({{{\left\{ {\left[ {{{\dot m}_l}\left( {\bm Z_i^T{\bm\beta _{l}^0} } \right){X_{il}}\bm J_l^{0T}{{\tilde {\bm Z}}_i}} \right]_{l = 1}^d} \right\}}^{ \otimes 2}}}  \right)} \right] \\
~~~~~~~~~~~~~~~~~~~~~~~  \le Ctr\left( \sum\limits_{i = 1}^n {{{\left\{ {\left[ {{{\dot m}_l}\left( {\bm Z_i^T {\bm \beta _{l}^0} } \right){X_{il}}\bm J_l^{0T}{{\tilde {\bm Z}}_i}} \right]_{l = 1}^d} \right\}}^{ \otimes 2}}} \right)\\
~~~~~~~~~~~~~~~~~~~~~~~= {O}\left( {n{p_n}} \right) \\
 \end{array}\]
by conditions (C3) and (C4). Thus, ${\bm S_n}\left( {\bm \beta _{ - 1}^0} \right) = {O_p}\left( {\sqrt {n{p_n}} } \right)$ and $\left| {{I_{n1}}} \right| \le \Delta {O_p}\left( {{p_n}} \right)$. Next, by Taylor expansion and Lemma \ref{lemma:2},
\[\begin{array}{l}
 {I_{n2}} =  - {\left( {{\bm \beta _{ - 1}} - \bm \beta _{ - 1}^0} \right)^T}\left[ {{\sum\limits_{i = 1}^nf(0\left| {\bm X_i, \bm Z_i} \right.){{\left\{ {\left[ {{{\dot m}_l}\left( {\bm Z_i^T {\bm \beta _{l}^0}} \right){X_{il}}\bm J_l^{0T}{{\tilde {\bm Z}}_i}} \right]_{l = 1}^d} \right\}}^{ \otimes 2}}} } \right] \\
~~~~~~~~  \times \left( {{\bm \beta _{ - 1}} - \bm \beta _{ - 1}^0} \right)\left\{ {1 + {o_p}\left( 1 \right)} \right\} \\
~~~~  \le  - C \rho_{min}\left( {\sum\limits_{i = 1}^n {{{\left\{ {\left[ {{{\dot m}_l}\left( {\bm Z_i^T {\bm \beta _{l}^0}} \right){X_{il}}\bm J_l^{0T}{{\tilde {\bm Z}}_i}} \right]_{l = 1}^d} \right\}}^{ \otimes 2}}} } \right)\left\| {{\bm \beta _{- 1}} - \bm \beta _{ - 1}^0} \right\|_2^2 \\
~~~~  \le  - C{\Delta ^2}{p_n}. \\
 \end{array}\]
For sufficiently large $\Delta$, $ {\left( {{\bm \beta _{ - 1}} - \bm \beta _{ - 1}^0} \right)^T}{{\rm{{\cal R}}}_{\tau nh}}\left( {{\bm \beta _{ - 1}}} \right) $ on $\left\{ \bm \beta_{-1}: \left\| {{\bm \beta _{ - 1}} - \bm \beta _{ - 1}^0} \right\|_2 = \Delta\sqrt {{{{p_n}} \mathord{\left/
 {\vphantom {{{p_n}} n}} \right.
 \kern-\nulldelimiterspace} n}} \right\}$ is dominated by $I_{n2}$, which is large and negative for all sufficiently large $n$. Thus, (\ref{A.23}) holds.

Next, we need prove (ii). Let ${\Gamma _{ni}} =\bm  e _n^T\bm M_n^{{{ - 1} \mathord{\left/
 {\vphantom {{ - 1} 2}} \right.
 \kern-\nulldelimiterspace} 2}}\left(\bm \beta_{-1}^{0}\right){\psi _{\tau }}\left\{ {{\varepsilon _i}} \right\}\left[ {{{\dot m}_l}\left( {\bm Z_i^T {\bm \beta _{l}^0}} \right){X_{il}}\bm J_l^{0T}{{\tilde {\bm Z}}_i}} \right]_{l = 1}^d$ and $\bm e _n^T\bm M_n^{{{ - 1} \mathord{\left/
 {\vphantom {{ - 1} 2}} \right.
 \kern-\nulldelimiterspace} 2}}\left( {\bm \beta _{ - 1}^0} \right){\bm S_n}\left( {\bm \beta _{ - 1}^0} \right) = \sum\limits_{i = 1}^n {{\ \Gamma _{ni}}} $. Since $E({\psi _{\tau }}\left\{ {{\varepsilon _i}} \right\})=0$ and $ Cov\left( \bm S_n \left(\bm \beta_{-1}^{0}\right) \right)= \tau (1-\tau)\bm M_n \left(\bm \beta_{-1}^{0}\right) $, we have $E\left(\bm e _n^T\bm M_n^{{{ - 1} \mathord{\left/
 {\vphantom {{ - 1} 2}} \right.
 \kern-\nulldelimiterspace} 2}}\left( {\bm \beta _{ - 1}^0} \right){\bm S_n}\left( {\bm \beta _{ - 1}^0} \right)\right) = 0$ and $Var\left( \bm e _n^T\bm M_n^{{{ - 1} \mathord{\left/
 {\vphantom {{ - 1} 2}} \right.
 \kern-\nulldelimiterspace} 2}}\left( {\bm \beta _{ - 1}^0} \right){\bm S_n}\left( {\bm \beta _{ - 1}^0} \right)\right) = \tau (1-\tau)$. To establish the asymptotic normality, it suffices to check the Lindberg condition, that is, $\forall \delta  > 0$, $\sum\limits_{i = 1}^n {E\left[ { \Gamma _{ni}^2I\left( {\left| {{\Gamma _{ni}}} \right| > \delta } \right)} \right]} \rightarrow 0$.
For any $\delta>0$,
\begin{gather}
\sum\limits_{i = 1}^n {E\left[ {\Gamma _{ni}^2I\left( {\left| {{\Gamma _{ni}}} \right| > \delta } \right)} \right]}  \le n{\left\{ {E\left( {\Gamma _{ni}^4} \right)} \right\}^{{1 \mathord{\left/
 {\vphantom {1 2}} \right.
 \kern-\nulldelimiterspace} 2}}}{\left\{ {P\left( {\left| {{\Gamma _{ni}}} \right| > \delta } \right)} \right\}^{{1 \mathord{\left/
 {\vphantom {1 2}} \right.
 \kern-\nulldelimiterspace} 2}}},
    \tag{A.24}\label{A.24}
\end{gather}
where
\[\begin{array}{l}
 \Gamma _{ni}^2 = \left( {\bm e_n^T\bm M_n^{{{ - 1} \mathord{\left/
 {\vphantom {{ - 1} 2}} \right.
 \kern-\nulldelimiterspace} 2}}\left( {\bm \beta _{ - 1}^0} \right){{\left\{ {\left[ {{{\dot m}_l}\left( {\bm Z_i^T {\bm \beta _{l}^0}} \right){X_{il}}\bm J_l^{0T}{{\tilde {\bm Z}}_i}} \right]_{l = 1}^d} \right\}}^{ \otimes 2}}\bm M_n^{{{ - 1} \mathord{\left/
 {\vphantom {{ - 1} 2}} \right.
 \kern-\nulldelimiterspace} 2}}\left( {\bm \beta _{ - 1}^0} \right){\bm e_n}} \right)\left( {{\psi _{\tau }}\left\{ {{\varepsilon _i}} \right\}} \right)^2 \\
 ~~~~ \triangleq{\gamma _{ni}}\left( {{\psi _{\tau }}\left\{ {{\varepsilon _i}} \right\}} \right)^2, \\
 \end{array}\]
 where $ \gamma _{ni}={\bm e_n^T\bm M_n^{{{ - 1} \mathord{\left/
 {\vphantom {{ - 1} 2}} \right.
 \kern-\nulldelimiterspace} 2}}\left( {\bm \beta _{ - 1}^0} \right){{\left\{ {\left[ {{{\dot m}_l}\left( {\bm Z_i^T {\bm \beta _{l}^0}} \right){X_{il}}\bm J_l^{0T}{{\tilde {\bm Z}}_i}} \right]_{l = 1}^d} \right\}}^{ \otimes 2}}\bm M_n^{{{ - 1} \mathord{\left/
 {\vphantom {{ - 1} 2}} \right.
 \kern-\nulldelimiterspace} 2}}\left( {\bm \beta _{ - 1}^0} \right){\bm e_n}}$. Next, we will prove that $\mathop {\max }\limits_{1 \le i \le n} {\gamma _{ni}} \to 0$ as $n\rightarrow\infty$. Note that
  \[{\gamma _{ni}} \le {\rho _{\max }}\left( {{{\left\{ {\left[ {{{\dot m}_l}\left( {\bm Z_i^T {\bm \beta _{l}^0}} \right){X_{il}}\bm J_l^{0T}{{\tilde{\bm  Z}}_i}} \right]_{l = 1}^d} \right\}}^{ \otimes 2}}} \right)\rho _{\min }^{ - 1}\left( {{\bm M_n}\left( {\bm \beta _{- 1}^0} \right)} \right).\]
Therefore, we have
\[\begin{array}{l}
{\gamma _{ni}} \le \frac{{{\rho _{\max }}\left( {{{\left\{ {\left[ {{{\dot m}_l}\left( {\bm Z_i^T {\bm \beta _{l}^0}} \right){X_{il}}\bm J_l^{0T}{{\tilde {\bm Z}}_i}} \right]_{l = 1}^d} \right\}}^{ \otimes 2}}} \right)}}{{ {{\rho _{\min }}\left(\sum\limits_{i = 1}^n {{{\left\{ {\left[ {{{\dot m}_l}\left( {\bm Z_i^T {\bm \beta _{l}^0}} \right){X_{il}}\bm J_l^{0T}{{\tilde {\bm Z}}_i}} \right]_{l = 1}^d} \right\}}^{ \otimes 2}}} \right)} }}
\le \frac{{tr\left( {{{\left\{ {\left[ {{{\dot m}_l}\left( {\bm Z_i^T {\bm \beta _{l}^0}} \right){X_{il}}\bm J_l^{0T}{{{\bm {\tilde Z}}}_i}} \right]_{l = 1}^d} \right\}}^{ \otimes 2}}} \right)}}{{ {{\rho _{\min }}\left(\sum\limits_{i = 1}^n {{{\left\{ {\left[ {{{\dot m}_l}\left( {\bm Z_i^T {\bm \beta _{l}^0}} \right){X_{il}}\bm J_l^{0T}{{\tilde {\bm Z}}_i}} \right]_{l = 1}^d} \right\}}^{ \otimes 2}}} \right)} }}
 \end{array},\]
it follows that $\mathop {\max }\limits_{1 \le i \le n} {\gamma _{ni}} \le O\left( {{{{p_n}} \mathord{\left/
 {\vphantom {{{p_n}} n}} \right.
 \kern-\nulldelimiterspace} n}} \right) $. Note that $\left( {{\psi _{\tau }}\left\{ {{\varepsilon _i}} \right\}} \right)^2$ is uniformly bounded, we have by Cauchy inequality
\begin{gather}
P\left( {\left| {{\Gamma _{ni}}} \right| > \delta } \right) \le \frac{{E\left( {\Gamma _{ni}^2} \right)}}{{{\delta ^2}}} = {O}\left( {\frac{{{p_n}}}{n}} \right),
 \tag{A.25}\label{A.25}
\end{gather}
and
\begin{gather}
\begin{array}{l}
 E\left( {\Gamma _{ni}^4} \right) \le C\rho _{\max }^2\left( {{{\left\{ {\left[ {{{\dot m}_l}\left( {\bm Z_i^T\bm \beta _l^0} \right){X_{il}}\bm  J_l^{0T}{\bm {\tilde Z}_i}} \right]_{l = 1}^d} \right\}}^{ \otimes 2}}} \right)\rho _{\min }^{ - 2}\left( {{\bm  M_n}\left( {\bm \beta _{ - 1}^0} \right)} \right) \\
  ~~~~~~~~~\le C{\left\{ {tr\left( {{{\left\{ {\left[ {{{\dot m}_l}\left( {\bm Z_i^T\bm \beta _l^0} \right){X_{il}}\bm J_l^{0T}{\bm {\tilde Z}_i}} \right]_{l = 1}^d} \right\}}^{ \otimes 2}}} \right)} \right\}^2}\rho _{\min }^{ - 2}\left( {{\bm M_n}\left( {\bm \beta _{ - 1}^0} \right)} \right) \\
    ~~~~~~~~~= {O}\left( {\frac{{p_n^2}}{{{n^2}}}} \right). \\
 \end{array}
 \tag{A.26}\label{A.26}
\end{gather}
Combining (\ref{A.24})--(\ref{A.26}) and condition $n^{-1}p_n^{3}=o(1)$, we have
\[\sum\limits_{i = 1}^n {E\left[ {\Gamma _{ni}^2I\left( {\left| {{\Gamma _{ni}}} \right| > \delta } \right)} \right]}  = O\left(  \frac{{p_n}^{3/2}}{n^{1/2}} \right) = o\left( 1 \right).\]
Therefore, the Lindberg condition has been verified, we have
\[\bm e _n^T\bm M_n^{{{ - 1} \mathord{\left/
 {\vphantom {{ - 1} 2}} \right.
 \kern-\nulldelimiterspace} 2}}\left( {\bm \beta _{ - 1}^0} \right){\bm S_n}\left( {\bm \beta _{ - 1}^0} \right)\mathop  \to \limits^d N\left( {0,\tau(1-\tau)} \right).\]
  Thus by Taylor series expansion, Lemma \ref{lemma:2} and the above result, we have
 \[\bm e _n^T\bm M_n^{{{ - 1} \mathord{\left/
 {\vphantom {{ - 1} 2}} \right.
 \kern-\nulldelimiterspace} 2}}\left( {\bm \beta _{ - 1}^0} \right){\bm S_n}\left( {\bm \beta _{ - 1}^0} \right)=  \bm e_n^T\bm  M_n^{-1/2}\left( {\bm \beta _{ - 1}^0} \right)\bm H_n\left( {\bm \beta _{ - 1}^0} \right)\left( {{\bm \beta _{ - 1}} - \bm \beta _{ - 1}^0} \right)\left\{ {1 + {o_p}\left( 1 \right)} \right\}.\]
  So we complete the proof. $\square$\\

\noindent \textbf{Proof of  Theorem~2.2} Since $\left\| {{{\bm{\hat \beta} }_{ - 1}} - \bm{\beta} _{ - 1}^0} \right\|_2 = {O_p}\left( \sqrt {{p_n}/n}  \right)$, Theorem 2.2 follows from this result and Lemma \ref{lemma:1}. $\square$\\

\noindent \textbf{Proof of  Theorem~3.1} Let $\alpha_n={O}\left( \sqrt {{p_n}} \left( {{n^{{{ - 1} \mathord{\left/
 {\vphantom {{ - 1} 2}} \right.
 \kern-\nulldelimiterspace} 2}}} + {a_n}} \right) \right)$. It suffices to show that $\forall \delta  > 0$, there exists a $\Delta>0$ such that for all $n$ sufficiently large,
\begin{gather}
P\left\{ {\mathop {\sup }\limits_{\left\| {{\bm \beta _{ - 1}} - \bm\beta _{ - 1}^0} \right\|_2 = \Delta \alpha_n } {{\left( {{\bm \beta _{ - 1}} - \bm \beta _{- 1}^0} \right)}^T}   [{\mathcal{R}_{\tau nh}}\left( {\bm\beta_{-1} } \right)-n{\bm b_{{\alpha _1}}}\left( {{\bm \beta _{ - 1}}} \right)]   < 0} \right\} \ge 1 - \delta .
 \tag{A.27}\label{A.27}
\end{gather}
From the proof of  Theorem 2.1 (i), we can obtain ${\left( {{\bm \beta _{- 1}} - \bm \beta _{ - 1}^0} \right)^T}{{\rm{{\cal R}}}_{\tau nh}}\left( {{\bm \beta _{- 1}}} \right)=O_p\left(-{ n\alpha_n^2}  {\Delta ^2} \right)$. Define
\[\bm {b_{\alpha _1}}(\bm\beta _{ - 1}^0) = \left[ {{\dot{p}_{{\alpha _1}}}\left( {\left| {\beta _{12}^{0}} \right|} \right){\rm{sgn}}\left( { {\beta _{12}^{0}} } \right),...,{\dot{p}_{{\alpha _1}}}\left( {\left| {\beta _{1p_n}^{0}} \right|} \right){\rm{sgn}}\left( { {\beta _{1p_n}^{0}}} \right),...,{\dot{p}_{{\alpha _1}}}\left( {\left| {\beta _{dp_n}^{0}} \right|} \right){\rm{sgn}}\left( { {\beta _{dp_n}^{0}} } \right)} \right]\]
 and
 \[\bm \Sigma _{\alpha _1} (\bm\beta _{ - 1}^0)= diag \left\{ {{ \ddot{p}_{{\alpha _1}}}\left( {\left| {\beta _{12}^{0}} \right|} \right),...,{\ddot{p}_{{\alpha _1}}}\left( {\left| {\beta _{1p_n}^{0}} \right|} \right),...,{\ddot{p}_{{\alpha _1}}}\left( {\left| {\beta _{dp_n}^{0}} \right|} \right)} \right\}.\]
  By the definition of ${\bm b_{{\alpha _1}}}\left( {\bm\beta _{ - 1}} \right)$ given in section 3, taking Taylor's explanation for ${\bm b_{{\alpha _1}}}\left( {\bm \beta _{ - 1}} \right)$ at ${\bm\beta _{ - 1}^0} $ gives
\[{\bm b_{{\alpha _1}}}\left( {{\bm\beta _{ - 1}}} \right) = {\bm b_{{\alpha _1}}}\left( {\bm\beta _{ - 1}^0} \right) + {\bm\Sigma _{{\alpha _1}}}\left( {\bm\beta _{ - 1}^0} \right)\left\{ {1 + o\left( 1 \right)} \right\}\left( {{\bm\beta _{ - 1}} -\bm \beta _{ - 1}^0} \right).\]
Moreover, by conditions (C9) and (C10), we have
\[\begin{array}{l}
 \left| {{{\left( {{\bm\beta _{ - 1}} -\bm \beta _{ - 1}^0} \right)}^T}n{\bm b_{{\alpha _1}}}\left( {{\bm\beta _{ - 1}}} \right)} \right| \\
  \le \left| {{{\left( {{\bm\beta _{ - 1}} - \bm\beta _{ - 1}^0} \right)}^T}n{\bm b_{{\alpha _1}}}\left( {\bm\beta _{ - 1}^0} \right)} \right|+ \left| {{{\left( {{\bm\beta _{ - 1}} - \bm\beta _{ - 1}^0} \right)}^T}n{\bm\Sigma _{{\alpha _1}}}\left( {\bm\beta _{ - 1}^0} \right)\left\{ {1 + o\left( 1 \right)} \right\}\left( {{\bm\beta _{ - 1}} - \bm\beta _{ - 1}^0} \right)} \right| \\
  \le n\sqrt {d {({p_n}-1) } } \mathop {\max }\limits_{l = 1,...,d,j = 2,...,{p_n}} {{\dot p}_{{\alpha _1}}}\left( {\left| {\beta _{lj}^0} \right|} \right){\left\| {{\bm\beta _{ - 1}} - \bm\beta _{ - 1}^0} \right\|_2} \\
  ~~+ n\mathop {\max }\limits_{l = 1,...,d,j = 2,...,{p_n}} {{\ddot p}_{{\alpha _1}}}\left( {\left| {\beta _{lj}^0} \right|} \right)\left\| {{\bm\beta _{ - 1}} - \bm\beta _{ - 1}^0} \right\|_2^2 \\
  \le n{a_n}\sqrt {d{({p_n}-1) } } {\left\| {{\bm\beta _{ - 1}} - \bm\beta _{ - 1}^0} \right\|_2} + n{b_n}\left\| {{\bm\beta _{ - 1}} - \bm\beta _{ - 1}^0} \right\|_2^2 \\
  = {O_p}\left( {n\alpha _n^2\Delta } \right) + {o_p}\left( {n\alpha _n^2{\Delta ^2}} \right). \\
 \end{array}\]
Therefore, for sufficiently large $\Delta$, $ {\left( {{\bm \beta _{ - 1}} - \bm \beta _{ - 1}^0} \right)^T}[{\mathcal{R}_{\tau nh}}\left( {\bm\beta_{-1} } \right)-n{\bm b_{{\alpha _1}}}\left( {{\bm \beta _{ - 1}}} \right)]   $ is dominated by ${\left( {{\bm \beta _{- 1}} - \bm \beta _{ - 1}^0} \right)^T}{{\rm{{\cal R}}}_{\tau nh}}\left( {{\bm \beta _{- 1}}} \right)$, which is large and negative for all sufficiently large $n$. Thus, (\ref{A.27}) holds. $\square$\\

\noindent \textbf{Proof of  Theorem 3.2}. From Theorem 3.1 for a sufficiently large $C$, ${\bm{\bar \beta }_{\alpha_1, - 1}}$ lies in the ball $\left\{ {\bm \beta _{ - 1}^0 + {\alpha _n}\bm u:{{\left\| \bm u \right\|}_2} \le C} \right\}$ with probability converging to 1, where ${\alpha _n} = \sqrt {{p_n}} \left( {{n^{{{ - 1} \mathord{\left/
 {\vphantom {{ - 1} 2}} \right.
 \kern-\nulldelimiterspace} 2}}} + {a_n}} \right)$.
From the proof of  Theorem 2.1, we have ${\left\| {{{\rm{{\cal R}}}_{\tau nh}}\left( {{\bm \beta _{ - 1}}} \right)} \right\|_2} = {O_p}\left( {\sqrt {n{p_n}} } \right)$. Using ${{\sqrt {{{{p_n}} \mathord{\left/
 {\vphantom {{{p_n}} n}} \right.
 \kern-\nulldelimiterspace} n}} } \mathord{\left/
 {\vphantom {{\sqrt {{{{p_n}} \mathord{\left/
 {\vphantom {{{p_n}} n}} \right.
 \kern-\nulldelimiterspace} n}} } {{\alpha _1}}}} \right.
 \kern-\nulldelimiterspace} {{\alpha _1}}} \to 0$ and condition (C8), for $l=1,...,d,$ $j=s_l+1,...,p_n$, we have
\begin{gather}
{{\rm{{\cal R}}}_{\tau nh}}\left( {{\beta _{lj, - 1}}} \right) - n{b_{{\alpha _1}}}\left( {{\beta _{lj, - 1}}} \right) = n{\alpha _1}\left\{ {{O_p}\left( {\frac{{\sqrt {{{{p_n}} \mathord{\left/
 {\vphantom {{{p_n}} n}} \right.
 \kern-\nulldelimiterspace} n}} }}{{{\alpha _1}}}} \right) - \frac{{{{\dot p}_{{\alpha _1}}}\left( {\left| {{\beta _{lj, - 1}}} \right|} \right)}}{{{\alpha _1}}}{\mathop{\rm sgn}} \left( {{\beta _{lj, - 1}}} \right)} \right\}.
   \tag{A.28}\label{A.28}
\end{gather}
Obviously the sign of $\beta _{lj, - 1}$ determines the sign of ${{\rm{{\cal R}}}_{\tau nh}}\left( {{\bm \beta _{ - 1}}} \right) - n{b_{{\alpha _1}}}\left( {{\bm \beta _{ - 1}}} \right)$. Hence, (\ref{A.28}) implies that
\[{{\rm{{\cal R}}}_{\tau nh}}\left( {{\beta _{lj, - 1}}} \right) - n{b_{{\alpha _1}}}\left( {{\beta _{lj, - 1}}} \right) = \left\{ \begin{array}{l}
  < 0,0 < {\beta _{lj, - 1}} < \epsilon_n \\
  > 0, -\epsilon_n  < {\beta _{lj, - 1}} < 0, \\
 \end{array} \right.\]
where $l=1,...,d,$ $j=s_l+1,...,p_n$. This implies ${{\bar \beta }_{\alpha _1lj, - 1}}=0$ with probability converging to 1 for $l=1,...,d,$ $j=s_l+1,...,p_n$. This completes the proof of (i). Next we need to prove (ii). As shown in Theorem 3.1, $\bm {\bar {\beta}}_{\alpha _1,-1}$ is $\sqrt {{p_n}} \left( {{n^{{{ - 1} \mathord{\left/
 {\vphantom {{ - 1} 2}} \right.
 \kern-\nulldelimiterspace} 2}}} + {a_n}} \right)$ consistent. By the proof of (i), each component of $\bar{\bm \beta}_{\alpha _1,-1}^{(1)}$ stays away from zero for a sufficiently large sample size $n$. At the same time $\bar{\bm \beta}_{\alpha _1,-1}^{(2)}=0$ with probability tending to 1. Thus, $\bm{\bar\beta}_{\alpha _1,-1}^{(1)}$ satisfies
\[n^{-1}{{\rm{{\cal R}}}_{\tau nh}}\left( {\bm \beta _{ - 1}^{(1)}} \right) - {\bm b_{{\alpha _1}}}\left( {\bm \beta _{ - 1}^{(1)}} \right) = 0,\]
where ${\bm b_{{\alpha _1}}}\left( {\bm \beta _{ - 1}^{(1)}} \right)=\left({{\dot p}_{{\alpha _1}}}\left( {\left| {{\beta _{lj}}} \right|} \right){\mathop{\rm sgn}} \left( {{\beta _{lj}}} \right)\right)_{l=1,j=2}^{d,s_l}$ is a $(s_1+...+s_l-d)\times 1$ vector. Applying the Taylor expansion to ${\bm b_{{\alpha _1}}}\left( {\bm \beta _{ - 1}^{(1)}} \right)$, we get that
\[{{\dot p}_{{\alpha _1}}}\left( {\left| {{\beta _{lj}}} \right|} \right){\mathop{\rm sgn}}\left( {{\beta _{lj}}} \right) = {{\dot p}_{{\alpha _1}}}\left( {\left| {\beta _{lj}^0} \right|} \right){\mathop{\rm sgn}}\left( {\beta _{lj}^0} \right) + {{\ddot p}_{{\alpha _1}}}\left( {\left| {\beta _{lj}^0} \right|} \right)\left\{ {1 + o\left( 1 \right)} \right\}\left( {{\beta _{lj}} - \beta _{lj}^0} \right)\]
for $l = 1,...,d,j = 2,...,{s_l}$. By the condition (C12), we can obtain ${{\dot p}_{{\alpha _1}}}\left( {\left| {\beta _{lj}^0} \right|} \right) = 0$ as $\alpha_1\rightarrow 0$, together with condition (C10), some simple calculations yields
\[n^{-1}{{\rm{{\cal R}}}_{\tau nh}}\left( {\bm \beta _{ - 1}^{(1)}} \right) - {{o_p}\left( {\bm\beta _{ - 1}^{(1)} - \bm\beta _{ - 1}^{0(1)}} \right)} = 0.\]
By Lemma \ref{lemma:2}, it is straightforward to prove that
\[{{\partial {{{\rm{{\cal R}}}_{\tau nh}}}\left( {\bm\beta _{ - 1}^{(1)}} \right)} \mathord{\left/
 {\vphantom {{\partial {R_{\tau nh}}\left( {\bm\beta _{ - 1}^{(1)}} \right)} {\partial \beta _{ - 1}^{(1)T}}}} \right.
 \kern-\nulldelimiterspace} {\partial \bm\beta _{ - 1}^{(1)T}}} =  - {\bm H_n}\left( {\bm \beta _{ - 1}^{0(1)}} \right)\left\{ {1 + {o_p}\left( 1 \right)} \right\} .\]
By Taylor expansion, Lemma \ref{lemma:2}, and the previous result,
\[\begin{array}{l}
 \sqrt n \left( {\bm\beta _{ - 1}^{(1)} -\bm \beta _{ - 1}^{0(1)}} \right) =  - {\left[ n^{-1}{{{\partial { {{{\rm{{\cal R}}}_{\tau nh}}}}\left( {\bm\beta _{ - 1}^{(1)}} \right)} \mathord{\left/
 {\vphantom {{\partial { {{{\rm{{\cal R}}}_{\tau nh}}}}\left( {\bm\beta _{ - 1}^{(1)}} \right)} {\partial \bm\beta _{ - 1}^{(1)T}}}} \right.
 \kern-\nulldelimiterspace} {\partial \bm\beta _{ - 1}^{(1)T}}}} \right]^{ - 1}}{n^{-1/2}{{{\rm{{\cal R}}}_{\tau nh}}}}\left( {\bm\beta _{ - 1}^{(1)}} \right)\left\{ {1 + {o_p}\left( 1 \right)} \right\} \\
~~~~~~~~~~~~~~~~~~~~~~~  ~= ({\mathbb{H} ^{ (1)}})^{-1} n^{-1/2}{\bm S_n}\left( {\bm\beta _{ - 1}^{0(1)}} \right) + {o_p}\left( 1 \right), \\
 \end{array}\]
 where ${\bm S_n}\left( {\bm \beta _{ - 1}^{0(1)}} \right)$ is $\sum\limits_{l = 1}^d {({s_l}-1)}  \times 1 $ sub-vector of ${\bm S _n}\left( {\bm \beta _{ - 1}^{0}} \right)$ corresponding to ${\bm\beta}_{-1}^{0(1)}$. Thus, by the central limit theorem and Slutsky's theorem, we can complete the proof of Theorem 3.2 (ii). $\square$\\

\noindent \textbf{Proof of  Theorem~4.1.} Let ${\bm \lambda ^0} = {\left\{ {{{\left( {\bm \lambda _l^0} \right)}^T}:1 \le l \le d} \right\}^T}$, ${k_n} = {J_n}{n^{{{ - 1} \mathord{\left/
 {\vphantom {{ - 1} 2}} \right.
 \kern-\nulldelimiterspace} 2}}} + J_n^{ - r + {1 \mathord{\left/
 {\vphantom {1 2}} \right.
 \kern-\nulldelimiterspace} 2}}$ and $\left\| \bm v \right\|_2 = C$, where $C$ is a sufficiently large constant. Our aim is to show that for
 any given $\epsilon$ there is a large constant $C$ such that, for a large $n$, we have
\begin{gather}
P\left\{ {\mathop {\inf }\limits_{{{\left\| \bm v \right\|}_2} = C} {{\rm{{\cal L}}}_{\tau n}^*}\left( {{\bm \lambda ^0} + { k_n}\bm v,\bm \beta ^0} \right) \ge {{\rm{{\cal L}}}_{\tau n}^*}\left( {{\bm \lambda ^0},\bm \beta ^0} \right)} \right\} \ge 1 -\epsilon.
\tag{A.29}\label{A.29}
\end{gather}
This implies with probability of at least $ 1 -\epsilon$ that there exists a local minimizer in the ball
$\left\{ {{\bm \lambda ^0} + {k_n}\bm v:{{\left\| \bm v \right\|}_2} \le C} \right\}$. This in turn implies that there exists a local minimizer $\bm{\bar \lambda}$ such that ${\left\| {\bm{\bar{  \lambda}}  - {\bm \lambda ^0}} \right\|_2} = {O_p}\left( {{k_n}} \right)$. We can find a constant
$c>0$ such that ${\left\| {\bm\lambda _l^0} \right\|_{{\bm D}}} = \sqrt {\bm\lambda _l^{0T}{\bm D}\bm\lambda _l^0} > c$ for $1\leq l\leq d_1$. Because $\left\| ({\bm\lambda _l^0 + {k_n}{\bm v_l}) -\bm \lambda _l^0} \right\|_2^2 = {o_p}\left( {{J_n}} \right)$ and that $\alpha_2=o(1)$, we have
\[P\left\{ {{p_{{\alpha _2}}}\left( {{{\left\| {\bm\lambda _l^0 + {k_n}{\bm v_l}} \right\|}_{{\bm D}}}} \right) = {p_{{\alpha _2}}}\left( {{{\left\| {\bm\lambda _l^0} \right\|}_{{\bm D}}}} \right)} \right\} \to 1~{\rm{if }}~ 1\leq l\leq d_1.\]
This fact implies that
\[n\sum\limits_{l = 1}^d {{p_{{\alpha _2}}}\left( {{{\left\| {\bm\lambda _l^0} \right\|}_{{\bm D}}}} \right)}  - n\sum\limits_{l = 1}^d {{p_{{\alpha _2}}}\left( {{{\left\| {\bm\lambda _l^0 + {k_n}{\bm v_l}} \right\|}_{{\bm D}}}} \right)}  \le 0\]
with probability tending to 1 since ${\left\| {\bm\lambda _l^0} \right\|}_{{\bm D}}={\bm {\lambda} _l^{0}{\bm D}\bm{\lambda }_l^0}=0$ for $ d_1+1\leq l\leq d $. By the proof of Lemma \ref{lemma:1} and the previous result,
\[\begin{array}{l}
 {{\rm{{\cal L}}}_{\tau n}^*}\left( {{\bm\lambda ^0},{\bm\beta ^0}} \right) - {{\rm{{\cal L}}}_{\tau n}^*}\left( {{\bm\lambda ^0} + {k_n}\bm v,{\bm\beta ^0}} \right) \\
  =  {{{\rm{{\cal L}}}_{\tau n}}\left( {{\bm\lambda ^0},{\bm\beta ^0}} \right) - {{\rm{{\cal L}}}_{\tau n}}\left( {{\bm\lambda ^0} + {k_n}\bm v,{\bm\beta ^0}} \right)} +  n{\sum\limits_{l = 1}^d {{p_{{\alpha _2}}}\left( {{{\left\| {\bm\lambda _l^0} \right\|}_{{\bm D}}}} \right)}  - n\sum\limits_{l = 1}^d {{p_{{\alpha _2}}}\left( {{{\left\| {\bm\lambda _l^0 + {k_n}{\bm v_l}} \right\|}_{{\bm D}}}} \right)} }  \\
 \leq   {{{\rm{{\cal L}}}_{\tau n}}\left( {{\bm\lambda ^0},{\bm\beta ^0}} \right) - {{\rm{{\cal L}}}_{\tau n}}\left( {{\bm\lambda ^0} + {k_n}\bm v,{\bm\beta ^0}} \right)}\\
  \leq 0\\
 \end{array}\]
 Hence, (\ref{A.29}) holds. This implies that ${\left\| {\bm{\bar{  \lambda}}  - {\bm \lambda ^0}} \right\|_2} = {O_p}\left( {{k_n}} \right)$. By the property (de Boor \cite{De01})
\[\begin{array}{l}
{C_1}{J_n}\left\| {\bm B{{\left( {\bm Z_i^T\bm\beta _l^0} \right)}^T}{\bm{\bar \lambda }_l} - \bm B{{\left( {\bm Z_i^T\bm\beta _l^0} \right)}^T}\bm{ \lambda} _l^0} \right\|_2^2 \le \left\| {{\bm{\bar \lambda }_l} - \bm{ \lambda }_l^0} \right\|_2^2 \\
~~~~~~~~~~~~~~~~~~~~~~~~~~~~~~~~~~~~~~~~~~~~~~~~~~~~\le {C_2}{J_n}\left\| {\bm B{{\left( {\bm Z_i^T\bm\beta _l^0} \right)}^T}{\bm{\bar \lambda }_l} - \bm B{{\left( {\bm Z_i^T\bm\beta _l^0} \right)}^T}\bm{ \lambda} _l^0} \right\|_2^2
 \end{array}\]
for some constants $C_1,C_2>0$. Thus, we have
\[\left| {{{\bar m}_l}\left( {{u_l},{\bm\beta ^0}} \right) - m_l^0\left( {{u_l}} \right)} \right| = {O_p}\left( {{{{k_n}} \mathord{\left/
 {\vphantom {{{k_n}} {\sqrt {{J_n}} }}} \right.
 \kern-\nulldelimiterspace} {\sqrt {{J_n}} }}} \right) = O_p\left( {\sqrt {{{{J_n}} \mathord{\left/
 {\vphantom {{{J_n}} n}} \right.
 \kern-\nulldelimiterspace} n}}  + J_n^{ - r}} \right)\]
uniformly for any $u_l\in [0,1]$ and $1\leq l\leq d$.  By the triangle inequality and (\ref{A.1}), we have $\left| {{{\bar m}_l}\left( {{u_l},{\bm\beta ^0}} \right) - m_l\left( {{u_l}} \right)} \right| = O_p\left( {\sqrt {{{{J_n}} \mathord{\left/
 {\vphantom {{{J_n}} n}} \right.
 \kern-\nulldelimiterspace} n}}  + J_n^{ - r}} \right)$ uniformly for any $u_l\in [0,1]$ and $1\leq l\leq d$. Theorem 4.1 follows from this result and Theorem 3.1. $\square$\\

\noindent \textbf{Proof of Theorem~4.2 .} Suppose for some $l=1+d_1,...,d$, ${\bm B}{\left( {\bm Z_i^T {{{\bar {\bm \beta} }_{\alpha _1l}}} } \right)^T}{{\bar {\bm \lambda } }_l}$ is not a linear function. Define $\bar{\bm \lambda} ^*$ to be the same as
$\bar{\bm \lambda} $
except that $\bar{\bm \lambda} _{l}$ is replaced by its projection onto the subspace
$\left\{ {{\bm \lambda _l}:{\bm B}{{\left( {\bm Z_i^T {{{\bar {\bm \beta} }_{\alpha _1l}}}} \right)}^T}{\bm \lambda _l}} ~{\rm{represents ~a ~linear~ function}} \right\}.$
 By definition of $\bar{\bm \lambda} $ and $\bar{\bm \lambda} ^*$, and ${\bm {\bar\lambda} _l^{*T}{\bm D}\bm{\bar \lambda }_l^*}=0$. we have
\begin{gather}
\begin{array}{l}
 0 \ge {{\rm{{\cal L}}}_{\tau n}^*}\left( {\bm{\bar \lambda} ,{{\bm{\bar\beta}_{\alpha_1} }}} \right) -
  {{\rm{{\cal L}}}_{\tau n}^*}\left( {{\bm{\bar \lambda} ^*},{\bm{\bar \beta }_{\alpha_1}}} \right) \\
 ~~ = \sum\limits_{i = 1}^n {{\rho _\tau }\left\{ {{Y_i} - \sum\limits_{l = 1}^d {{\bm B}{{\left( {\bm Z_i^T {{{\bm{\bar \beta} }_{\alpha _1l}}}} \right)}^T}{{\bm{\bar \lambda} }_l}{X_{il}}} } \right\}} - \sum\limits_{i = 1}^n {{\rho _\tau }\left\{ {{Y_i} - \sum\limits_{l = 1}^d {{\bm B}{{\left( {\bm Z_i^T {{{\bm{\bar \beta} }_{\alpha _1l}}}} \right)}^T}\bm{\bar \lambda} _l^*{X_{il}}} } \right\}}  \\
 ~~~~+ n {{p_{{\alpha _2}}}\left( {\sqrt {\bm{\bar \lambda}_l^T{\bm D}{{\bm{\bar \lambda} }_l}} } \right)}\\
   ~~=I-II+III.
 \end{array}
    \tag{A.30}\label{A.30}
\end{gather}
Let ${\Lambda _i} = \sum\limits_{l = 1}^d {\bm B{{\left( {\bm Z_i^T{\bm{\bar \beta }_{{\alpha _1}l}}} \right)}^T}\left( {\bm{\bar \lambda} _l^* - \bm\lambda _l^0} \right){X_{il}}}  + \sum\limits_{l = 1}^d {\left( {\bm B{{\left( {\bm Z_i^T{\bm{\bar \beta }_{{\alpha _1}l}}} \right)}^T}\bm\lambda _l^0 - {m_l}\left( {\bm Z_i^T\bm\beta _l^0} \right)} \right){X_{il}}} $.
Since ${\rho _\tau }\left( u \right) - {\rho _\tau }\left( v \right) \ge 2\left( {\tau  - I\left( {v \le 0} \right)} \right)\left( {u - v} \right)$ for any $u,v\in \mathbb{R}$, we have
\begin{gather}
\begin{array}{l}
I-II\ge  - \sum\limits_{i = 1}^n {\left( {2\tau  - 2I\left( {{\varepsilon _i} \le 0} \right)} \right)} \sum\limits_{l = 1}^d {\bm B{{\left( {\bm Z_i^T{\bm{\bar \beta }_{{\alpha _1}l}}} \right)}^T}\left( {{\bm{\bar \lambda }_l} - \bm{\bar \lambda} _l^*} \right){X_{il}}}  \\
  ~~~~~~~~~~~~- \sum\limits_{i = 1}^n {\left( {2I\left( {{\varepsilon _i} \le 0} \right) - 2I\left( {{\varepsilon _i} \le {\Lambda _i}} \right)} \right)} \sum\limits_{l = 1}^d {\bm B{{\left( {\bm Z_i^T{\bm{\bar \beta }_{{\alpha _1}l}}} \right)}^T}\left( {{\bm{\bar \lambda }_l} - \bm{\bar \lambda} _l^*} \right){X_{il}}}  \\
    ~~~~~~~~~\ge \left( { - {{\left\| {\sum\limits_{i = 1}^n {\left( {2\tau  - 2I\left( {{\varepsilon _i} \le 0} \right)} \right)} \sum\limits_{l = 1}^d {\bm B\left( {\bm Z_i^T{\bm{\bar \beta }_{{\alpha _1}l}}} \right){X_{il}}} } \right\|}_2}} \right. \\
   ~~~~~~~~~~~~\left. { - {{\left\| {\sum\limits_{i = 1}^n {\left( {2I\left( {{\varepsilon _i} \le 0} \right) - 2I\left( {{\varepsilon _i} \le {\Lambda _i}} \right)} \right)} \sum\limits_{l = 1}^d {\bm B\left( {\bm Z_i^T{\bm{\bar \beta }_{{\alpha _1}l}}} \right){X_{il}}} } \right\|}_2}} \right){\left\| {{\bm{\bar \lambda }_l} - \bm{\bar \lambda} _l^*} \right\|_2} \\
    ~~~~~~~~~= {O_p}\left\{-\left[ { \sqrt {{n \mathord{\left/
 {\vphantom {n {{J_n}}}} \right.
 \kern-\nulldelimiterspace} {{J_n}}}}+n{J_n^{{{ - 1} \mathord{\left/
 {\vphantom {{ - 1} 2}} \right.
 \kern-\nulldelimiterspace} 2}}} \left( {\sqrt {{{{J_n}} \mathord{\left/
 {\vphantom {{{J_n}} n}} \right.
 \kern-\nulldelimiterspace} n}}  + J_n^{ - r}} \right)} \right]\right\}{\left\| {\bm{\bar \lambda}  - {\bm{\bar \lambda }^*}} \right\|_2}\\
     ~~~~~~~~~={O_p}\left\{-\left[\sqrt n  + nJ_n^{ - r{{ - 1} \mathord{\left/
 {\vphantom {{ - 1} 2}} \right.
 \kern-\nulldelimiterspace} 2}}\right]\right\}{\left\| {\bm{\bar \lambda}  - {\bm{\bar \lambda }^*}} \right\|_2} \\
 \end{array}
 \tag{A.31}\label{A.31}
\end{gather}
since ${\max _{1 \le i \le n}}{\Lambda _i}=O_p\left( {\sqrt {{{{J_n}} \mathord{\left/
 {\vphantom {{{J_n}} n}} \right.
 \kern-\nulldelimiterspace} n}}  + J_n^{ - r}} \right)$. On the other hand,
\[\begin{array}{l}\sqrt {\bm{\bar \lambda} _l^T{\bm D}{{\bm{\bar \lambda} }_l}}  = \sqrt {{{\left( {{{\bm{\bar \lambda} }_l} - \bm \lambda _l^0} \right)}^T}{\bm D}\left( {{{\bm{\bar \lambda} }_l} - \bm \lambda _l^0} \right)}  = {O_p}\left( {{{\left( {{{{J_n}} \mathord{\left/
 {\vphantom {{{J_n}} n}} \right.
 \kern-\nulldelimiterspace} n}} \right)}^{{1 \mathord{\left/
 {\vphantom {1 2}} \right.
 \kern-\nulldelimiterspace} 2}}} + J_n^{ - r}} \right) = o\left( {{\alpha _2}} \right)
   \end{array}\]
   by the proof of Theorem~4.1 and assumption ${\left( {{{\left( {{{{J_n}} \mathord{\left/
 {\vphantom {{{J_n}} n}} \right.
 \kern-\nulldelimiterspace} n}} \right)}^{{1 \mathord{\left/
 {\vphantom {1 2}} \right.
 \kern-\nulldelimiterspace} 2}}} + J_n^{ - r}} \right)^{ - 1}}{\alpha _2} \to \infty $, we have
\begin{gather}
{p_{{\alpha _2}}}\left( {\sqrt {\bm{\bar { \lambda}} _l^T{\bm D}{\bm {\bar {\lambda} }_l}} } \right) = {\alpha _2}\sqrt {\bm{\bar { \lambda }_l}^T{\bm D}{\bm{\bar \lambda }_l}}
\tag{A.32}\label{A.32}
\end{gather}
 with probability tending to 1, by the definition of the SCAD penalty function. It remains to show that
\begin{gather}
{\left\| {\bar {\bm \lambda}  - {{\bar {\bm \lambda} }^*}} \right\|_2} = {\left\| {\bm{\bar  \lambda} _l -\bm{\bar  \lambda} _l^*} \right\|_2} =O_p\left( \sqrt {{J_n}\bar {\bm \lambda }_l^T{\bm D}{\bm{\bar \lambda }_l}} \right).
 \tag{A.33}\label{A.33}
\end{gather}
Then, combine (\ref{A.30}), (\ref{A.32}) and (\ref{A.33}), it follows that $III = n{\alpha _2}{J_n^{{{ - 1} \mathord{\left/
 {\vphantom {{ - 1} 2}} \right.
 \kern-\nulldelimiterspace} 2}}}{\left\| {\bm{\bar \lambda}  - {\bm{\bar \lambda }^*}} \right\|_2}$. Then, by the condition ${\left( {{{\left( {{{{J_n}} \mathord{\left/
 {\vphantom {{{J_n}} n}} \right.
 \kern-\nulldelimiterspace} n}} \right)}^{{1 \mathord{\left/
 {\vphantom {1 2}} \right.
 \kern-\nulldelimiterspace} 2}}} + J_n^{ - r}} \right)^{ - 1}}{\alpha _2} \to \infty $, it is easy to show that $III$ is of higher order than ${O_p}\left\{\left[\sqrt n  + nJ_n^{ - r{{ - 1} \mathord{\left/
 {\vphantom {{ - 1} 2}} \right.
 \kern-\nulldelimiterspace} 2}}\right]\right\}{\left\| {\bm{\bar \lambda}  - {\bm{\bar \lambda }^*}} \right\|_2}$. Thus, it follows that $I-II$ is dominated by $III$, which is large and positive for all sufficiently large $n$. Hence, we have a contradiction if ${\bm{\bar\lambda} _l^T{\bm D}\bm{\bar  \lambda}_l }>0$.

 Next, we prove (\ref{A.33}). Note that ${\bm {\bar\lambda} _l^{*T}{\bm D}\bm{\bar \lambda }_l^*}=0$. Furthermore,
since $\bm {\bar{\lambda}}_l ^*$ is the projection of $\bm {\bar{\lambda}}_l $ onto $\left\{ {{\bm \lambda _l}:  {\bm \lambda _l^T}{\bm D} {\bm \lambda _l}} =0 \right\}$, $\bm {\bar{\lambda}}_l -\bm {\bar{\lambda}}_l ^*$ is orthogonal to this projection space. Thus $\bm{\bar\lambda}  _l^T{\bm D}{\bm{\bar \lambda }_l} = {\left( {{\bm{\bar \lambda }_l} - \bm{\bar \lambda} _l^*} \right)^T}{\bm D}\left( {{\bm{\bar \lambda }_l} - \bm{\bar \lambda} _l^*} \right)$ and ${{{{\left( {{\bm{\bar \lambda }_l} - \bm{\bar \lambda} _l^*} \right)}^T}{\bm D}\left( {{\bm{\bar \lambda }_l} - \bm{\bar \lambda} _l^*} \right)} \mathord{\left/
 {\vphantom {{{{\left( {{\bm{\bar \lambda }_l} - \bm{\bar \lambda} _l^*} \right)}^T}{\bm D}\left( {{\bm{\bar \lambda }_l} - \bm{\bar \lambda} _l^*} \right)} {\left\| {{\bm{\bar \lambda }_l} - \bm{\bar \lambda} _l^*} \right\|_2^2}}} \right.
 \kern-\nulldelimiterspace} {\left\| {{\bm{\bar \lambda }_l} - \bm{\bar \lambda} _l^*} \right\|_2^2}}$ lies between the minimum and the maximum positive eigenvalues of $\bm D$, which is of order $1/{J_n}$. $\square$\\
